\documentclass[11pt]{article}

\usepackage[T1]{fontenc}
\usepackage[latin1]{inputenc}
\usepackage{mathtools,amsthm,amssymb}

\usepackage{mathrsfs}
\usepackage{xcolor}
\usepackage{authblk}

\usepackage{hyperref}
\usepackage[margin=2.5cm]{geometry}

\newtheorem{theo}{Theorem}[section]
\newtheorem{cor}[theo]{Corollary}
\newtheorem{lem}[theo]{Lemma}
\newtheorem{prop}[theo]{Proposition}

\numberwithin{equation}{section}

\def\tr{\mbox{\rm Tr}}
\def \PP{\mathbb{P}}
\def \RR{\mathbb{R}}
\def \EE{\mathbb{E}}

\def \LL{\mathbb{L}}
\def \BB{\mathbb{B}}

\newcommand{\MM}{\mathbb{M}}
\newcommand{\GG}{\mathbb{G}}
\newcommand{\Aa}{ {\cal A }}
\newcommand{\Ba}{ {\cal B }}
\newcommand{\Da}{ {\cal D }}
\newcommand{\Ea}{ {\cal E }}
\newcommand{\Sa}{ {\cal S }}
\newcommand{\Fa}{ {\cal F }}

\newcommand{\Xa}{ {\cal X }}
\newcommand{\Ma}{ {\cal M }}

\newcommand{\Pa}{ {\cal P }}
\newcommand{\Ta}{ {\cal T }}
\newcommand{\Ya}{ {\cal Y }}
\newcommand{\Wa}{ {\cal W }}
\newcommand{\Va}{ {\cal V }}
\newcommand{\Za}{ {\cal Z }}

\newcommand{\vertiii}[1]{{\left\vert\kern-0.25ex\left\vert\kern-0.25ex\left\vert #1
    \right\vert\kern-0.25ex\right\vert\kern-0.25ex\right\vert}}

\makeatletter
\DeclareRobustCommand\frownotimes{\mathbin{\mathpalette\frown@otimes\relax}}
\newcommand{\frown@otimes}[2]{
  \vbox{
    \ialign{##\cr
      \hidewidth$\m@th#1{}_\frown$\kern-\scriptspace\hidewidth\cr
      \noalign{\nointerlineskip\kern-.1pt}
      $\m@th#1\otimes$\cr
    }
  }
}
\makeatother

\begin{document}

\title{On the stability of matrix-valued Riccati diffusions}

\author[$1$]{Adrian N. Bishop}
\author[$2$]{Pierre Del Moral}
\affil[$1$]{{\small University of Technology Sydney (UTS); and CSIRO, Australia}}
\affil[$2$]{{\small INRIA, Bordeaux Research Center, France}}
\date{}

\maketitle

\begin{abstract}
The stability properties of matrix-valued Riccati diffusions are investigated. The matrix-valued Riccati diffusion processes considered in this work are of interest in their own right, as a rather prototypical model of a matrix-valued quadratic stochastic process. Under rather natural observability and controllability conditions, we derive time-uniform moment and fluctuation estimates and exponential contraction inequalities. Our approach combines spectral theory with nonlinear semigroup methods and stochastic matrix calculus. This analysis seem to be the first of its kind for this class of matrix-valued stochastic differential equation. This class of stochastic models arise in signal processing and data assimilation, and more particularly in ensemble Kalman-Bucy filtering theory. In this context, the Riccati diffusion represents the flow of the sample covariance matrices associated with McKean-Vlasov-type interacting Kalman-Bucy filters. The analysis developed here applies to filtering problems with unstable signals. 

\end{abstract}

{\footnotesize
\setcounter{tocdepth}{2}
\tableofcontents
}

\section{Introduction}

We introduce some matrix notation needed from the onset. Let $\Ma_{r}$ be the set of $(r\times r)$ real matrices with $r\geq 1$. Let $\Sa_r\subset \Ma_{r}$ be the subset of symmetric matrices, and $\Sa^0_r$, and $\Sa^+_r$ the subsets of positive semi-definite and definite matrices respectively. We write $A \geq B$ when $A-B\in \Sa^0_r$; and $A > B$ when $A-B\in \Sa^+_r$. We denote by $0$ and $I$ the null and identity matrices, for any $r\geq 1$. Given $R\in \partial \Sa_r^+:= \Sa_r^0-\Sa_r^+$ we denote by $R^{1/2}$ a (non-unique) symmetric square root of $R$. When $R\in\Sa_r^+$ we choose the unique symmetric square root. We write $A^{\prime}$ the transpose of $A$, and $A_{\mathrm{sym}}=(A+A^{\prime})/2$ its symmetric part. We denote by $\mathrm{Absc}(A):=\max{\left\{\mbox{\rm Re}(\lambda)\,:\,\lambda\in \mathrm{Spec}(A)\right\}} $ its spectral abscissa. We also denote by $\tr(A)$ the trace. When $A\in\Sa_r$ we let $\lambda_1(A)\geq \ldots\geq \lambda_r(A)$ denote the ordered eigenvalues of $A$. We equip $\Ma_{r}$ with the spectral norm $\Vert A \Vert=\Vert A \Vert_2=\sqrt{\lambda_{1}(AA^{\prime})}$ or the Frobenius norm $\Vert A \Vert=\Vert A \Vert_{\mathrm{Frob}}=\sqrt{\tr(AA^{\prime})}$. Let $\mu(A)=\lambda_{1}(A_{\mathrm{sym}})$ denote the (2-)logarithmic ``norm'' (which can be $<0$). We have $\mu(\cdot)\geq\mathrm{Absc}(\cdot)$.

\subsection{Description of the Model}

We associate with some given matrices $(A,R,S)\in (\Ma_{r}\times\Sa^0_r\times \Sa^0_r)$ the Riccati drift function $\Theta$ from $\Sa_r$ into itself defined by the matrix concave function
\begin{equation}\label{def-Riccati-drift}
	\Theta(P) :=(A-PS)P+P(A-PS)^{\prime}+\Sigma_{1,0}(P)\quad\mbox{\rm with}\quad \Sigma_{1,0}(P):=R+PSP
\end{equation}	
which may be written in canonical form, $\Theta(P)=AP + PA^{\prime} + R - PSP$. 

If the matrix pair $(A,R^{1/2})$ is stabilisable, and the pair $(A,S^{1/2})$ is detectable \cite{anderson-moore,anderson79,Lancaster1995}, then there exists a unique matrix: 
\begin{equation}\label{fixed-point-Ricc}
	\mathscr{P}_{\infty}\in\Sa^0_r\quad\mbox{\rm s.t.} \quad\Theta(\mathscr{P}_{\infty})=0\quad \mbox{\rm and $(A-\mathscr{P}_{\infty}S)$ is stable, i.e.}~{\mathrm{Absc}(A-\mathscr{P}_{\infty}S) <0}
\end{equation}
If $(A,R^{1/2})$ is controllable, then $\mathscr{P}_{\infty}\in\Sa^+_r$. See \cite{kucera72,Molinari77,Lancaster1995}.

The matrix-valued Riccati diffusions discussed in this article are defined by the stochastic model
\begin{equation}\label{f21}
	dQ_t~=~\Theta(Q_t)\,dt \,+\, \epsilon\,dM_t \end{equation}
with $t\in[0,\infty[$, $Q_0=Q\in\Sa^0_r$, and some noise parameter $\epsilon\geq 0$. The matrix-valued martingale is defined by
\begin{equation}
dM_t~:=~\left[Q_t^{1/2}~d\Wa_t~\Sigma^{1/2}_{\kappa,\varpi}\left(Q_t\right)\right]_{\mathrm{sym}} \label{def-matrix-martingale}
\end{equation}
where throughout $\Wa_t$ denotes an $(r\times r)$-matrix with independent Brownian entries. The non-negative map $\Sigma_{\kappa,\varpi}\,: \Sa_r\rightarrow \Sa_r^0$ is defined by
\begin{equation} \label{cond-UV-ref-intro}
	~~\Sigma_{\kappa,\varpi}(P)~:=~R+\kappa\,(P+\varpi I)\,S\,(P+\varpi I) \\
\end{equation}
for some finite $\varpi\geq 0$ and some (binary) parameter $\kappa\in \{0,1\}$.

For example, if $\kappa=0$, then $\Sigma_{0,\varpi}=\Sigma_{0,0}=R$ and thus $dM_t = [Q_t^{1/2}d\Wa_t\,R^{1/2}]_{\mathrm{sym}}$ or, explicitly
$$
	dQ_t \,=\, (AQ_t + Q_tA^{\prime} + R - Q_tSQ_t)\,dt \,+\, \frac{\epsilon}{2}\left[Q_t^{1/2}d\Wa_t\,R^{1/2} + R^{1/2}\,d\Wa_t^{\prime}\,Q_t^{1/2}\right]
$$
This special case ($\kappa=0$) defines, in some sense, a minimal prototype of a forward-in-time matrix-valued Riccati diffusion in the space of symmetric positive (semi-)definite matrices.

We let $\phi^{\epsilon}_t(Q):=Q_t$ be the stochastic flow of the matrix diffusion equation (\ref{f21}) with $Q_0=Q$. Whenever it exists, the inverse stochastic flow of (\ref{f21}) is denoted by $\phi^{-\epsilon}_t(Q):=Q_t^{-1}$. For any $0\leq s\leq t$, we let $\mathcal{E}_{s,t}^{\epsilon}(Q)$ be the transition semigroup associated with the flow of random matrices $\left[A-\phi^{\epsilon}_t(Q)S\right]$, i.e. the solution of the forward and backward equations
\begin{equation}\label{def-Est}
 \partial_t \mathcal{E}_{s,t}^{\epsilon}(Q)\,=\,\left[A-\phi^{\epsilon}_t(Q)S\right]\,\mathcal{E}^{\epsilon}_{s,t}(Q)\qquad\mbox{\rm and}\qquad
\partial_s \mathcal{E}_{s,t}^{\epsilon}(Q)\,=\,-\mathcal{E}_{s,t}^{\epsilon}(Q)\,\left[A-\phi^{\epsilon}_s(Q)S\right]
\end{equation}
with $\mathcal{E}^{\epsilon}_{t,t}(Q)=I$. When $s=0$ we write $\mathcal{E}^{\epsilon}_{t}(Q)$ instead of $\mathcal{E}^{\epsilon}_{0,t}(Q)$. We write $\phi_{t}(Q)$ and  $\mathcal{E}_{s,t}(Q)$ instead of $\phi_{t}^0(Q)$, and $\mathcal{E}^0_{s,t}(Q)$, to denote the flow of the deterministic matrix Riccati differential equation when $\epsilon=0$, and the exponential semigroup defined via $\phi_{t}(Q)$.

\subsection{Background and Motivation}

The main concern in this article is the matrix-valued Riccati diffusion in (\ref{f21}) and its, time-uniform, moment boundedness and fluctuation behaviour, along with its stability and contraction properties. 

Positive semi-definite matrix diffusions in a specialized form of (\ref{f21}) also arise in multivariate statistics,  econometrics and financial mathematics. For example, the Wishart process considered in~\cite{bru1991wishart,Cox/Ingersoll/Ross:1985a} corresponds to the choice of parameters $\Sigma_{0,0} = R$, $S=0$ and {\em $A$ stable}. In financial mathematics, Wishart-type processes are used to model multivariate stochastic volatility in equity and fixed income models. The article~\cite{cuc} also considers a general class of affine processes in the cone of positive semi-definite matrices. These processes combine Wishart diffusions and pure jump processes with a compensator of affine-type. The Feller properties of the transition semigroup of affine processes are developed in the articles~\cite{cuc,keller}. The main feature of these processes is that the characteristic functions and the moment generating functions are explicitly known. For more details on the mathematical analysis of affine processes we refer to~\cite{cuc,keller,mayerhofer2011strong}, and the references therein. To the best of our knowledge, whenever they exist, such explicit formulae are unknown for general matrix-valued Riccati diffusions of the form (\ref{f21}), as soon as $S\not=0$.

In the context of filtering of conditionally Gaussian signal-observation models \cite{Liptser2001,Haussmann1988}, the associated conditional covariance matrix of the posterior (filtering) distribution is random and satisfies a type of Riccati diffusion equation, see \cite{Haussmann1988}.

We remark that different models involving backward matrix Riccati diffusions arise in linear-quadratic optimal control problems with random coefficients; see e.g. \cite{bismut1976linear,hu2003indefinite,Kohlmann2003}. Another class of random Riccati equations, different from the diffusion equation (\ref{f21}), arises in network control and filtering with random observation losses; see for instance \cite{sinopoli2004kalman}. The details of these works are beyond the scope of the forward-in-time Riccati diffusions considered herein.

\subsubsection{Ensemble Kalman-Bucy-Type Filters}

The stochastic Riccati equations defined by \eqref{f21} may be motivated by applications in signal processing and data assimilation problems, and more particularly in the stochastic analysis of ensemble Kalman-Bucy-type filters (abbreviated {\tt EnKF}) \cite{evensen03,sakov2008a}.  

 In this context, up to a change of probability space, the matrix-valued Riccati diffusion \eqref{f21} describes the evolution of the sample covariance associated with these filters. With this application, the general form of (\ref{cond-UV-ref-intro}) accommodates two popular {\tt EnKF} filter models (determined by the binary switch $\kappa\in\{0,1\}$) \cite{evensen03,sakov2008a}; as well as accommodating a class of inflation-based regularization methods (determined by $\varpi\geq0$)~\cite{Bishop/DelMoral/Pathiraja:2017}. The case $\varpi=0$ in (\ref{cond-UV-ref-intro}) corresponds to non-regularized models. We refer to Section~\ref{EnKF-sec} for further discussion on these particle-type filters and on inflation-regularisation.
 
In the context of state estimation, the difference between the {\tt EnKF} sample mean
and the true signal state (i.e. the estimation error) is described by a stochastic Ornstein-Uhlenbeck-type vector-valued process of the form,
\begin{equation}\label{Stoch-OU}
dX_t\,=\,\left(A-Q_tS\right)X_t\,dt+(\Sigma^{\overline{\epsilon} }_{\kappa,\varpi}(Q_t))^{1/2}\,dW_t
\end{equation}
where $W_t$ is an $r$-dimensional Wiener process independent of $\Wa_t$, and $(\Sigma^{\overline{\epsilon}}_{\kappa,\varpi})^{1/2}$ denotes the square root of the non-negative map $\Sigma^{\overline{\epsilon}}_{\kappa,\varpi}$ from $\Sa_r$ into $\Sa_r^0$ defined by,
\begin{equation}\label{ref-inf-reg-EnKF}
	\Sigma^{\overline{\epsilon}}_{\kappa,\varpi} \,:=\, \Sigma_{1,\varpi} + \overline{\epsilon}^2\, \Sigma_{\kappa,\varpi}
\end{equation}
where $\overline{\epsilon}=\epsilon/\sqrt{4+\epsilon^2} \,<\, \epsilon$.

Whenever ${\epsilon}=0=\varpi$, the diffusion process (\ref{Stoch-OU}) resumes to the difference (i.e. error) between the classical Kalman-Bucy filter \cite{Bishop/DelMoral:2016} and the true signal state of an auxiliary linear-Gaussian process with drift matrix $A$ and diffusion matrix $R$. In this case we have
\begin{equation}\label{Stoch-OU-R}
dX_t=(A-\mathscr{P}_tS)\,X_t\,dt+\Sigma^{1/2}_{1,0}(\mathscr{P}_t)\,dW_t
~\quad\mbox{\rm with the Riccati equation}\quad~\partial_t\mathscr{P}_t=\Theta(\mathscr{P}_t)
\end{equation}
Note that $\mathscr{P}_t=\EE(X_tX_t^{\prime})$ coincides with the covariance matrix of the state estimation error defined by the Ornstein-Uhlenbeck process in the l.h.s. of (\ref{Stoch-OU-R}) and $\phi_{t}(Q)=\mathscr{P}_t$ when $\mathscr{P}_0=Q$. 

Under appropriate controllability and observability conditions, one of the main features of the Kalman-Bucy filter is that it delivers a {\em stable} state estimate of the underlying signal with {\em unstable drift matrix $A$} and uniformly w.r.t. the time parameter. In particular, when the pair of matrices $(A,R^{1/2})$ is stabilisable, and  $(A,S^{1/2})$ is detectable, then the error $X_t$ is a stable process in the sense that $(A-\mathscr{P}_tS)$ delivers a uniformly exponentially stable linear (time-varying) system \cite{Bishop/DelMoral:2016}. Moreover, $\phi_{t}(Q)\rightarrow_{t\rightarrow\infty}\mathscr{P}_{\infty}$ exponentially fast for any $Q\in\Sa_r^0$, where $\mathscr{P}_{\infty}\in\Sa^0_r$ is the stabilising fixed point defined in (\ref{fixed-point-Ricc}). See \cite{kucera72,Molinari77,Lancaster1995}, and the convergence results in \cite{Kwakernaak72,callier81}. In this situation, the second moments of the  diffusion process (\ref{Stoch-OU-R}), as well as the solution of the deterministic Riccati equation are uniformly bounded w.r.t. the time horizon. For further discussion on these stability properties we refer to~\cite{Bishop/DelMoral:2016,bd-CARE} and the references therein. We point to \cite{anderson-moore,anderson79,Lancaster1995} for precise definitions of controllable, stabilizable, and their duals, observable, detectable. 

In some cases (later), we may ask for a stronger stability property; i.e. $\mu(A-\mathscr{P}_{\infty}S) <0$. We claim this condition requires an even stronger notion of observability and controllability. A crude, yet sufficient, example here is to suppose a change of basis such that $A$ is symmetric and $\Sa\propto I$.

The stability of $(A-\mathscr{P}_tS)$, and of $\phi_{t}(Q)\rightarrow_{t\rightarrow\infty}\mathscr{P}_{\infty}$, is directly related \cite{Bishop/DelMoral:2016,bd-CARE} to the contractive properties of the deterministic semigroup $\mathcal{E}_{s,t}(Q)$. The stochastic equation (\ref{Stoch-OU}) implies that the stability properties of {\tt EnKF} state estimators with $\varpi=0$ depend on the stability properties of the stochastic exponential semigroup $\mathcal{E}_{s,t}^{\epsilon}(Q)$. 

Section~\ref{EnKF-sec} examines the stability of these ensemble filtering methods in detail (with any $\varpi\geq 0$, $\kappa\in \{0,1\}$), and in relation to our main results on Riccati diffusions of the form (\ref{f21}). We also relate the contractive properties of $\mathcal{E}_{s,t}^{\epsilon}(Q)$ to the stability of the state estimation error flow $X_t$ in (\ref{Stoch-OU}).

\subsection{General Statements of the Main Results}

We make the standing assumption throughout that $(A,R^{1/2})$ is stabilizable and $(A,S^{1/2})$ is detectable; see \cite{Kwakernaak72,anderson79,Lancaster1995}. The main concern in this article is the general matrix-valued Riccati diffusion in (\ref{f21}). This includes as a special case the minimal prototype for this type of matrix-valued quadratic stochastic differential equation that arises when $\kappa=0$. In particular, we are interested in the stability of the flow $\phi^{\epsilon}_t(Q)$ and the contraction properties of the associated exponential semigroup $\mathcal{E}_{s,t}^{\epsilon}(Q)$. We also consider moment estimates on $\phi^{\epsilon}_t(Q)$ and the fluctuation properties of $\phi^{\epsilon}_t(Q)$ about the deterministic Riccati flow $\phi_t(Q)$. Later we also consider applications of this work to {\tt EnKF} theory, and we consider the stability of the associated Ornstein-Uhlenbeck-type flow $X_t$.

The analysis of Riccati diffusions of the form (\ref{f21}) with the parameters $(\kappa,\varpi)=(1,0)$ has been started in \cite{Bishop/DelMoral/Niclas:2017,DelMoral/Tugaut:2016}. In these articles, the authors provide several uniform convergence results when $S$ is proportional to the identity and when {\em $A$ is a stable matrix}. In \cite{Bishop/DelMoral/Niclas:2017} we also provide a complete Taylor-type stochastic expansion of the Riccati flow with estimates given at any order with bounded remainder terms, and with a fluctuation analysis considered over the entire path space of the matrix-valued stochastic Riccati flow. We remark that in the scalar case we also have a complete time-uniform fluctuation and stability analysis of one-dimensional Riccati diffusions in \cite{2017arXiv171110065B}. Nevertheless, the understanding of the long time behaviour of \emph{matrix-valued} Riccati diffusions with arbitrary matrices $A$ seems lacking, and requires the development of new mathematical techniques. 

To address this problem, we develop a novel stability and fluctuation analysis for (\ref{f21}) combining spectral theory with nonlinear semigroup techniques in matrix spaces. The present article thus complements the recent article~\cite{2017arXiv171110065B} dedicated to the stability and fluctuation properties of general one-dimensional Riccati diffusions. And this article extends \cite{Bishop/DelMoral/Niclas:2017,DelMoral/Tugaut:2016} in a number of directions, and corrects some claims in \cite{DelMoral/Tugaut:2016}. 

The main contributions of this work are listed succinctly below and discussed throughout the remainder. For more precise statements, we refer to the series of theorems stated in Section~\ref{statement-sec}.

$\bullet$~~As shown in~\cite{2017arXiv171110065B}, for one-dimensional models the equation (\ref{f21}) has a unique strong solution in $\Sa^0_1=[0,\infty[$. In addition, the origin is repellent as soon as $R>0$, for any $\epsilon\geq 0$. To the best of our knowledge, the extension of this result in the multivariate case is unknown.

In the present article, we show that the existence and uniqueness of a weak solution in $\Sa^0_r$ of equation (\ref{f21}) is ensured for any time horizon and any fluctuation parameter $\epsilon\geq 0$; see Theorem~\ref{theo-existence-s-ric}.

Up to a change of probability space, the sample covariance matrices of ensemble Kalman-Bucy filters with $N+1$ particles satisfies  (\ref{f21}) with $\epsilon\propto 1/\sqrt{N}$. The rank of a sample covariance matrix $Q_t$ is thus at most $N$, and with $N<r$, it follows that $Q_t\in \partial \Sa_r^+$ is a unique weak solution of equation (\ref{f21}). We refer to Section~\ref{EnKF-sec} for a more precise discussion of ensemble Kalman-Bucy filters.

$\bullet$~~Whenever $S\in \partial \Sa_r^+$, without additional regularity properties, the solution of (\ref{f21}) may blow up when the matrix $A$ is unstable as the time horizon $t\rightarrow\infty$. Nevertheless, when the pair of matrices $(A,R^{1/2})$ is stabilisable, and $(A,S^{1/2})$ is detectable \cite{Lancaster1995}, for any $t\geq 0$ and any $\epsilon\geq 0$, we prove the following uniform under bias estimate
\begin{equation}\label{control-intro-bias}
\EE\left[\phi^{\epsilon}_t(Q)\right]\,\leq~ \phi_{t}\left(Q\right)\,\leq~ c\,(1+\Vert Q\Vert)\,I
\end{equation}
for a finite constant $c<\infty$ that doesn't depend on the time horizon. In addition, the above estimate does not depend on $\|Q\|$ when $t\geq \upsilon$ for any $\upsilon>0$ and with some parameter $c$ dependent on $\upsilon$. The l.h.s. under bias estimate is a consequence of the inequality (\ref{ref-phi-1-est-proof}); see also Theorem 1.3 in~\cite{Bishop/DelMoral/Niclas:2017}. The proof of the r.h.s. uniform estimate is in~\cite{Bishop/DelMoral:2016,bd-CARE}; e.g. it is easy to verify $\|\phi_{t}(Q)\| \leq c\left(\|\mathscr{P}_\infty\|\vee\|Q\|\right)$. See also the refined uniform bias estimates in Theorem \ref{theo-3-intro}.

The uniform moment estimates (\ref{control-intro-bias}) ensure that the stochastic Riccati diffusion (\ref{f21}) is uniformly tight. By Prohorov's theorem, this implies that the distributions of the random states $(Q_t)_{t\geq 0}$ is relatively compact so there exists 
at least one limiting invariant distribution $\Gamma^{\epsilon}_{\infty}$ on $\Sa_r^0$ and a sequence of random times $t_n$ such that 
$$
\Pi^{\epsilon}_{t_n}(P,dQ)~~\stackrel{ weakly}{-\!\!\!-\!\!\!-\!\!\!-\!\!\!-\!\!\!-\!\!\!\longrightarrow}_{n\rightarrow\infty}~~\Gamma^{\epsilon}_{\infty}(dQ)
$$
where $\Pi^{\epsilon}_t$ denotes a Markov semigroup for $Q_t$, defined more formally later. We remark however, that at this level of generality, it is difficult to ensure the uniqueness of the invariant measure and the stability properties of matrix Riccati diffusions.

$\bullet$~~One central question towards this goal is to analyze the regularity properties of the transition semigroup associated with the Riccati diffusion (\ref{f21}). 
Firstly, observe that the positive map $\Sigma_{\kappa,\varpi}$ defined in (\ref{cond-UV-ref-intro}) satisfies the inequality
\begin{equation}\label{cond-UV-ref-intro-uv}
\Sigma_{\kappa,\varpi}(P)\,\leq\, U+PVP~\quad \mbox{\rm with}~\quad  (U,V)\,:=\,\big(R+\kappa\,\varpi\, S(S+\varpi I), \,\kappa(S+\varpi I)\big)~~
\end{equation}
Notice that when $\varpi=0$, the estimate (\ref{cond-UV-ref-intro-uv}) resumes to the formula 
\begin{equation}\label{cond-UV-ref-intro-uv-ref}
\Sigma_{\kappa,0}(P) \,=\, U+PVP ~\qquad\Longrightarrow\qquad (U,V)\,=\,(R,\kappa S) \,\in\, \{(R,S),(R,0)\}
\end{equation}
Note implicitly that $(U,V)$ are $(\kappa, \varpi)$-indexed. The introduction of these matrices allows us to control the positive diffusion map $\Sigma_{\kappa,\varpi}(P)$ in terms of a single quadratic-type form on $P$. From (\ref{cond-UV-ref-intro-uv-ref}), this control trivially holds when $\varpi=0$.

Now we introduce a fluctuation parameter of the form,
\begin{equation}\label{Hyp-RS-repsilon}
 \varepsilon_0 \,:=\,\sup{\left\{\epsilon\geq 0~:~R^{\epsilon}:=R-\frac{\epsilon^{2}}{4}\,(r+1)\,U\geq 0\quad\mbox{\rm and}\quad S^{\epsilon}:=S-
\frac{\epsilon^{2}}{4}\,(r+1)\,V\geq 0\right\}}
\end{equation}
with the matrices $(U,V)$ defined in (\ref{cond-UV-ref-intro-uv}). Notice that this condition may simplify significantly,
$$
\varpi=0 \quad\Longrightarrow\quad (U,V)\in \{(R,S), (R,0)\}  \quad\Longrightarrow\quad \varepsilon_0:=2/\sqrt{r+1}
$$
With $\epsilon\leq \varepsilon_0$, we prove that the matrix Riccati diffusion (\ref{f21}) has a unique strong solution in $\Sa^+_r$ and that it never hits the boundary $\partial \Sa^+_r$ on any positive time horizon.

 In addition, the transition semigroup of $Q_t$ is strongly Feller, and admits a smooth positive density w.r.t. the Lebesgue measure on $\Sa^+_r$; thus, it is irreducible  (cf. Theorem~\ref{theo-existence-s-ric}).  

The uniqueness of the invariant measure $\Gamma^{\epsilon}_{\infty}$ now follows via the fact the semigroup transitions are mutually absolutely continuous. This also ensures that $\Gamma^{\epsilon}_{\infty}$ has a positive density on $\Sa^+_r$.

$\bullet$~~To quantify the convergence to equilibrium we need to quantify in more detail the moments of the Riccati flow and its inverse. 
We need the additional fluctuation parameters,
\begin{equation} \label{Hyp-UV-eps}
\begin{split}
\varepsilon_n(V)~:=&~\sup{\left\{\epsilon\geq 0~:~\frac{\epsilon^2}{2}\,r\,(n-1)\,\lambda_1(V) \,<\,\lambda_r(S)\right\}} \\
\varepsilon_n(U,V)~:=&~\sup{\left\{\epsilon\in [0,\varepsilon_0]~:~\frac{\epsilon^2}{2}\,\left[
 (1+nr)\,\lambda_1(U)+\,\frac{\lambda_1(V)}{4}\,r\right]\,<\,\lambda_r\left(R\right)\right\}}
\end{split} 
\end{equation}
Observe that $S\in \Sa^+_r\Longrightarrow\varepsilon_{1}(V)=\infty$. If $\kappa=0$ and $S\in \Sa^+_r$ then $\varepsilon_n(V)=\infty$ for all $n\geq1$. Actually, we have $\varepsilon_n(V)>0$ if and only if $S\in \Sa^+_r$ and $\varepsilon_n(U,V)>0$ if and only if $R\in \Sa^+_r$.

When $\epsilon\leq\varepsilon_n(V)\wedge\varepsilon_n(U,V)$ we prove that the $n$-th moments of $Q_t$ and its inverse matrix $Q^{-1}_t$ are uniformly bounded  w.r.t. the time horizon. In addition these moments are uniformly bounded w.r.t. the initial state for strictly positive time horizons; see Theorem~\ref{theo-existence-s-ric-proof-bis}. 

$\bullet$~~When $\epsilon\leq\varepsilon_1(V)\wedge\varepsilon_1(U,V)$ we also show that the function $\Lambda(P):=\Vert P\Vert_2+\Vert P^{-1}\Vert_2$ is a Lyapunov function on $\Sa^+_r$ with compact level sets. In this situation, the distribution of $Q_t$ converges exponentially fast to the unique invariant probability measure $\Gamma^{\epsilon}_{\infty}$; see Theorem~\ref{theo-stab-intro}.

$\bullet$~~This article is also concerned with uniform fluctuation estimates of $\phi^{\epsilon}_t(Q)$ about the limiting object $\phi_t(Q)$, and w.r.t. the time horizon.  For instance, when $S\in \Sa^+_r$ and $\kappa=0$ in (\ref{cond-UV-ref-intro}), then for {\em any fluctuation parameter $\epsilon\geq 0$} and any $n\geq 1$ we have the uniform estimates
$$
\sup_{t\geq 0}\,{\EE\left[\Vert \phi^{\epsilon}_t(Q)-\phi_{t}(Q)\Vert^n \right]^{1/n}} \,\leq\, c_{n}(Q)\,\epsilon
$$
for some constant $c_{n}(Q)$ whose values only depends on $n$ and $Q$; see Theorem~\ref{theo-3-intro}. 

$\bullet$~~These latter uniform fluctuation estimates allow one to quantity with some precision the exponential decay of the exponential semigroups $\mathcal{E}_{s,t}^{\epsilon}(Q)$ as discussed further in Section~\ref{sec-stability}.

To get some intuition on the complexity of matrix Riccati diffusion models, we mention that  the evolution model of the eigenvalues of $Q_t$ is generally not closed, in the sense that it also depends on the random eigenvectors of the Riccati diffusion. We note however, that for diagonal matrices $(A,R,S)$, the solution $\mathscr{P}_t$ of the deterministic matrix Riccati differential equation in (\ref{Stoch-OU-R}) is diagonal as soon as $\mathscr{P}_0$ is diagonal. In addition, when $(A,R,S)$ are proportional to the identity, we have $\mathscr{P}_t$ proportional to the identity with a time-varying proportionality constant that solves a naturally associated univariate (scalar) Riccati differential equation. Even in this simplified (identity proportional) setting, these elementary properties fail for the matrix diffusion (\ref{f21}) as soon as $\epsilon>0$. To be more precise, when $r>1$ the evolution of the eigenvalues of $Q_t$ is associated with {\em an additional repulsion force} that prevents the collision of eigenvalues. These logarithmic Coulomb repulsion forces are dictated by the second order Hadamard variational formula, the strength of repulsion is inversely proportional to their separation. The interacting diffusion model discussed above is closely related to the Dyson-Brownian motion model that represents the evolution of the eigenvalues of Gaussian orthogonal ensembles $(\Wa_t+\Wa^{\prime}_t)$. For example, in the simple setting $A=R=S=I$ and $\varpi=\kappa=0$, and $\epsilon\leq \varepsilon_0$, the ordered eigenvalues  $0< \lambda_r(t)< \ldots<\lambda_1(t)$ of the Riccati diffusion matrix $Q_t$ satisfy the Dyson-type diffusion equation 
\begin{equation}\label{dyson-mod-intro-ref}
\begin{array}{l}
\displaystyle
d\lambda_i(t)
=
\left[2\lambda_i(t)+1-\lambda_i(t)^2+\frac{\epsilon^2}{4}~\sum_{j\not=i}~\frac{\lambda_i(t)+\lambda_j(t)}{\lambda_i(t)-\lambda_j(t)}\right]dt+
\epsilon~\sqrt{\lambda_i(t)}~dW^{i}_t
\end{array}
\end{equation}
for some sequence $W^{i}_t$ of independent Brownian motions. We refer to Section~\ref{dyson-section} for a more detailed and general discussion on these Dyson-type equations. For background details on Dyson-Brownian motions we refer to~\cite{dyson1962brownian,mehta2004random,anderson2010introduction,tao2012topics}. Most of these studies are primarily concerned with the behaviour of eigenvalues for isotropic-type Gaussian models, when $r\rightarrow\infty$. The literature on positive semidefinite matrix diffusions is also mainly concerned with the existence and numerical approximation schemes on finite time horizon. In contrast with these works, the present article is concerned with the fluctuation and the stability analysis of these models over long time horizons $t\rightarrow\infty$.

\subsubsection{Article Organisation}

The main contributions of the article are presented in Section~\ref{statement-sec}. In Section~\ref{EnKF-sec} we illustrate the impact of our results in the context of ensemble Kalman-Bucy filters, including inflation-regularization methodologies. Section~\ref{sec-Riccati-flows} presents some pivotal results concerning Riccati flows, including a characterisation of inverse matrix-valued Riccati diffusions, a matrix-comparison lemma and Liouville determinant-type formulae for Riccati diffusion flows. The end of the section is concerned with the derivation of the Dyson-type equations associated with the evolution of the eigenvalues of this class of matrix diffusions. Section~\ref{proof-main-theo} is dedicated to the proof of the main results stated in Section~\ref{statement-sec}.

\subsection{Some Basic Notation}

This section presents some basic notation and preliminary results necessary for the statement of our main results. 

Throughout, we write $c,c_{n},c_{\upsilon},c_{\upsilon,n}, c_{\upsilon,n}(x),c_{\upsilon,n}(Q),c_{\upsilon,n}(Q,x)\ldots$ for some positive universal constants whose values may vary from line to line, but which only depend on some parameters $n,\upsilon,x,Q$, etc, as well as on the parameters of the Riccati process $(A,R,S,U,V)$. Importantly, these constants do not depend on the time horizon $t$, nor on the fluctuation parameters $(\epsilon,\overline{\epsilon})$.

Given a suitably regular matrix-valued stochastic process $t\mapsto A_t\in \Ma_r$, for any $t\geq 0$ and $n\geq 1$ we set
$$
   \vertiii{A_t}_n=\EE\left[\Vert A_t\Vert^n\right]^{1/n}
$$

We denote by $\Pi^{\epsilon}_t$ the Markov semigroup of $Q_t$ defined for any bounded measurable function $F\in\Ba(\Sa_r)$ 
and $Q\in \Sa_r^0$ by
$$
\Pi_t^{\epsilon}(F)(Q):=\EE\left[F(\phi^{\epsilon}_t(Q))\right]
$$

We consider the symmetric tensor products on $\Sa_r^0$ defined by
\begin{eqnarray*}
P_1\otimes_sP_2&:=&\frac{1}{2} \left(P_1\,\otimes\, P_2+P_2\,\otimes \,P_1\right)\\
P_1\,\overline{\otimes}_s\,P_2&:=&\frac{1}{2} \left(P_1\,\overline{\otimes}\, P_2+P_2\,\overline{\otimes} \,P_1\right) \\
P_1\,\frownotimes \,P_2&:=&\frac{1}{2}~(P_1\,\otimes_s\,P_2+P_1\,\overline{\otimes}_s\,P_2)~\geq~ 0
\end{eqnarray*}
with the tensor products
$$
(P_1\,\overline{\otimes} \,P_2)((i,j),(k,l)) \,:=\,(P_1\otimes P_2)((i,j),(l,k))\,=\,P_1(i,l)P_2(j,k)
$$
In this notation, we have
$
(I\,\overline{\otimes}\,I)(H)=H^{\prime}
$.
In addition, the angle bracket of the matrix-valued martingale $M_t=(M_t(i,j))_{1\leq i,j\leq r}$ defined in (\ref{f21}) is given by the formula
\begin{equation}\label{ref-M-frownotimes}
\partial_t\langle M(i,j)~\vert~M(k,l) \rangle_t=\left(Q_t\frownotimes \Sigma_{\kappa,\varpi}(Q_t)\right)((i,j),(k,l))
\end{equation}

We set $\overline{r}:=r(r+1)/2$ and we equip the product space $\RR^{\overline{r}}$  with the inner product
$$
\langle x,y\rangle_{\overline{r}} \,=\, \sum_{1\leq i\leq r}x_{i,i}y_{i,i}+2\sum_{1\leq i<j\leq r}x_{i,j}y_{i,j}
$$
where $x,y\in\RR^{\overline{r}}$ and where we index these vectors via $x=(x_{i,j})_{1\leq  i\leq j\leq r}$ and $y=(y_{i,j})_{1\leq  i\leq j\leq r}$. We equip $\RR^{\overline{r}}$ with the rescaled Lebesgue measure,
$$ 
	\gamma_{\overline{r}}(dx) \,:=\, 2^{-r(r-1)/4}~\prod_{1\leq i\leq j\leq r} dx_{i,j}
$$
Let $E_{i,j}$ with $1\leq i,j\leq r$ be the $(r\times r)$-matrices with entries $E_{i,j}(k,l)=1_{(i,j)=(k,l)}$. For any $H\in \Sa_r$
we have,
$$
H\,=\sum_{1\leq i\leq r}~H_{i,i}~E_{i,i}^s+\sum_{1\leq i<j\leq r}~\sqrt{2}H_{i,j}~E^{s}_{i,j}\,=\sum_{1\leq i\leq j\leq r}\langle H,E^{s}_{i,j}\rangle_{\mathrm{Frob}}~E_{i,j}^s
$$
with the orthonormal basis of $\Sa_r$ given by 
$$
 E^{s}_{i,i}\,=\,E_{i,i}\qquad \mbox{\rm and}\qquad E^{s}_{i,j}\,:=\,\frac{E_{i,j}+E_{j,i}}{\sqrt{2}}~1_{i<j}
$$
The above decomposition yields the  isomorphism $\varsigma:\left(\Sa_r,\langle\cdot,\cdot\rangle_{\mathrm{Frob}}\right)\mapsto\left( \RR^{\overline{r}},\langle\cdot,\cdot\rangle_{\overline{r}}\right)$ defined by
$$
\begin{array}{l}
\left(\varsigma (H)\right)^{\prime}=\left(\left(H_{1,1},\sqrt{2}H_{1,2},\ldots,\sqrt{2}H_{1,r}\right),\ldots\left(H_{i,i},\sqrt{2}H_{i,i+1},\ldots,\sqrt{2}H_{i,r}\right),
\ldots,H_{r,r}\right)\qquad\qquad\\
\\
\qquad\qquad\qquad\Longrightarrow\qquad \langle H_1,H_2\rangle_{\mathrm{Frob}}\,:=\,\tr(H_1H_2)\,=\,\left(\varsigma (H_1)\right)^{\prime}\left(\varsigma (H_2)\right)\,:=\,\langle
\varsigma (H_1),\varsigma (H_2)\rangle_{\overline{r}}
\end{array}$$
Note the set $\Da_{\overline{r}}:=\varsigma(\Sa_r^+)$ is an open smooth manifold embedded in $ \RR^{\overline{r}}$ with boundary $\partial \Da_{\overline{r}}=\varsigma(\partial \Sa_r^+)$ of {\em $\gamma_{\overline{r}}$-null measure}  on $\RR^{\overline{r}}$ and parametrised by the equation $\det(\varsigma^{-1}(\cdot))=0$.

We define the Lebesgue measure on $\Sa_r$ using $\varsigma^{-1}$ and $\gamma_{\overline{r}}$ according to the natural relationship $\Gamma_{r}:=\gamma_{\overline{r}}\circ\varsigma^{-1}$. The Markov semigroup $\pi_t^{\epsilon}(p,dq)$ of the process $q_t:=\varsigma(Q_t)$ is defined for any bounded measurable function $f\in\Ba(\RR^{\overline{r}})$ and any $q\in \varsigma(\Sa^0_r)\subset\RR^{\overline{r}}$ by the formula
$$
\pi_t^{\epsilon}(f)(q)\,:=\,\Pi_t^{\epsilon}(f\circ\varsigma)(\varsigma^{-1}(q)) \qquad\Longleftrightarrow\qquad
\Pi_t^{\epsilon}(F)(Q)\,=\,\pi_t^{\epsilon}(F\circ\varsigma^{-1})(\varsigma(Q))
$$

The symmetric tensor product $P_1\otimes_{ s} P_2$ can be identified with the matrix $\left\{P_1\otimes_{ s} P_2\right\}\in \RR^{\overline{r}\times\overline{r}}$ defined by
\begin{equation}\label{eq-ap-1}
\left\{P_1\otimes_{ s} P_2\right\}:=\varsigma \circ(P_1\otimes_{ s} P_2)\circ \varsigma^{-1}
\quad\Longrightarrow\quad
\left\{P_1\otimes_{ s} P_2\right\}^{1/2}=\varsigma \circ(P_1\otimes_{ s} P_2)^{1/2}\circ \varsigma^{-1}
\end{equation}
and we have the estimate
\begin{equation}\label{eq-ap-2}
\lambda_r(P_1)\,\lambda_r(P_2) ~I~\leq \left\{P_1\otimes_{ s} P_2\right\}\leq \lambda_1(P_1)\,\lambda_1(P_2)~I~
\end{equation}
The proof of the above tensor product formulae are provided in the Appendix.

Finally, define the optimal matching distance between the spectrum of matrices $A,B\in\Ma_r$ by
\begin{equation}\label{optimal-match-d}
d\left(\mathrm{Spec}(A),\mathrm{Spec}(B)\right)=\min_{\mathrm{perm(\cdot)}}\,{\max_{1\leq i\leq r}\vert \lambda_i(A)-\lambda_{\mathrm{perm}(i)}(B)\vert}
 \end{equation}
where the minimum is taken over the set of $r!$ permutations of $\{1,\ldots,r\}$. Recall also the Krause \cite{KRAUSE199473} and Friedland \cite{friedland1982variation} inequalities,
\begin{equation}\label{krause-ref}
 d\left(\mathrm{Spec}(A),\mathrm{Spec}(B)\right)\vee \vert\mbox{\rm det}(A)-\mbox{\rm det}(B)\vert^{1/r}~\leq~ c\,\left[
 \Vert A\Vert\vee\Vert B\Vert \right]^{1-1/r}\,\Vert A-B\Vert^{1/r}
\end{equation}
for any $A,B\in\Ma_{r}$. For any $A,B\in\Sa_r$ we also have \cite{fiedler2008special} the Hoffman-Wielandt inequality 
\begin{equation}\label{hw-ref}
\sum_{1\leq i\leq r}\left(\lambda_i(A)-\lambda_i(B)\right)^2~\leq~ \Vert A-B\Vert^2
\end{equation}

\section{Formal Statement of the Main Results: Regularity and Stability}\label{statement-sec}

Recall that throughout we make the standing assumption that $(A,R^{1/2})$ is stabilizable and $(A,S^{1/2})$ is detectable; see \cite{Kwakernaak72,anderson79,Lancaster1995}.

\subsection{Regularity Properties and Fluctuation Estimates}

\begin{theo}\label{theo-existence-s-ric}
For any $\epsilon\geq 0$ the Riccati diffusion (\ref{f21}) has an unique weak solution on $\Sa^0_r$. For
 $\epsilon\leq \varepsilon_0$ there exists an unique strong solution on $\Sa^+_r$. 
In this situation, we have
 \begin{equation}\label{equivalence-intro}
dQ_t\stackrel{ law}{=}\Theta(Q_t)\,dt+\epsilon\,(Q_t\otimes_s \Sigma_{\kappa,\varpi}(Q_t))^{1/2}~d\Va_{t,\mathrm{sym}}
\end{equation}
where $\Va_{t,\mathrm{sym}}$ denotes a symmetric Brownian matrix with entries $$\Va_{t,\mathrm{sym}}(i,i)=\mathcal{W}_t(i,i)\quad\mbox{and} \quad\Va_{t,\mathrm{sym}}(i,j)=\mathcal{W}_t(i,j)/\sqrt{2}\quad
\mbox{ for any $i<j$.}
$$
When $\epsilon\leq \varepsilon_0$, the process $q_t:=\varsigma(Q_t)\in \Da_{\overline{r}}$ satisfies the $\overline{r}$-dimensional diffusion equation
 \begin{equation}\label{equivalence-intro-2}
dq_t=\theta(q_t)\,dt+\epsilon\,\sigma(q_t)\,dv_t\quad \mbox{with}\quad  \theta=\varsigma \circ\Theta\circ\varsigma^{-1}\quad\mbox{and}\quad
 \sigma(q):= 
\left\{\varsigma^{-1}(q)\otimes_s \Sigma_{\kappa,\varpi}(\varsigma^{-1}(q))\right\}^{1/2}
\end{equation}
where $v_t$ denotes an $\overline{r}$-dimensional Brownian motion. In addition, there exists a smooth positive density $\rho^{\epsilon}\in C^{\infty}(]0,\infty[\times\Da^2_{\overline{r}})$ such that for any $t>0$ and $p\in  \Da_{\overline{r}}$ we have 
 \begin{equation}\label{density-intro}
\pi^{\epsilon}_t(p,dq) \,=\, \rho_t^{\epsilon}(p,q)\,\gamma_{\overline{r}}(dq)
\end{equation}
  \end{theo}
 
  The proof of the above Theorem is provided in Section~\ref{theo-existence-s-ric-proof}. Using (\ref{density-intro}) we check that $\pi_t^{\epsilon}(p,dq)$ and $\Pi^{\epsilon}_t(P,dQ)$ are strongly Feller and irreducible semigroups. Thus, they have  an unique invariant probability measure $\gamma_{\infty}^{\epsilon}$ and $\Gamma_{\infty}^{\epsilon}$ on $\Da_{\overline{r}}$ and $\Sa^+_r$. In addition $\gamma_{\infty}^{\epsilon}$ and $\Gamma_{\infty}^{\epsilon}$ have a positive density w.r.t. $\gamma_{\overline{r}}$ and $\Gamma_{r}$.

The next theorem concerns some time-uniform moment estimates on the stochastic Riccati flow in (\ref{f21}) itself and on its inverse flow. 

\begin{theo}\label{theo-existence-s-ric-proof-bis}
Assume that $S\in\Sa_r^+$. In this situation, for any $n\geq 1$, $t\geq 0$, and any $\epsilon_1< \varepsilon_n(V)$ and $\epsilon_2< \varepsilon_n(U,V)$ we have the uniform estimates
\begin{equation}\label{trace-Phi-inverse}
\vertiii{\phi^{\epsilon_1}_t(Q)}_n
\,\leq\, c_{n}\,(1+\Vert Q\Vert)
\quad\mbox{and}\quad
\vertiii{\phi^{-\epsilon_2}_t(Q)}_n\,\leq\,  c_{n}\,\left(1+
\Vert Q\Vert+ \Vert Q^{-1}\Vert\right)
\end{equation}
Furthermore, for any time horizon $t\geq \upsilon>0$ we also have the uniform estimates
\begin{equation}\label{trace-Phi-inverse-bis}
\vertiii{\phi^{\epsilon_1}_t(Q)}_n\,\leq\, c_{\upsilon,n}~
\quad\mbox{and}\quad
\vertiii{\phi^{-\epsilon_2}_t(Q)}_n\,\leq\,  c_{\upsilon,n}
\end{equation}
In addition, if $\kappa=0$, then for any $\epsilon\geq 0$, $t\geq 0$ and any $s\geq \upsilon>0$  we have the refined estimates,
\begin{equation}\label{ref-trace-unif-intro}
\vertiii{\phi^{\epsilon}_t(Q)}_n \,\leq\, c\,(1+ \Vert Q\Vert)\,(1+\epsilon~\sqrt{n})
\quad\mbox{and}\quad
\vertiii{\phi^{\epsilon}_s(Q)}_n \,\leq\, c_{\upsilon}\,(1+\epsilon~\sqrt{n})
\end{equation}

\end{theo}

The proof of the above Theorem is provided in Section~\ref{theo-existence-s-ric-bis-proof}. A more precise description of the parameters $c_{n},c_{\upsilon,n},c,c_{\upsilon}$  are provided in (\ref{trace-estimates})
and in (\ref{trace-estimates-inverse}). The first estimates stated in (\ref{trace-Phi-inverse}) also hold if $S=V=0$ when $\mu(A)<0$. The proof of this Theorem is based on a reduction of (\ref{f21}) to a scalar Riccati diffusion, a novel representation of its $n$-th powers, and a comparison of its moments to a judiciously designed deterministic scalar Riccati equation. We note the proof is conservative by nature (due to the scalar reduction and comparison).

From (\ref{trace-Phi-inverse-bis}), for any $t\geq \upsilon>0$ there exists some matrices $\overline{\Phi}^\epsilon_{\upsilon},\underline{\Phi}^\epsilon_{\upsilon}>0$ such that
$$
\epsilon \,\leq\, \varepsilon_{1}(V)\wedge\varepsilon_{1}(U,V) \quad\Longrightarrow\quad
 \underline{\Phi}^\epsilon_{\upsilon} \,\leq\, \EE\left(\phi^{\epsilon}_t(Q)\right)\leq \overline{\Phi}^\epsilon_{\upsilon}
 $$
This estimate in a sense generalises the well-known bounds $\underline{\Phi}_{\upsilon}\leq \phi_t(Q)\leq \overline{\Phi}_{\upsilon}$
for some $\overline{\Phi}_{\upsilon},\underline{\Phi}_{\upsilon}>0$ and $t\geq \upsilon>0$; see e.g. \cite{Bishop/DelMoral:2016,bd-CARE}.

Note that the uniform estimates independent of the initial condition stated throughout, involve some arbitrarily small, positive time parameter $\upsilon$, which is related to the notion of a so-called observability/controllability interval; for further details on this topic we refer to \cite{Bishop/DelMoral:2016}.

Now we turn to quantifying the fluctuations of the matrix Riccati diffusions around their limiting values when the diffusion parameter $\epsilon$ tends to $0$. The next theorem extends (in some directions) the uniform fluctuation estimates obtained in~\cite{Bishop/DelMoral/Niclas:2017}. In some results in \cite{Bishop/DelMoral/Niclas:2017} time-uniform estimates were obtained only with $A$ stable, whereas here we accommodate more general matrix models with possibly unstable modes.

\begin{theo}\label{theo-3-intro}
We have that (\ref{control-intro-bias}) holds. Assume further that $S\in\Sa_r^+$. In this situation, for any time horizon $t\geq 0$ and $\epsilon\leq \varepsilon_{10}(V)$ we also have the refined uniform bias estimates
\begin{equation}\label{ref-vp-max}
0~\leq\, \phi_{t}\left(Q\right)- \EE\left[\phi^{\epsilon}_t(Q)\right]\,\leq\, c\,\epsilon^2\, (1+\Vert Q\Vert^5)\,(\lambda_1(U)+\lambda_1(V)\,\Vert Q\Vert^2)\, I 
\end{equation}
Furthermore, for any $n\geq 1$  and $\epsilon\leq \varepsilon_{10n}(V)$ we have the uniform estimates
\begin{equation}
\vertiii{\phi^{\epsilon}_t(Q)-\phi_{t}(Q)}_n \,\leq\, c_n\,\epsilon\,(1+ \Vert Q\Vert^7)
\label{ref-phi-1-est}
\end{equation}
 In addition, if $\kappa=0$, then for any $\epsilon\geq 0$ and $n\geq 1$  we have
 \begin{equation}\label{ref-V-0-1}
 \vertiii{\phi^{\epsilon}_t(Q)-\phi_{t}(Q)}_n
\,\leq\, c\, \epsilon\,(1+ \Vert Q\Vert^5)\,(1+\epsilon~\sqrt{n})^5
\end{equation}
\end{theo}
The proof of the preceding Theorem is provided in Section~\ref{proof-theo-3-intro}. Recall that the proof of (\ref{control-intro-bias}) is an easy consequence of the inequality (\ref{ref-phi-1-est-proof}), see also Theorem 1.3 in~\cite{Bishop/DelMoral/Niclas:2017}, and the uniform estimates in~\cite{Bishop/DelMoral:2016,bd-CARE}. The rest of the proof is based simply on a second-order expansion of the stochastic flow $\phi_t^\epsilon$ about the deterministic flow $\phi_t$ and then an appropriate bounding of the first and second order stochastic terms; see also \cite{Bishop/DelMoral/Niclas:2017} for details on the decomposition of $\phi_t^\epsilon$ in terms of $\phi_t$ plus stochastic terms of any desired order (in $\epsilon$).

Observe that the condition $S\in\Sa_r^+$ ensures that for any $n\geq 1$, the $n$-th moments of the trace of the Riccati diffusion are uniformly bounded w.r.t. the time horizon (when the fluctuation parameter is small enough) even when {\em the matrix $A$ is unstable}. This condition may be thought of as a strengthening of the detectability/observability condition.

Several spectral estimates can be deduced from the estimates (\ref{ref-vp-max}), (\ref{ref-phi-1-est}) and (\ref{ref-V-0-1}). For example, let $\kappa=0$ and $\epsilon\in [0,\varepsilon_0]$, then combining (\ref{ref-V-0-1}) with the $n$-version of the Hoffman-Wielandt inequality (\ref{hw-ref}) we have the uniform estimates
$$
\sup_{1\leq i\leq r}\,\vertiii{\lambda_i\left( \phi^{\epsilon}_t(Q)\right)-\lambda_i\left(\phi_t(Q)\right)}_{n} \,\leq\, c\,\epsilon~(1+ \Vert Q\Vert^5)\,(1+\epsilon~\sqrt{n})^5
$$

The uniform estimates (\ref{ref-phi-1-est}) can also be used to analyse the fluctuation of the inverse flow $\phi^{-\epsilon}_t(Q)$.
For instance, for any $n\geq 1$ and any $\epsilon\leq \varepsilon_{2n}(U,V)\wedge \varepsilon_{20n}(V)$ we have the uniform 
 \begin{equation}\label{ref-app-1}
\vertiii{\phi^{-\epsilon}_t(Q)- \phi_{t}\left(Q\right)^{-1}}_n \,\leq\, c_n\,\epsilon\,(1+ \Vert Q\Vert+\Vert Q^{-1}\Vert)^8
\end{equation}
A proof of this estimate is provided in the Appendix.

\subsection{Stability Estimates}

We set $\Lambda(P):=\Vert P\Vert_2+\Vert P^{-1}\Vert_2$ and we consider the collection of 
$\Lambda$-norms on the set of probability measures $\Gamma_1,\Gamma_2\in \Pa(\Sa_r^+)$ on $\Sa_r^+$, indexed by $\iota>0$, and defined by
 $$
 \Vert \Gamma_1-\Gamma_2\Vert_{\iota,\Lambda}:=\sup{\left\{\vert \Gamma_1(F)-\Gamma_2(F)\vert~:~F\in\Ba(\Sa^+_r)~~s.t.~~\Vert F\Vert_{\Lambda}:=\sup_{P\in \Sa_r^+}\frac{\vert F(P)\vert}{1+\iota\, \Lambda(P)}\leq 1\right\}}
 $$
  \begin{theo}\label{theo-stab-intro}
  Assume that the fluctuation parameter  $\epsilon\leq \varepsilon_{1}(V)\wedge\varepsilon_{1}(U,V)$. Then, there exists some parameters $\alpha<\infty$ and $\beta,\iota>0$ such that for any $t\geq 0$ and probability measures $\Gamma_1,\Gamma_2\in \Pa(\Sa_r^+)$ we have the $\Lambda$-norm contraction inequality
 \begin{equation}\label{expo-decays}
 \Vert \Gamma_1\,  \Pi^{\epsilon}_{t}-\Gamma_2\,  \Pi^{\epsilon}_{t}\Vert_{\iota,\Lambda} ~\leq~ \alpha\,e^{-\beta\, t}\,\,\Vert \Gamma_1-\Gamma_2\Vert_{ \iota,\Lambda}
 \end{equation}
 \end{theo}
 The proof of the above theorem is provided in Section~\ref{theo-stab-intro-proof} and is based on matrix-valued Lyapunov (choosing the function $\Lambda(\cdot)$) and minorisation conditions. 
 
 For one-dimensional models, the article~\cite{2017arXiv171110065B} provides explicit analytical expressions for the reversible measure of $Q_t$ in terms of the model parameters. As expected, heavy tailed reversible measures arise when $\kappa=1$, and weighted Gaussian distributions when $\kappa=0$; see the examples in \cite{2017arXiv171110065B}. The article \cite{2017arXiv171110065B} also provides sharp exponential decay rates to equilibrium, in the sense that the decay rates tend to those of the limiting deterministic Riccati equation when $\epsilon$ tends to $0$.

\subsubsection{Contraction Properties of Exponential Semigroups}\label{sec-stability}

The stochastic flow of the matrix Riccati diffusion (\ref{f21}) is given implicitly by
\begin{equation}\label{closed-form-ricc-diff}
	\phi^\epsilon_t(Q) \,=\, \mathcal{E}_{s,t}^{\epsilon}(Q)\,\phi^\epsilon_{s}(Q)\,\mathcal{E}_{s,t}^{\epsilon}(Q)' + \int_{s}^t\,\mathcal{E}^{\epsilon}_{u,t}(Q)\,\Sigma_{1,0}(\phi^{\epsilon}_{u}(Q))\,\mathcal{E}^{\epsilon}_{u,t}(Q)'\,du +\epsilon\int_s^t~\mathcal{E}_{u,t}^{\epsilon}(Q)dM_u\mathcal{E}_{u,t}^{\epsilon}(Q)^{\prime}
\end{equation}
for any $s\leq t$. This formula intuitively illustrates that the regularity properties of the matrix Riccati diffusion (\ref{f21}) are also intimately connected to the contraction properties of $\mathcal{E}_{s,t}^{\epsilon}(Q)$.

The stability properties of the deterministic ($\epsilon=0$) semigroups $\mathcal{E}_{s,t}(Q)$ and $\phi_t(Q)$ are rather well understood; e.g. see \cite{Bishop/DelMoral:2016,bd-CARE}. We begin this section with a brief review on this topic and some key contraction inequalities. Firstly, for any $t\geq 0$ and $Q\in  \Sa^0_r$ we have
\begin{equation}\label{ref-E-1}
\Vert \mathcal{E}_{t}(Q)\Vert \,\leq\, c\,(1+\Vert Q\Vert)\,\Vert \mathcal{E}_{t}(\mathscr{P}_{\infty})\Vert \quad \mbox{\rm and}\quad \Vert \mathcal{E}_{t}(\mathscr{P}_{\infty})\Vert \,\leq\, \alpha\,e^{-\beta\, t}
\quad\mbox{for some  $\alpha,\beta>0$.}
\end{equation}
with $\mathscr{P}_\infty$ from (\ref{fixed-point-Ricc}). In addition, there exists some parameter $\upsilon> 0$ such that for any $s\geq 0$ and any $t\geq \upsilon>0$ we have the uniform estimates,
\begin{equation}\label{ref-krause-inf}
\Vert \mathcal{E}_{s,s+t}(Q)\Vert \,\leq\,c_\upsilon\,\Vert \Ea_{s}(\mathscr{P}_{\infty})\Vert 
\end{equation}
Proof of the above inequalities follows from~\cite{Bishop/DelMoral:2016,bd-CARE}.

Let $P_1,P_2\in \Sa^0_r$. Then, for any $t\geq 0$, using (\ref{ref-E-1}) we find 
\begin{equation}\label{ref-phi-1-stability}
\Vert \phi_{t}(P_1) -  \phi_{t}(P_2)\Vert \,\leq\, c\,(1+\Vert P_1\Vert^2+\Vert P_2\Vert^2)\,\Vert \mathcal{E}_{t}(\mathscr{P}_{\infty})\Vert^2\,\Vert P_1 -  P_2\Vert
\end{equation}
and similarly, using (\ref{ref-krause-inf}), for any $s\geq 0$ and any $t\geq \upsilon>0$, we have 
\begin{equation}\label{ref-phi-1-stability-upsilon}
\Vert \phi_{s,s+t}(P_1) -  \phi_{s,s+t}(P_2)\Vert \,\leq\, c_\upsilon\,\Vert \mathcal{E}_{t}(\mathscr{P}_{\infty})\Vert^2\,\Vert P_1 -  P_2\Vert
\end{equation}
Note that both (\ref{ref-phi-1-stability}) and (\ref{ref-phi-1-stability-upsilon}) imply immediately that $\phi_{t}(Q)\rightarrow_{t\rightarrow\infty}\mathscr{P}_\infty$ for any $Q\in \Sa^0_r$. Again, the proof of these estimates follows from~\cite{Bishop/DelMoral:2016,bd-CARE}. We emphasise that in the deterministic case ($\epsilon=0$), stability of the matrix-valued Riccati differential equation, e.g. as in (\ref{ref-phi-1-stability}), follows directly from the contraction properties of $\mathcal{E}_{s,t}(Q)$ in (\ref{ref-E-1}); see \cite{Bishop/DelMoral:2016,bd-CARE} for the derivation.

We come now to the contractive properties of $\mathcal{E}_{s,t}^{\epsilon}(Q)$. Firstly, we remark that if $S\in\Sa_r^+$, then up to a change of basis we can always assume that $S=I$. Moreover, for any $s,t\in[0,\infty[$ we immediately have the rather crude almost sure estimate
\begin{equation}\label{change-basis-intro-formula}
\mu\left(A\right)<0 \qquad\Longrightarrow\qquad
\Vert \mathcal{E}_{s,s+t}^{\epsilon}(Q)\Vert_2 \,\leq\, \exp{\left[t\,\mu\left(A)\right)\right]}~\longrightarrow_{t\rightarrow \infty}~0
\end{equation}
In general, asking for $A$ to be stable in this form is a very strong and restrictive condition. We typically seek contraction results on $\mathcal{E}_{s,t}^{\epsilon}(Q)$ that accomodate arbitrary $A\in\Ma_r$ matrices. To this end, fix the matrix $Q\in \Sa^0_r$ and consider the process $\Aa^{\epsilon}$ defined by
\begin{equation}\label{def-Aa}
 \Aa^{\epsilon}~:~t\in[0,\infty[\,\mapsto\, \Aa_t^{\epsilon}\,:=\,A-\phi_t^{\epsilon}(Q)S
\end{equation}
We write $\Aa$ for the analogous process driven by $\phi_t(Q)$, i.e. with $\epsilon=0$.

In this notation, for example when $\kappa=0$, combining (\ref{ref-trace-unif-intro}) (\ref{ref-V-0-1}) and (\ref{ref-krause-inf})  with Krause's inequality (\ref{krause-ref}) for any $nr\geq 1$ we also have the uniform estimate
\begin{equation}\label{krause-ref-c}
 \vertiii{\,d
 \left(\mathrm{Spec}(\Aa_t),\mathrm{Spec}(\Aa_t^{\epsilon})\right)\,}_{nr}~\leq~ c_{n}(Q)\,\epsilon
\end{equation}
In addition, for any $t\geq \upsilon>0$, using the Lipschitz estimates discussed above we also have
\begin{equation}\label{mean-spectral-abscissa-proof-ref}
\Aa_{\infty}\,:=\,A-\mathscr{P}_{\infty}S~~~\Longrightarrow~~~
d\left(\mathrm{Spec}(\Aa_t),\mathrm{Spec}(\Aa_{\infty})\right)~\leq~ c\,\exp{\left[-2\beta t/r\right]}\,\Vert Q-\mathscr{P}_{\infty}\Vert^{1/r}
\end{equation}
with the parameter $\beta$ coming from (\ref{ref-E-1}). These spectral estimates are of interest on their own, but are not immediately usable for controlling the contraction properties of the exponential semigroups.

By Theorem~\ref{theo-existence-s-ric-proof-bis} and Theorem~\ref{theo-3-intro}, with $S\in\Sa_r^+$, the collection of processes $(\Aa,\Aa^{\epsilon})$ introduced in (\ref{def-Aa}) satisfy the following regularity properties:
\begin{itemize}
\item {\em Case $\kappa\in\{1,0\}$}: 
For any $n\geq 1$ and $\epsilon\leq \varepsilon_{10n}(V)$ and any $t\geq 0$ we have the uniform estimates
$$
\epsilon^{-1}\,\vertiii{\Aa_t-\Aa^{\epsilon}_t}_{n} \,\leq\, c_n\,(1+\|Q\|^7)
$$
\item {\em Case $\kappa=0$}:  

For any $n\geq 1$ and any $\epsilon\geq 0$ and any $t\geq 0$ we have the uniform estimates
$$
\epsilon^{-1}\,\vertiii{\Aa_t-\Aa^{\epsilon}_t}_{n} \,\leq\, c\,(1+ \|Q\|^5)\,(1+\epsilon\,\sqrt{n})^5~\quad\mbox{\rm and}\quad~
\vertiii{\Aa^{\epsilon}_t}_{n} \,\leq\, c\,(1+ \|Q\|)\,(1+\epsilon\,\sqrt{n})
$$
\end{itemize}
with the parameter $\kappa\in\{0,1\}$ introduced in (\ref{cond-UV-ref-intro}).

The stability properties of stochastic semigroups associated with a collection of stochastic flows $(\Aa,\Aa^{\epsilon})$ 
satisfying the above properties have been developed in our prior work~\cite{Bishop/DelMoral/STV2018}. Several local-type contraction estimates can be derived. For instance, the stochastic semigroup $\mathcal{E}_{s,t}^{\epsilon}(Q)$ exhibits the following stability properties derived as immediate corollaries of our work in \cite{Bishop/DelMoral/STV2018}:

\begin{cor}\label{cor-H1-semigroups}
 Let $\kappa\in\{1,0\}$. The semigroup $\mathcal{E}_{s,t}^{\epsilon}(Q)$ is asymptotically stable in probability if $\mu(\Aa_{\infty})<0$. That is, for any increasing sequence of times $0\leq s \leq t_k\uparrow_{k\rightarrow\infty}\infty$, the probability of the following event
 \begin{equation}\label{EA-epsilon-H1}
\limsup_{k\rightarrow\infty}\frac{1}{t_k}\log{\Vert \Ea_{s,t_k}^{\epsilon}(Q)\Vert} \,< \, \frac{1}{2}\,\mu(A-\mathscr{P}_{\infty}S)\quad \mbox{is greater than $1-\nu$}
 \end{equation}
for any $\nu\in]0,1[$, as soon as $\epsilon^n \leq c_n\,\nu$ for some $n\geq 1$.
\end{cor}

This log-Lyapunov estimate (\ref{EA-epsilon-H1}) immediately implies the semigroup $\Ea_{s,t_k}^{\epsilon}(Q)$ may be exponentially contracting with a high probability; given strong observability and controllability conditions that imply $\mu(A-\mathscr{P}_\infty S)<0$. A number of reformulations of this result that shed insight individually are worth highlighting:
\begin{itemize}
\item Let $\kappa\in\{1,0\}$. For any $0\leq s \leq t_{k_1}\uparrow_{{k_1}\rightarrow\infty}\infty$, there exists a sequence $\epsilon_{k_2}\downarrow_{{k_2}\rightarrow\infty} 0$ such that we have the almost sure Lyapunov estimate
 \begin{equation}\label{EA-epsilon-H1-dd}
\limsup_{{k_2}\rightarrow\infty}\limsup_{{k_1}\rightarrow\infty}\frac{1}{t_{k_1}}\,\log{\Vert \Ea_{s,s+t_{k_1}}^{\epsilon_{k_2}}(Q))\Vert}\,<\, \frac{1}{2}\,\mu(A-\mathscr{P}_{\infty}S)
 \end{equation} 
\item Let $\kappa\in\{1,0\}$. Then, for any increasing sequence of times $0\leq s \leq t_k\uparrow_{k\rightarrow\infty}\infty$, the probability of the following event,   
 \begin{align}
 \left\{\begin{array}{l}
 \forall 0<\nu_2\leq 1~~~ \exists l\geq 1 ~~~\mbox{such~that}~~~ \forall k\geq l~~~\mbox{it~holds~that~} \\~
 \\
 \qquad\qquad\qquad\qquad\qquad\qquad\qquad\displaystyle\frac{1}{t_k}\log{\Vert \Ea_{s,t_k}^{\epsilon}(Q)\Vert} \,\leq\,   \frac{1}{2}\,(1-\nu_2)\,\mu(A-\mathscr{P}_{\infty}S)
 \end{array}\right\} \label{EA-epsilon-H1-cor}
\end{align}
is greater than $1-\nu_1$, for any $\nu_1\in]0,1[$, as soon as $\epsilon^n \leq c_n\,\nu$ for some $n\geq 1$.
 
\item Let $\kappa\in\{1,0\}$. Consider any $s\geq 0$, any increasing sequence of time horizons $t_k\uparrow_{{k_1}\rightarrow\infty}\infty$, and any sequence $\epsilon_{k_2}\downarrow_{{k_2}\rightarrow\infty} 0$ such that $\sum_{{k_2}\geq 1}\epsilon_{k_2}^{\,n}<\infty$ for some $n\geq 1$. Then, we have the almost sure Lyapunov estimate, 
\begin{align}
 \left\{\begin{array}{l}
 \forall 0<\nu\leq 1~~~ \exists l_1,l_2\geq 1 ~~~\mbox{such~that}~~~ \forall k_1\geq l_1,~\forall k_2\geq l_2~~~\mbox{it~holds~that~} \\~
 \\
 \qquad\qquad\qquad\qquad\qquad\qquad\quad\displaystyle \frac{1}{t_{k_1}}\log{\Vert \Ea_{s,s+t_{k_1}}^{\epsilon_{k_2}}(Q)\Vert} \,\leq\,   \frac{1}{2}\,(1-\nu)\,\mu(A-\mathscr{P}_{\infty}S)
\end{array}\right\} \label{EA-epsilon-H1-cor-v2}
\end{align}
\end{itemize}

The first dot-point result captured by (\ref{EA-epsilon-H1-dd}) is derived from (\ref{EA-epsilon-H1}) in Corollary \ref{cor-H1-semigroups} via the Borel-Cantelli lemma. The next two dot-point results provide some reformulation of the supremum limit estimates (\ref{EA-epsilon-H1}) and (\ref{EA-epsilon-H1-dd}) in terms of random relaxation time horizons and random relaxation-type fluctuation parameters. The first reformulation in (\ref{EA-epsilon-H1-cor}) shows that with a high probability, the semigroup $\Ea_{s,t}^{\epsilon}(Q)$ is stable after some possibly random relaxation time horizon, as soon as $\epsilon$ is chosen sufficiently small and $\mu(\Aa_{\infty})<0$. The last reformulation in (\ref{EA-epsilon-H1-cor-v2}) underlines the fact that after some random time (i.e. determined by $l_1$), and given some randomly sufficiently small perturbation (determined by $l_2$) the semigroup $\Ea_{s,t}^{\epsilon}(Q)$ is exponentially contractive. We have no direct control over the parameters $l_1$ and $l_2$ in (\ref{EA-epsilon-H1-cor-v2}) which depend on the randomness in any realisation.

Additional results are applicable if we restrict $\kappa=0$. We have the following immediate corollary of our prior work in \cite{Bishop/DelMoral/STV2018}:

\begin{cor}\label{cor-H2-semigroups}
 Assume $\kappa=0$. If $\mu(\Aa_{\infty})<0$, then the semigroup $\mathcal{E}_{s,t}^{\epsilon}(Q)$ is asymptotically $\LL_n$-stable for any $n\geq 1$ over time horizons with lengths controlled by $\epsilon$. More specifically, for any $n\geq 1$, $s\geq0$, there exists some time horizons $\mathfrak{t}_n<\mathfrak{t}_n^\epsilon\longrightarrow_{\epsilon\rightarrow 0}\infty$ such that for any $\mathfrak{t}_n\leq t\leq \mathfrak{t}_n^\epsilon$ we have
 \begin{equation}\label{EA-epsilon-H2}
 \frac{1}{t}\log{\EE\left(\Vert \mathcal{E}_{s,s+t}^{\epsilon}(Q)
\Vert^{n}\right)} \,\leq\, 
 \frac{n}{4}\,\mu(A-\mathscr{P}_{\infty}S) 
  \end{equation}
whenever $\epsilon\leq\varepsilon_{n,\mathfrak{t}}$ where here $\varepsilon_{n,\mathfrak{t}}$ is the largest parameter $\epsilon$ such that $\mathfrak{t}_n^\epsilon>\mathfrak{t}_n$; see \cite{Bishop/DelMoral/STV2018} for more details on these time parameters.
\end{cor}

Importantly, in this last result we have $\mathfrak{t}_n^\epsilon\longrightarrow_{\epsilon\rightarrow 0}\infty$ and thus we can control the horizon on which the semigroup $\mathcal{E}_{s,t}^{\epsilon}(Q)$ is asymptotically $\LL_n$-stable for any $n\geq 1$ when $\kappa=0$. In other words, the estimate (\ref{EA-epsilon-H2}) ensures that the stochastic semigroup $\mathcal{E}_{s,t}^{\epsilon}(Q)$ is stable {\em on arbitrary long finite time horizons}, as soon as $\kappa=0$, and when the perturbation parameter is chosen sufficiently small. We have the following fact immediate from Corollary \ref{cor-H2-semigroups}:

\begin{itemize}
\item Assume $\kappa=0$. For any $n\geq 1$, $s\geq0$, we have
$$
\limsup_{\epsilon\rightarrow0}\,
\frac{1}{\mathfrak{t}_n^\epsilon}\,\log{\EE\left(\Vert \Ea_{s,s+\mathfrak{t}_n^\epsilon}^{\epsilon}(Q)\Vert^{n}\right)} \,\leq \,
 \frac{n}{4}\,\mu(A-\mathscr{P}_{\infty}S) 
$$
\end{itemize}

Finally, we also have the following new result which extends the exponential decay results for one-dimensional models presented in~\cite{2017arXiv171110065B} to the determinant of the matrix-valued Riccati diffusions considered herein.

\begin{theo}[Stochastic Liouville Formula]\label{det-E-theo}
For any $Q\in \Sa_r^+$ consider the parameters $n>1$ and $
\epsilon\leq \varepsilon_{2n r}(V)$ such that
 \begin{equation}\label{def-RS-repsilon}
R^\epsilon_{n} \,:=\, R^{\epsilon}-n\frac{\epsilon^{2}}{2}U \,>\, 0
\quad \mbox{and}\quad
S^\epsilon_{n} \,:=\, S^{\epsilon}-n\frac{\epsilon^{2}}{2}V \,>\, 0
\end{equation}
where $(U,V)$ and $(R^\epsilon,S^\epsilon)$ are defined as in (\ref{cond-UV-ref-intro-uv}) and (\ref{Hyp-RS-repsilon}). Then, we have the exponential decay estimate
\begin{equation}\label{ars-0}
\EE\left[\mbox{\rm det}(\Ea_t^{\epsilon}(Q))^{n}\right]^{1/n}\,=\,\EE\left(\exp{\left[n\int_0^t~\tr(A-\phi_s^{\epsilon}(Q)S)\,ds\right]}\right)^{1/n}
\displaystyle\,\leq\, c_{n}(Q)\,
\exp{\left[
-t\,\sqrt{\tr\left(R^\epsilon_{n}S^\epsilon_{n}\right)}
\right]}
\end{equation}
In addition,  the exists some function $\lim_{\epsilon\rightarrow 0}\hbar_{n}(\epsilon)=0$ such that
\begin{equation}\label{ars}
\EE\left[\mbox{\rm det}(\Ea_t^{\epsilon}(Q))^{n}\right]^{1/n}
\displaystyle\,\leq\,  c_{n}(Q)\,
\exp{\left[
-t\,\sqrt{\tr(A)^2+\tr(RS)}(1-\hbar_{n}(\epsilon))
\right]}
\end{equation}
\end{theo}
The proof of this theorem is given in Section~\ref{det-E-theo-proof}. In the one-dimensional case, $r=1$, this result collapses to capture the strong exponential contraction results presented in \cite{2017arXiv171110065B} on the semigroups $\Ea_t^{\epsilon}$ associated with a scalar-valued Riccati diffusion. In the scalar case, strong stability results on the stochastic Riccati flow $\phi^\epsilon_t$ analogous to the deterministic setting, e.g. (\ref{ref-phi-1-stability}), also follow; see \cite{2017arXiv171110065B}.

\section{Ensemble Kalman-Bucy Filters}\label{EnKF-sec}

Because of their practical importance, this section is dedicated to the illustration of our main results within the {\tt EnKF} framework \cite{evensen03,sakov2008a}. Consider a time-invariant linear-Gaussian filtering model of the following form,
\begin{equation}\label{lin-Gaussian-diffusion-filtering}
d\mathscr{X}_t=A\,\mathscr{X}_t~dt+R_1^{1/2}\,d\mathscr{W}_t\quad\mbox{\rm and}\quad
d\mathscr{Y}_t=B\,\mathscr{X}_t~dt+R_2^{1/2}\,d\mathscr{V}_{t},
\end{equation}
where $(\mathscr{W}_t,\mathscr{V}_t)$ is an $(r+\overline{r})$-dimensional Brownian motion, $\mathscr{X}_0$ is an $r$-dimensional Gaussian random variable with mean and variance $(\EE(\mathscr{X}_0),P_0)$
(independent of $(\mathscr{W}_t,\mathscr{V}_t)$), $(A,B)\in (\Ma_{r}\times \Ma_{r,\overline{r}})$,  $(R_1,R_2)\in (\Sa^+_{r}\times 
\Sa^+_{\overline{r}})$, $\mathscr{Y}_0=0$. To simplify, and relate the notation here to our general analysis, set 
$$
R:=R_1\quad \mbox{\rm and} \quad S:=B^{\prime}R_2^{-1}B$$

We let $\Ya_t=\sigma\left(\mathscr{Y}_s,~s\leq t\right)$ be the $\sigma$-algebra filtration generated by the observations.
The conditional distribution $\eta_t=\mbox{\rm Law}\left(\mathscr{X}_t~|~\Ya_t\right)$ of the signal internal states $\mathscr{X}_t$ given $\Ya_t$ is a Gaussian distribution with a conditional mean
and a conditional variance  given by
$$
\mathscr{M}_t:=\EE\left(\mathscr{X}_t~|~\Ya_t\right)\quad \mbox{\rm and}\quad
\mathscr{P}_t:=\EE\left(\left[\mathscr{X}_t-\EE\left(\mathscr{X}_t~|~\Ya_t\right)\right]\left[\mathscr{X}_t-\EE\left(\mathscr{X}_t~|~\Ya_t\right)\right]^{\prime}~|~\Ya_t\right).
$$

\subsection{McKean-Vlasov Interpretations}

Ensemble Kalman-Bucy filters can be interpreted as a (non-unique) mean field particle approximation of the Kalman-Bucy filtering equation. To describe with some precision these models we need to introduce some terminology. For any probability measure $\eta$ on $\RR^r$ we let $\Pa_{\eta}$ be the $\eta$-covariance matrix
$$
\Pa_{\eta}:=\eta\left([e-\eta(e)][e-\eta(e)]^{\prime}\right)\quad\mbox{\rm with}\quad\quad e(x)=x
$$
We now consider two \cite{evensen03,sakov2008a} different classes of conditional nonlinear McKean-Vlasov-type diffusion processes
\begin{equation} \label{Kalman-Bucy-filter-nonlinear-ref}
\begin{split}
(1)\qquad d\overline{\Xa}_t~=&~A~\overline{\Xa}_t~dt~+~R^{1/2}~d\overline{W}_t+\Pa_{\overline{\eta}_t}~B^{\prime} R_2^{-1}~\left[
dY_t-\left(
B~\overline{\Xa}_t~dt+R_2^{1/2}~
d\overline{V}_{t}\right)\right]; \\
(2)\qquad  d\overline{\Xa}_t~=&~A~\overline{\Xa}_t~dt~+~R^{1/2}~d\overline{W}_t+\Pa_{\overline{\eta}_t}~B^{\prime} R_2^{-1}~\left[dY_t-B~\left(
\frac{\overline{\Xa}_t+
\overline{\eta}_t(e)}{2}\right)~dt\right]
\end{split}
\end{equation}
In all cases  $(\overline{W}_t, \overline{V}_t,\overline{\Xa}_0)$ are independent copies of $(\mathscr{W}_t, \mathscr{V}_t,\mathscr{X}_0)$ (thus independent of
 the signal and the observation path) and
 \begin{equation}\label{def-nl-cov}
\overline{\eta}_t=
\mbox{\rm Law}(\overline{X}_t~|~\Ya_t).
\end{equation}
These diffusions are time-varying Ornstein-Uhlenbeck processes \cite{DelMoral/Tugaut:2016} and consequently $ \overline{\eta}_t$ is Gaussian; see also \cite{Bishop/DelMoral:2016}.
These Gaussian distributions have the same conditional mean $\mathscr{M}_t=\overline{\eta}_t(e)$ and the same conditional variance $\mathscr{P}_t=\Pa_{\eta_t}=\Pa_{\overline{\eta}_t}$. They satisfy
the Kalman-Bucy filter
 \begin{equation}\label{def-nl-cov-evol}
d\mathscr{M}_t=(A-\mathscr{P}_tS)~\mathscr{M}_t~dt+\mathscr{P}_t~B^{\prime} R_2^{-1}~d\mathscr{Y}_t\quad\mbox{\rm with the Riccati equation}\quad
\partial_t\mathscr{P}_t=\Theta\left(\mathscr{P}_t\right).
\end{equation}
Ensemble Kalman-Bucy filters coincide with the mean-field particle interpretation of the nonlinear diffusion processes defined in \eqref{Kalman-Bucy-filter-nonlinear-ref}. To be more precise, let $(\overline{W}^i_t,\overline{V}^i_t,\overline{X}_0^i)_{1\leq i\leq N+1}$ be $(N+1)$ independent copies of $(\overline{W}_t,\overline{V}_t,\overline{X}_0)$.

The {\tt EnKF} associated with the nonlinear processes $\overline{\Xa}_t$ defined in (\ref{Kalman-Bucy-filter-nonlinear-ref}) are given by the Mckean-Vlasov-type interacting diffusion processes
\begin{equation}\label{fv1-3}
\begin{split}
(1)\qquad d\overline{\Xa}_t^i~=&~A~\overline{\Xa}_t^i~dt+R^{1/2}~d\overline{W}_t^i+\widehat{\mathscr{P}}_t~B^{\prime} R_2^{-1}
\left[d\mathscr{Y}_t-\left(B ~\overline{\Xa}_t^i~ dt+\Sigma^{1/2}~d\overline{V}^i_{t}\right)\right] \\
(2)\qquad  d\overline{\Xa}_t^i~=&~A~\overline{\Xa}_t^i~dt+R^{1/2}~d\overline{W}_t^i+\widehat{\mathscr{P}}_t~B^{\prime} R_2^{-1}\left[
d\mathscr{Y}_t-2^{-1}~B~\left(\overline{\Xa}_t^i+\overline{\eta}^N_t(e)\right) ~dt\right]\qquad
\end{split}
\end{equation}
with $1\leq i\leq N+1$, $N\geq1$, and the rescaled particle variance
\begin{equation}\label{fv1-3-2}
\begin{array}{l}
\displaystyle \widehat{\mathscr{P}}_t:=\left(1+N^{-1}\right)\,\Pa_{\overline{\eta}^{N}_t}
~~\quad\mbox{\rm with}\quad~~
\displaystyle\overline{\eta}^{N}_t:=(N+1)^{-1}\sum_{1\leq i\leq N+1}\delta_{\overline{X}_t^i}.
\end{array}\end{equation}
Following the arguments as those provided in the beginning of Section~\ref{theo-existence-s-ric-proof}, we can check that the interacting diffusions discussed above have a unique weak solution on $\RR^r$ for any time horizon.

Let $\widehat{\mathscr{M}}_t=\overline{\eta}^{N}_t(e)$ be the particle estimate of the conditional mean $\mathscr{M}_t$.
From~\cite{DelMoral/Tugaut:2016}, and via the representation Theorem (e.g. Theorem 4.2~\cite{karatzas}; see also~\cite{doob}), there exists a filtered probability space enlargement under which we have
\begin{equation} \label{EnKF-1}
\begin{split}
d\widehat{\mathscr{P}}_t~=&~\Theta(\widehat{\mathscr{P}}_t)~dt+\epsilon~\left(
\widehat{\mathscr{P}}_t^{1/2}~d\Wa_t~\Sigma^{1/2}_{\kappa,0}(\widehat{\mathscr{P}}_t)\right)_{\mathrm{sym}}
 \\
d\widehat{\mathscr{M}}_t~=&~(A-\widehat{\mathscr{P}}_tS)~\widehat{\mathscr{M}}_t~dt+\widehat{\mathscr{P}}_t~B^{\prime} R_2^{-1}~d\mathscr{Y}_t+ \overline{\epsilon}~\Sigma^{1/2}_{\kappa,0}(\widehat{\mathscr{P}}_t)~d\widehat{\Wa}_t,
 \quad
\end{split}
\end{equation}
with the parameters
$$
\epsilon:=\frac{2}{\sqrt{N}}~~\Longrightarrow~~
 \overline{\epsilon}:=\frac{\epsilon}{\sqrt{\epsilon^2+4}}=\frac{1}{\sqrt{N+1}} $$
 and from (\ref{cond-UV-ref-intro}) the function
 $$
 \Sigma_{\kappa,0}(Q)\,:=\,R+\kappa\,QSQ\qquad\mbox{\rm with}\quad \kappa=\left\{
\begin{array}{rl}
1&\mbox{\rm in case (1)}\\
0&\mbox{\rm in case (2)}
 \end{array}\right..
$$
In the above display $\Wa_t$ and $\widehat{\Wa}_t$ denotes an $(r\times r)$ and an $r$-dimensional Wiener process respectively, with $\Wa_t$ is independent of $(\mathscr{V}_t,\mathscr{W}_t,\widehat{\Wa}_t)$. Observe that
\begin{eqnarray}
d(\widehat{\mathscr{M}}_t-\mathscr{X}_t)&=&
(A-\widehat{\mathscr{P}}_tS)\,\left(\widehat{\mathscr{M}}_t-\mathscr{X}_t\right)dt+
\widehat{\mathscr{P}}_t~B^{\prime} R_2^{-1/2}~d\mathscr{V}_t-R^{1/2}\,d\mathscr{W}_t+ \overline{\epsilon}~\Sigma^{1/2}_{\kappa,0}(\widehat{\mathscr{P}}_t)~d\widehat{\Wa}_t \nonumber\\
&\stackrel{ law}{=}&(A-\widehat{\mathscr{P}}_tS)\,\left(\widehat{\mathscr{M}}_t-\mathscr{X}_t\right)dt+(\Sigma^{\overline{\epsilon}}_{\kappa,0}(\widehat{\mathscr{P}}_t))^{1/2}\,dW_t  \label{EnKF-2}
\end{eqnarray}
for some $r$-dimensional Wiener process $W_t$ independent of $\Wa_t$ and with $\Sigma^{\overline{\epsilon}}_{\kappa,0}$, as in (\ref{ref-inf-reg-EnKF}), inheriting the parameter $\kappa$ from $\Sigma_{\kappa,0}$ in (\ref{EnKF-1}) in the following manner,
$$
	\Sigma^{\overline{\epsilon}}_{\kappa,0} \,:=\, \Sigma_{1,0} + \overline{\epsilon}^2\, \Sigma_{\kappa,0}
$$

Note that one of the equations in (\ref{fv1-3}) and the sample mean and covariance in (\ref{fv1-3-2}) constitute the (continuous-time) ensemble Kalman-Bucy filter ({\tt EnKF}) methodology for state estimation in (\ref{lin-Gaussian-diffusion-filtering}). The Riccati diffusion in (\ref{EnKF-1}) itself does not explicitly appear in the method/algorithm as applied \cite{evensen03,sakov2008a}; nor does the diffusion describing the flow of the sample mean in (\ref{EnKF-1}). 

In typical applications of the {\tt EnKF}, the dimension $r$ is very large, while $N$ is rather moderately sized for computational reasons \cite{evensen03}. Typically also the signal (and perhaps observation) process in (\ref{lin-Gaussian-diffusion-filtering}) is nonlinear. This nonlinearity requires a straightforward methodological modification in (\ref{fv1-3}), see \cite{evensen03}. The sample covariance, for example, in this latter case will not satisfy a Riccati-type diffusion equation like in (\ref{EnKF-1}) and the analysis of the {\tt EnKF} behaviour in that case is more delicate; e.g. see \cite{dkt17}. We point to the introduction in our prior work \cite{2017arXiv171110065B, Bishop/DelMoral/Niclas:2017} for further literature pointers on the analysis of the {\tt EnKF} and its variants. 

We conclude that $(\widehat{\mathscr{P}}_t,(\widehat{\mathscr{M}}_t-\mathscr{X}_t))$ coincides with the processes $(Q_t,X_t)$ introduced in (\ref{f21}) and (\ref{Stoch-OU}) with (\ref{cond-UV-ref-intro}) and (\ref{ref-inf-reg-EnKF}) given in the forms noted above (with $\varpi=0$ and $\kappa\in\{0,1\}$, which switches if case (1) or case (2) is considered). Thus, all the (fluctuation and stability) estimates  in Section~\ref{statement-sec} on $\widehat{\mathscr{P}}_t$, and the exponential semigroup generated by $(A-\widehat{\mathscr{P}}_tS)$, apply immediately to this class of state estimator. We consider the stability of $(\widehat{\mathscr{M}}_t-\mathscr{X}_t)$ more explicitly later. We also underline that
$$
(\widehat{\mathscr{M}}_t-\mathscr{X}_t)-(\mathscr{M}_t-\mathscr{X}_t)=\widehat{\mathscr{M}}_t-\mathscr{M}_t
$$

In this non-regularised ($\varpi=0$) {\tt EnKF} context, the condition $\epsilon \leq\varepsilon_0$ in (\ref{Hyp-RS-repsilon}) resumes to the (almost) natural condition $N \geq (r+1)$.

\subsection{Regularized Ensemble Kalman-Bucy Filters}

This section is concerned with some applications of the results developed in the article to the analysis of the regularized {\tt EnKF} filters discussed in~\cite{Bishop/DelMoral/Pathiraja:2017}. We only consider an inflation-type regularisation often discussed in the {\tt EnKF} literature \cite{anderson99,hamill01,evensen03}. In its simplest form, the inflation regularisation method involves replacing in (\ref{Kalman-Bucy-filter-nonlinear-ref}) and (\ref{fv1-3}) the covariance matrices $\Pa_{\overline{\eta}_t}$ and $\widehat{\mathscr{P}}_t$ with some inflated matrices $\Pa_{\overline{\eta}_t}+~\varpi I$ and $\widehat{\mathscr{P}}_t+~\varpi I$ for some judiciously chosen parameter $\varpi> 0$. This inflation translates into (\ref{def-nl-cov-evol}) and into (\ref{EnKF-1}) and (\ref{EnKF-2}).

\subsubsection{McKean-Vlasov Diffusions of Type (1)}

In this case $\kappa=1$. Following the proof of Theorem 3.1 in~\cite{DelMoral/Tugaut:2016} we check that the sample covariance matrix $\widehat{\mathscr{P}}_{t,\varpi}$ associated with the regularised interacting particle system (\ref{fv1-3}), satisfies (\ref{f21}) with $\kappa=1$ in (\ref{cond-UV-ref-intro}), and with the replacement of $(A, \Sigma_{1,0})$  in (\ref{def-Riccati-drift}) by $(A_{1,\varpi},\Sigma_{1,\varpi})$ with
$$
A_{1,\varpi}:=(A-\varpi S)
$$
used in in (\ref{f21}). Note that in this case we have
$$
(A_{1,\varpi}-PS)P+P(A_{1,\varpi}-PS)^{\prime}+\Sigma_{1,\varpi}(P)=AP+PA^{\prime}+(R+\varpi^2S)-PSP
$$
Since we accommodate arbitrary matrices $A\in\Ma_r$ in (\ref{f21}), the latter replacement is covered. Thus again, all the estimates presented in Section~\ref{statement-sec} immediately apply to the sample covariance matrix $\widehat{\mathscr{P}}_{t,\varpi}$ for this class of inflated {\tt EnKF} model.

We can comment on the effect of inflation regularisation on the contractive properties of $\Ea_{s,t}^{\epsilon}$; i.e. specifically with the replacement of $A\leftarrow A_{1,\varpi}$ in (\ref{f21}) and in the definition of $\Ea_{s,t}^{\epsilon}$ in (\ref{def-Est}). Arguing as in (\ref{change-basis-intro-formula}), when $S\in\Sa_r^+$, then up to a change of basis we can always assume that $S=I$. Then,
$$
\mu(A)<\varpi~~\Longrightarrow~~~~
\Vert {\Ea}_{s,t}^{\epsilon}(Q)\Vert_2 \,\leq\, \exp{\left[(\mu(A)-\varpi) (t-s)\right]}~\longrightarrow_{(t-s)\rightarrow\infty}~0
$$
which illustrates the added stabilising effects of $\varpi\,I$ in the extreme case in which $\widehat{\mathscr{P}}_{t,\varpi}\,S$ has no stabilising effect at all. Contrast this with (\ref{change-basis-intro-formula}). It is also worth noting, given the contraction estimates in Section \ref{sec-stability}, that,
$$
\mu(A_{1,\varpi} - P\,S)~=~\mu((A-\varpi S) - P\,S)~\leq~\mu(A - P\,S)
$$
for any fixed matrix $P\in\Sa^0_r$ and $S\in\Sa^0_r$.

Note also that the form of (\ref{Stoch-OU}) with (\ref{ref-inf-reg-EnKF}) is immediately applicable in this case. That is, following the proof of Theorem 3.1 in~\cite{DelMoral/Tugaut:2016}, we can check that the evolution of the sample mean $\widehat{\mathscr{M}}_{t,\varpi}$ associated with this class of inflated {\tt EnKF} model is given by,
$$
\begin{array}{l}
d\widehat{\mathscr{M}}_{t,\varpi}=\left[A-(\widehat{\mathscr{P}}_{t,\varpi}+\varpi I)\,S\right]\widehat{\mathscr{M}}_{t,\varpi}\,dt+\left(\widehat{\mathscr{P}}_{t,\varpi}+\varpi~I\right)\,B^{\prime} R_2^{-1}\,d\mathscr{Y}_t+ \overline{\epsilon}\,\Sigma^{1/2}_{1,\varpi}(\widehat{\mathscr{P}}_{t,\varpi})\,d\widehat{\Wa}_t\\
\\
\qquad\quad\Longrightarrow\qquad
d(\widehat{\mathscr{M}}_{t,\varpi}-\mathscr{X}_t)~\stackrel{ law}{=}~\left[A_{1,\varpi}-\widehat{\mathscr{P}}_{t,\varpi}S\right]\left(\widehat{\mathscr{M}}_{t,\varpi}-\Xa_t\right)dt+(\Sigma^{\overline{\epsilon}}_{1,\varpi}(\widehat{\mathscr{P}}_{t,\varpi}))^{1/2}\,dW_t
\end{array}
$$
and thus $(\widehat{\mathscr{M}}_{t,\varpi}-\mathscr{X}_t)$ corresponds exactly with the general form of (\ref{Stoch-OU}) with (\ref{ref-inf-reg-EnKF}) with $\kappa=1$. Later we consider the stability of this process $(\widehat{\mathscr{M}}_{t,\varpi}-\mathscr{X}_t)$.

\subsubsection{McKean-Vlasov Diffusions of Type (2)}

In this case $\kappa=0$. Following the proof of Theorem 3.1 in~\cite{DelMoral/Tugaut:2016} we check that the sample covariance matrix $\widehat{\mathscr{P}}_{t,\varpi}$ associated with the interacting particle system (\ref{fv1-3}), satisfies (\ref{f21}) with $\kappa=0$ in (\ref{cond-UV-ref-intro}), and with the replacement of $A$ given by
$$
	A~\longleftarrow~A_{0,\varpi}:=(A-2^{-1}\varpi S)
$$
used in in (\ref{f21}). Since we accommodate arbitrary matrices $A\in\Ma_r$ in (\ref{f21}), the latter replacement is again already covered. All the estimates presented in Section~\ref{statement-sec} immediately apply to the sample covariance matrix $\widehat{\mathscr{P}}_{t,\varpi}$ of this class of inflated {\tt EnKF} model.

We highlight that in this case, $\Sigma_{0,\varpi}=\Sigma_{0,0}=R$ and thus if $\kappa=0$ then $\varpi$ has no effect in terms of the diffusion matrix. Nevertheless, we may repeat the commentary as in case (1) on the effect of inflation regularisation on the contractive properties of $\Ea_{s,t}^{\epsilon}$; i.e. specifically in this case with the replacement of $A\leftarrow A_{0,\varpi}$ in (\ref{f21}) and in $\Ea_{s,t}^{\epsilon}$.

Now the evolution of the sample mean $\widehat{\mathscr{M}}_{t,\varpi}$ associated with this class of inflated {\tt EnKF} model is given by,  
$$
\begin{array}{l}
d\widehat{\mathscr{M}}_{t,\varpi}=\left[A-(\widehat{\mathscr{P}}_{t,\varpi}+\varpi I)\,S\right]\widehat{\mathscr{M}}_{t,\varpi}\,dt+\left(\widehat{\mathscr{P}}_{t,\varpi}+\varpi~I\right)\,B^{\prime} R_2^{-1}\,d\mathscr{Y}_t+ \overline{\epsilon}\,R^{1/2}\,d\widehat{\Wa}_t\\
\\
\qquad\quad\Longrightarrow\qquad
d(\widehat{\mathscr{M}}_{t,\varpi}-\mathscr{X}_t)~\stackrel{ law}{=}~\left[A_{1,\varpi}-\widehat{\mathscr{P}}_{t,\varpi}\,S\right]\left(\widehat{\mathscr{M}}_{t,\varpi}-\Xa_t\right)dt+(\Sigma^{\overline{\epsilon}}_{0,\varpi}(\widehat{\mathscr{P}}_{t,\varpi}))^{1/2} dW_t
\end{array}
$$
We remark further in this case, that the drift matrix $A_{1,\varpi}$ in the flow of $(\widehat{\mathscr{M}}_{t,\varpi}-\mathscr{X}_t)$ is different to the drift matrix $A\leftarrow A_{0,\varpi}$ in the Riccati diffusion (describing the flow of the sample covariance) in the presence of (non-zero $\varpi>0$) inflation regularisation.

\subsection{Ensemble Filtering Stability Properties}

In this section we consider the stability of the flow of $\psi^{\epsilon}_t(Q,x) := (\widehat{\mathscr{M}}_{t,\varpi}-\mathscr{X}_t)$ in both type (1) and type (2) (possibly regularized) {\tt EnKF} models with $\widehat{\mathscr{P}}_{0,\varpi} = Q\in\Sa_r^0$ and $(\widehat{\mathscr{M}}_{0,\varpi}-\mathscr{X}_0)=x\in\mathbb{R}^r$. This flow may be related to the Ornstein-Uhlenbeck process (\ref{Stoch-OU}) and can be written more generally as, 
\begin{equation}\label{closed-form-OU}
\psi^{\epsilon}_t(Q,x) \,=\, \mathscr{E}_{s,t}^{\epsilon}(Q)\,\psi^{\epsilon}_s(Q,x)+\int_s^t~\mathscr{E}_{u,t}^{\epsilon}(Q)\left(\Sigma^{\overline{\epsilon}}_{\kappa,\varpi}(\widehat{\mathscr{P}}_{u,\varpi})\right)^{1/2}\,dW_u
\end{equation}
Here, $\mathscr{E}_{s,t}^{\epsilon}(Q)$ is a transition matrix associated with the flow of matrices $[A_{1,\varpi}-\widehat{\mathscr{P}}_{t,\varpi}S]$ (defined similarly to $\Ea_{s,t}^{\epsilon}$ in (\ref{def-Est})). This implies the stability properties of the flow of $(\widehat{\mathscr{M}}_{t,\varpi}-\mathscr{X}_t)$ depend on the long time behaviour and contraction properties of the random transition matrices $\mathscr{E}_{s,t}^{\epsilon}(Q)$. If $\varpi=0$ then $A_{1,\varpi}=A$ and $\mathscr{E}_{s,t}^{\epsilon}(Q)=\mathcal{E}_{s,t}^{\epsilon}(Q)$ for any $\kappa\in\{0,1\}$.

Now from preceding results on $\Ea_{s,t}^{\epsilon}$ in Section \ref{sec-stability}, we may comment on the stability of the flow $\psi^{\epsilon}_t(Q,x)= (\widehat{\mathscr{M}}_{t,\varpi}-\mathscr{X}_t)$ in (\ref{closed-form-OU}). Indeed, from our prior work in \cite{Bishop/DelMoral/STV2018}, see also Section \ref{sec-stability}, we have the following stability estimates:

\begin{itemize}
\item Let $\kappa\in\{1,0\}$. For any increasing sequence of time horizons $t_k\uparrow_{k\rightarrow\infty}\infty$ and any $x_1\not=x_2$ and any $Q\in\Sa_r^0$, the probability of the following event
 \begin{equation}\label{EA-epsilon-H1-cor-v3}
\limsup_{k\rightarrow\infty}\frac{1}{t_k}\log{\Vert \psi^{\epsilon}_{t_k}(Q,x_1) - \psi^{\epsilon}_{t_k}(Q,x_2) \Vert} \,<\, \frac{1}{2}\,\mu(A_{1,\varpi}-\mathscr{P}_{\infty,\varpi}\,S)\quad \mbox{is greater than $1-\nu$}
 \end{equation}
 for any $\nu\in]0,1[$, as soon as $\epsilon^n \leq c_n\,\nu$ for some $n\geq 1$. Here, $\mathscr{P}_{\infty,\varpi}$ denotes the unique fixed point satisfying the Riccati matrix map 
\[
	A_{\kappa,\varpi}\mathscr{P}_{\infty,\varpi} + \mathscr{P}_{\infty,\varpi}A_{\kappa,\varpi}^{\prime} + R - \mathscr{P}_{\infty,\varpi}S\mathscr{P}_{\infty,\varpi} = 0
\]

\item Assume $\kappa=0$. For any $n\geq 1$, $\epsilon\leq \varepsilon_{n,\mathfrak{t}}$ and any time horizon $t$ such that $\mathfrak{t}_{n}\leq t\leq \mathfrak{t}_{n}^{\epsilon}$, we have the contraction inequality,
\begin{equation}\label{uniformly-Lipschitz-moments-OU-Lip}
{\EE\left( \Vert \psi^{\epsilon}_t(Q,x_1) - \psi^{\epsilon}_t(Q,x_2) \Vert^n \right)}^{1/n} \,\leq\, \exp{\left[\frac{1}{4}\,t\,\mu(A_{1,\varpi}-\mathscr{P}_{\infty,\varpi}\,S)\right]} \Vert x_1-x_2\Vert
\end{equation}
\end{itemize}

These results concern the flow of the estimation error $(\widehat{\mathscr{M}}_{t,\varpi}-\mathscr{X}_t)$. See \cite{Bishop/DelMoral/STV2018} for further discussion. Note that (\ref{EA-epsilon-H1-cor-v3}) is analogous to (\ref{EA-epsilon-H1}) in Corollary \ref{cor-H1-semigroups} but at the level of the process (\ref{closed-form-OU}) itself. Analogous reformulations as in (\ref{EA-epsilon-H1-dd}), (\ref{EA-epsilon-H1-cor}), and \ref{EA-epsilon-H1-cor-v2}, but on the process (\ref{closed-form-OU}), also follow.

We can comment on the effect of inflation regularisation on the contraction properties of $\mathscr{E}_{s,t}^{\epsilon}(Q)$, as compared e.g. to ${\Ea}_{s,t}^{\epsilon}(Q)$. Arguing as in (\ref{change-basis-intro-formula}), when $S\in\Sa_r^+$, then up to a change of basis we can always assume that $S=I$. We then have,
$$
\mu(A)<\varpi~~\Longrightarrow~~~~
\Vert \mathscr{E}_{s,t}^{\epsilon}(Q) \Vert_2 \,\leq\, \exp{\left[(\mu(A)-\varpi) (t-s)\right]}~\longrightarrow_{(t-s)\rightarrow\infty}~0
$$
which illustrates the added stabilising effects of $\varpi I$ in the extreme case in which $\widehat{\mathscr{P}}_{t,\varpi}\,S$ has no stabilising effect at all. Contrast this with (\ref{change-basis-intro-formula}). In practice, $\widehat{\mathscr{P}}_{t,\varpi}\,S$ will also act to stabilise the filter, see e.g. (\ref{uniformly-Lipschitz-moments-OU-Lip}). Indeed, in the classical Kalman filtering setting (\ref{Stoch-OU-R}) with $\epsilon=0=\varpi$, the time-varying matrix $(A-\mathscr{P}_tS)$ is stabilising \cite{Bishop/DelMoral:2016} for any $A\in\Ma_r$, even $A$ unstable. In the {\tt EnKF}, we know that $\widehat{\mathscr{P}}_{t}$ will fluctuate about $\mathscr{P}_t$, e.g. see Theorem \ref{theo-3-intro}. Therefore, the stabilisation properties of $(A-\widehat{\mathscr{P}}_{t}S)$ are unclear; indeed the study of ${\Ea}_{s,t}^{\epsilon}(Q)$ in the preceding Section \ref{sec-stability} is concerned with precisely this issue. Now the above implies that the addition of $\varpi I$ can act to counter the negative effects of this fluctuation (and directly add a stabilising effect on the state estimation error).

\section{Matrix Riccati Diffusion Flows}\label{sec-Riccati-flows}

In this section we present some general properties and high-level results concerning the matrix Riccati diffusion (\ref{f21}). The results in this section are of interest on their own (and are also used later in the proof of our main results). We still suppose $(A,R^{1/2})$ is stabilizable and $(A,S^{1/2})$ is detectable throughout the remainder.

\subsection{Inverse Matrix Riccati Diffusion Flows}

In our prior discussion and main results we characterise the moments and behaviour of the inverse stochastic flow of (\ref{f21}), which we denoted by $\phi^{-\epsilon}_t(Q):=Q_t^{-1}$. Characterising this flow is important as it allows us to lower bound, in a positive definite sense, moments of the actual stochastic flow (\ref{f21}). This is further required for our main stability results; in an analogous manner to the fact that the inverse deterministic Riccati flow $\phi^{-1}_t(Q)$ is used to study the stability of the deterministic Riccati flow $\phi_t(Q)$; e.g. see \cite{Bishop/DelMoral:2016}.

 Here we characterise the inverse matrix Riccati diffusion and its general structure. This is also likely of interest on its own (e.g. it characterises the so-called flow of the sample information matrix for the ensemble Kalman-Bucy filters).

\begin{lem}
When $\epsilon\leq \varepsilon_0$,  the inverse stochastic flow $Q_t^{-1}$
 satisfies the diffusion equation
\begin{equation}
dQ_t^{-1}\,\stackrel{ law}{=} \, \Theta^\epsilon_{-}(Q_t^{-1})~dt+\epsilon~dM_{t,-}
 \quad\mbox{with}\quad dM_{t,-}:=\left[Q_t^{-1/2}~d\Wa_t~\Sigma_{\kappa,\varpi,-}\left(Q_t^{-1}\right)^{1/2}\right]_{\mathrm{sym}}\label{eq-inverse}
\end{equation}
and where 
$$
\Sigma_{\kappa,\varpi,-}(Q) \,:= \,Q\,\Sigma_{\kappa,\varpi}(Q^{-1})\,Q ~\leq~ V+QUQ
$$
with $(U,V)$ defined as in (\ref{cond-UV-ref-intro-uv}). Here, $\Theta^\epsilon_{-}$ denotes the collection of drift functions satisfying the following inequality
\begin{equation}\label{eq-inverse-drift}
\begin{array}{l}
\displaystyle\Theta^\epsilon_{-}(Q) \,\leq\,-QA-A^{\prime}~Q+S^\epsilon_{-}-QR^\epsilon_{-}Q+\frac{\epsilon^2}{4}\,\left(\tr\left(QU\right)+\tr\left(VQ^{-1}\right)\right)Q
\end{array}
\end{equation}
with the collection of matrices $(R^\epsilon_{-},S^\epsilon_{-})$ defined by
$$
R^\epsilon_{-}\,:=\,R-\frac{\epsilon^2}{4}\,\left(r+2\right)\,U\quad\mbox{and}\quad
S^\epsilon_{-}\,:=\,S+\frac{\epsilon^2}{4}\,\left(r+2\right)\,V
$$
\end{lem}

Note that the specific (equality) form of $\Theta^\epsilon_{-}(\cdot)$ is given in the proof below.

\proof
Note that setting $F(Q):=Q^{-1}$ implies that
$$
\nabla F(Q)\cdot H=-Q^{-1}~H~Q^{-1}\quad \mbox{\rm and}\quad
\frac{1}{2}~ \nabla^2 F(Q)\cdot (H,H)=Q^{-1}~H~Q^{-1}~H~Q^{-1}
$$
Using the Ito differential calculus for stochastic matrix diffusions developed in~\cite{Bishop/DelMoral/Niclas:2017}, with a slight abuse of notation we obtain the formula
\begin{eqnarray*}
dQ^{-1}_t&=&-Q^{-1}_t~\Theta(Q_t)~Q^{-1}_t~dt+\epsilon^2~Q^{-1}_t~dM_t~Q^{-1}_t~dM_t~Q^{-1}_t-\epsilon ~Q^{-1}_t~dM_t~Q^{-1}_t\\
&=&\bigg(\left[-Q^{-1}_tA-A^{\prime}~Q^{-1}_t+S-Q^{-1}_tRQ^{-1}_t\right]~\\
&&\hskip.7cm +\frac{\epsilon^2}{4}~(r+2)~Q^{-1}_t\Sigma_{\kappa,\varpi}\left(Q_t\right)Q^{-1}_t+\frac{\epsilon^2}{4}\tr\left(Q^{-1}_t\Sigma_{\kappa,\varpi}(Q_t)\right)Q^{-1}_t\bigg)\,dt-\epsilon ~Q^{-1}_t~dM_t~Q^{-1}_t
\end{eqnarray*}
The last assertion comes from the decomposition
$$
\begin{array}{l}
4\,\left[Q_t^{1/2}\,d\Wa_t\,\Sigma^{1/2}_{\kappa,\varpi}\left(Q_t\right)\right]_{\mathrm{sym}}Q^{-1}_t\left[Q_t^{1/2}~d\Wa_t~\Sigma^{1/2}_{\kappa,\varpi}\left(Q_t\right)\right]_{\mathrm{sym}}\\
\\
~~=\left[Q_t^{1/2}\,d\Wa_t\,\Sigma^{1/2}_{\kappa,\varpi}\left(Q_t\right)Q_t^{-1/2}+\Sigma^{1/2}_{\kappa,\varpi}\left(Q_t\right)\,d\Wa^{\prime}\right]\left[d\Wa_t\,\Sigma^{1/2}_{\kappa,\varpi}\left(Q_t\right)+Q_t^{-1/2}\,\Sigma^{1/2}_{\kappa,\varpi}\left(Q_t\right)\,d\Wa^{\prime}\,Q_t^{1/2}\right]\\
\\
~~=\left((r+2)~\Sigma_{\kappa,\varpi}\left(Q_t\right)+\tr\left(Q^{-1}_t\Sigma_{\kappa,\varpi}(Q_t)\right)Q_t\right)\,dt
\end{array}
$$
For a more rigorous derivation of the angle bracket of matrix valued martingales we refer the reader to Section 3 in~\cite{Bishop/DelMoral/Niclas:2017}. On the other hand, we have
$$
2~Q^{-1}_t~dM_t~Q^{-1}_t=\left[Q_t^{-1/2}~d\Wa_t~\Sigma^{1/2}_{\kappa,\varpi}\left(Q_t\right)Q^{-1}_t+Q^{-1}_t\Sigma^{1/2}_{\kappa,\varpi}\left(Q_t\right)~d\Wa^{\prime}~Q_t^{-1/2}\right]
$$
Also observe that
\begin{align*}
\left(\Sigma^{1/2}_{\kappa,\varpi}\left(Q_t\right)Q^{-1}_t\right)^{\prime}\left(\Sigma^{1/2}_{\kappa,\varpi}\left(Q_t\right)Q^{-1}_t\right)\,=&\,Q^{-1}_t~\Sigma_{\kappa,\varpi}(Q)Q^{-1}_t \\
:=&\,\Sigma_{\kappa,\varpi,-}\left(Q^{-1}_t\right)\\
=&\,\Sigma^{1/2}_{\kappa,\varpi,-}\left(Q^{-1}_t\right)\Sigma^{1/2}_{\kappa,\varpi,-}\left(Q^{-1}_t\right)\\
\leq &\,Q^{-1}_t~(U+Q_tVQ_t)Q^{-1}_t \,=\,V+Q^{-1}_tUQ^{-1}_t
\end{align*}
This implies that
$$
\Sigma^{1/2}_{\kappa,\varpi}\left(Q_t\right)Q^{-1}_t \Sigma_{\kappa,\varpi,-}^{-1/2}\left(Q^{-1}_t\right)\quad \mathrm{is~an~orthogonal~matrix}
$$
Using the invariance of the matrix Brownian motion by orthogonal transformation this implies that
$$
Q^{-1}_t~dM_t~Q^{-1}_t \,\stackrel{ law}{=} \,\left[Q_t^{-1/2}~d\Wa_t~\Sigma^{1/2}_-\left(Q^{-1}_t\right)\right]_{\mathrm{sym}}
$$
We also have
$$
Q^{-1}_t\Sigma_{\kappa,\varpi}\left(Q_t\right)Q^{-1}_t\leq Q^{-1}_tUQ^{-1}_t+V\quad\mbox{\rm and}\quad
\tr\left(Q^{-1}_t\Sigma_{\kappa,\varpi}(Q_t)\right)\leq 
\tr\left(Q^{-1}_tU\right)+\tr\left(VQ_t\right)
$$
This shows that the drift term $\Theta^\epsilon_{-}(Q^{-1}_t)$ of $Q^{-1}_t$ is given by
$$
\begin{array}{l}
\displaystyle\Theta^\epsilon_{-}(Q^{-1}_t)\\
\\
\,\displaystyle:=-Q^{-1}_tA-A^{\prime}Q^{-1}_t+S-Q^{-1}_tRQ^{-1}_t
+\frac{\epsilon^2}{4}(r+2)Q^{-1}_t\Sigma_{\kappa,\varpi}\left(Q_t\right)Q^{-1}_t+\frac{\epsilon^2}{4}\tr\left(Q^{-1}_t\Sigma_{\kappa,\varpi}(Q_t)\right)Q^{-1}_t\\
\\
\,\,\displaystyle\leq -Q^{-1}_tA-A^{\prime}\,Q^{-1}_t+\left(S+\frac{\epsilon^2}{4}\,(r+2)\,V\right)-Q^{-1}_t\left(R-\frac{\epsilon^2}{4}\,(r+2)\,U\right)Q^{-1}_t\\
\displaystyle\hskip9cm\displaystyle+\frac{\epsilon^2}{4}\,\left(\tr\left(Q^{-1}_tU\right)+\tr\left(VQ_t\right)\right)Q^{-1}_t
\end{array}
$$
This ends the proof of (\ref{eq-inverse}).
\qed

\subsection{A Comparison Lemma}\label{comp-sec-ref}

Here we provide a basic comparison lemma which is useful for deriving moment bounds. For example, we will show subsequently that the left hand side under bias estimate in (\ref{control-intro-bias}), see also (\ref{ref-vp-max}), is a simple consequence of the next lemma. 

\begin{lem}[Comparison Formulae]\label{Lemma-1}
Assume that the flow $t\mapsto \varphi_t(Q)$ satisfies a matrix Riccati-type inequality of the form
$$
\partial_t\varphi_t(Q) \,\leq\, \Theta\left(\varphi_t(Q)\right)
$$
for any $t\geq 0$ any $Q\in \Sa^0_r$. Then, for any times $s\leq t$ and any $P_1,P_2\in \Sa^0_r$ we have the estimate
\begin{equation}\label{comparison-1}
\varphi_t(P_1) \,\leq\, \phi_{t}(P_2)+\mathcal{E}_{s,t}(P_2)\left[\varphi_s(P_1)-\phi_{s}(P_2)\right]\mathcal{E}_{s,t}(P_2)^{\prime}
\end{equation}
and for any $Q\in \Sa^0_r$ we also have the reverse estimate
$$
\partial_t\varphi_t(Q) \,\geq\, \Theta\left(\varphi_t(Q)\right)\quad\Longrightarrow\quad \varphi_t(Q)\,\geq\, \phi_{t}(Q)
$$
\end{lem}

\proof
We recall the polarization-type formulae 
\begin{equation}\label{polarization-formulae}
\begin{array}{l}
\Theta(P_1)-\Theta(P_2) \qquad\\
\\
\qquad\qquad=\left[A-\frac{1}{2}(P_1+P_2)S\right](P_1-P_2)+(P_1-P_2)\left[A-\frac{1}{2}(P_1+P_2)S\right]^{\prime}\\
\\
\qquad\qquad=(A-P_2S)(P_1-P_2)+(P_1-P_2)(A-P_2S)^{\prime}-(P_1-P_2)S(P_1-P_2)\\
\end{array}
\end{equation}
We set 
$$
	\Delta_t\,:=\,\varphi_t(P_1)-\phi_{t}(P_2)
$$
Assume that $\partial_t\varphi_t(Q)\leq \Theta\left(\varphi_t(Q)\right)$. Using the second polarization formula, we have
$$
\partial_t\Delta_t \,\leq\,
\left(A-\phi_t(P_2)S\right)\Delta_t+\Delta_t\left(A-\phi_t(P_2)S\right)^{\prime}-\Delta_tS\Delta_t
$$
On the other hand, for any $s\leq t$ we have
$$
\partial_t \mathcal{E}_{s,t}(P_2)^{-1}=-\mathcal{E}_{s,t}(P_2)^{-1}\,\partial_t \mathcal{E}_{s,t}(P_2)~\mathcal{E}_{s,t}(P_2)^{-1}=-\mathcal{E}_{s,t}(P_2)^{-1}\left(A-\phi_t(P_2)S\right)
$$
This implies that
$$
\begin{array}{l}
\partial_t\left(\mathcal{E}_{s,t}(P_2)^{-1}\Delta_t\,\left(\mathcal{E}_{s,t}(P_2)^{-1}\right)^{\prime}\right)
\,\leq\, -~\mathcal{E}_{s,t}(P_2)^{-1}\Delta_tS\Delta_t\left(\mathcal{E}_{s,t}(P_2)^{-1}\right)^{\prime}
\end{array}
$$
from which we conclude that
$$
\mathcal{E}_{s,t}(P_2)^{-1}\Delta_t\,\left(\mathcal{E}_{s,t}(P_2)^{-1}\right)^{\prime} \,\leq\, \Delta_s
 -\int_s^t~\mathcal{E}_{s,u}(P_2)^{-1}\Delta_uS\Delta_u\left(\mathcal{E}_{s,u}(P_2)^{-1}\right)^{\prime}\,du
$$
$$
\Delta_t \,\leq\, \mathcal{E}_{s,t}(P_2)\Delta_s\,\mathcal{E}_{s,t}(P_2)^{\prime}-\int_s^t\mathcal{E}_{u,t}(P_2)\Delta_uS\Delta_u\mathcal{E}_{u,t}(P_2)^{\prime}~du
$$
This ends the proof of the first assertion. 

We further assume that 
$$
\partial_t\varphi_t(Q) \,\geq\, \Theta\left(\varphi_t(Q)\right)\qquad\mbox{\rm and we let}\qquad
\Delta_t: \,=\,\varphi_t(Q)-\phi_{t}(Q)
$$
Using the first polarization formula, we have
$$
\partial_t\Delta_t \,\geq\,  A_t(Q)~\Delta_t+\Delta_t~A_t(Q)^{\prime}
\quad\mbox{\rm with}\quad
A_t(Q)\,:=\,A-\frac{1}{2}(\varphi_t(Q)+\phi_{t}(Q))S
$$
Let $\widetilde{\Ea}_{s,t}(Q)$ denote the state transition matrix associated with the matrix flow $u\mapsto A_u(Q)$.
Arguing as above, we have
$$
\partial_t \widetilde{\Ea}_{s,t}(Q)^{-1}\,=\,-\widetilde{\Ea}_{s,t}(Q)^{-1}\,\partial_t \widetilde{\Ea}_{s,t}(Q)~\widetilde{\Ea}_{s,t}(Q)^{-1} \,=\,-\widetilde{\Ea}_{s,t}(Q)^{-1}A_t(Q)
$$
This implies that
$$
\begin{array}{l}
\partial_t\left(\widetilde{\Ea}_{s,t}(Q)^{-1}\Delta_t\,\left(\widetilde{\Ea}_{s,t}(Q)^{-1}\right)^{\prime}\right)
\,=\,\widetilde{\Ea}_{s,t}(Q)^{-1}\left[\partial_t\Delta_t-\left(A_t(Q)\Delta_t+\Delta_tA_t(Q)^{\prime}\right)\right]\left(\mathcal{\Ea}_{s,t}(Q)^{-1}\right)^{\prime}
\,\geq\, 0
\end{array}$$
from which we conclude that 
$$
\begin{array}{l}
\displaystyle\widetilde{\Ea}_{s,t}(Q)^{-1}\Delta_t\,\left(\widetilde{\Ea}_{s,t}(Q)^{-1}\right)^{\prime} \,\geq\, \Delta_s
\displaystyle \qquad\Longrightarrow\qquad
\Delta_t \,\geq\, \widetilde{\Ea}_{s,t}(Q)\Delta_s~\widetilde{\Ea}_{s,t}(Q)^{\prime} \,\geq\, 0
\end{array}$$
This ends the proof of the lemma. \qed

We illustrate the impact of the above lemma with a simple proof of the l.h.s. under bias estimate in (\ref{control-intro-bias}), see also the refined estimate in (\ref{ref-vp-max}). Note in (\ref{control-intro-bias}) we do not ask for $S\in\Sa^+_r$ and only require stabilisability and detectability of the model. For any symmetric matrix valued random variable $Q\in\Sa^0_r$ we have
$$
\EE([Q-\EE(Q)]\,S\,[Q-\EE(Q)]) \,\geq\, 0 \qquad\Longleftrightarrow\qquad
\EE(Q\,S\,Q) \,\geq\, \EE(Q)\,S\,\EE(Q)
$$
Using (\ref{comparison-1}), this implies that
\begin{equation}\label{ref-phi-1-est-proof}
\partial_t\EE\left(Q_t\right) \,\leq\, \Theta(\EE\left(Q_t\right)) \qquad\Longrightarrow\qquad \EE\left(\phi^{\epsilon}_t(P_1)\right) \,\leq\,  
\phi_{t}(P_2)+\mathcal{E}_{t}(P_2)~(P_1-P_2)~\mathcal{E}_{t}(P_2)^{\prime}
 \end{equation}
 which immediately implies the left hand side under bias in (\ref{control-intro-bias}). The polarization formula (\ref{polarization-formulae}) also yields the monotone property 
 $$
 P_1 \,\leq\, P_2 \quad\Longrightarrow\quad \EE\left(\phi^{\epsilon}_t(P_1)\right) \,\leq\, \EE\left(\phi^{\epsilon}_t(P_2)\right)
 $$
Using the polarization formula (\ref{polarization-formulae}) we also have
$$
 \partial_t(\EE(Q_t)-\mathscr{P}_\infty)
 \,\leq\, (A-\mathscr{P}_\infty S)(\EE(Q_t)-\mathscr{P}_\infty)+(\EE(Q_t)-\mathscr{P}_\infty)(A-\mathscr{P}_\infty S)^{\prime}
$$
which yields the formulae 
$$
\EE\left[Q_t\right]\,\leq\, \mathscr{P}_\infty+\mathcal{E}_{t}(\mathscr{P}_\infty)\,(Q-\mathscr{P}_\infty)\,\mathcal{E}_{t}(\mathscr{P}_\infty)^{\prime}
$$

 \subsection{A Liouville Formula}
 
 This section is concerned with a stochastic version of the Liouville formula connecting the determinant with the trace of the logarithm of 
 the stochastic exponential semigroup  $\mathcal{E}_{s,t}^{\epsilon}(Q)$. This result of its own interest is also pivotal in proof of Theorem~\ref{det-E-theo}
 provided in in section~\ref{det-E-theo-proof}.
  
 We recall the trace formula
\begin{equation}\label{trace-det-formula}
\log{\mbox{\rm det}\left(\mathcal{E}_{s,t}^{\epsilon}(Q)\right)} \,=\, \tr\left(\log{\mathcal{E}_{s,t}^{\epsilon}(Q)}\right) \,=\, \int_s^t\tr(A-\phi^{\epsilon}_u(Q)S)\,du
\end{equation}
which is valid for any $\epsilon\in [0,\varepsilon_0]$.
By Jacobi's formula we have
\begin{eqnarray*}
\partial_t\mbox{\rm det}(\phi_t(Q))&=&\mbox{\rm det}(\phi_t(Q))\,\tr(Q^{-1}\partial_t\phi_t(Q))\\
&=&\mbox{\rm det}(\phi_t(Q))\,\left[2(A-\phi_t(Q)S)+\phi_t(Q)^{-1}R+\phi_t(Q)S\right]
\end{eqnarray*}
Using (\ref{trace-det-formula}), this implies that
\begin{eqnarray*}
 \log{\left[\mbox{\rm det}(\phi_t(Q)Q^{-1})\right]}&=&\int_0^t\,
\left[2~\tr(A-\phi_s(Q)S)+\tr\left(\phi_s(Q)^{-1}R+\phi_s(Q)S\right)
\right]\,ds\\
&=& \log{\left[\mbox{\rm det}(E_t(Q)E_t(Q)^{\prime})\right]}+\int_0^t\,
\tr\left(\phi_s(Q)^{-1}R+\phi_s(Q)S\right)\,ds
\end{eqnarray*}
for any $Q\in\Sa_r^+$.
In particular choosing $Q=\mathscr{P}_\infty$ we have the exponential decay
$$
\mbox{\rm det}(E_t(\mathscr{P}_\infty)) \,=\,\exp{\left[t~\tr\left(A-\mathscr{P}_{\infty}S\right)\right]}\,=\,\exp{\left[-\frac{t}{2}~
\tr\left(\mathscr{P}_{\infty}^{-1}R+\mathscr{P}_{\infty}S\right)\right]} \,\leq\, \exp{\left[-t\sqrt{\tr(RS)} \right]}
$$
To find the last inequality, we used the fact that \cite{bernstein2005matrix}
\begin{equation}\label{QQ-1}
\forall P\in\Sa^+_r,~~~\forall R,S\in\Sa^+_0,\qquad
\tr(P^{-1}R+PS) \,\geq\, 2\,\tr\left(\left[S^{1/2}~R~S^{1/2}\right]^{1/2}\right) \,\geq\, 2\sqrt{\tr(RS)} 
\end{equation}

\begin{lem}[Liouville Formula]\label{lem-det}
For any time horizon $t\geq 0$ and any $Q_0\in\Sa^+_r$ we have the log-determinant formula
\begin{equation}\label{trace-detQt}
\begin{array}{l}
\log{\left[\mbox{\rm det}(Q_tQ_0^{-1})\right]}\\
\\
\displaystyle ~~~~=~\int_0^t\,
\left[2~\tr(A-Q_sS)+\tr\left(Q^{-1}_s\left(\Sigma_{1,0}-\frac{\epsilon^2}{2}\frac{r+1}{2}~\Sigma_{\kappa,\varpi}\right)(Q_s)\right)
\right]\,ds+\epsilon\,\int_0^t~\tr\left(Q^{-1}_sdM_s\right)\\
\\
\displaystyle ~~~~\geq~ \int_0^t\,
\left[2~\tr(A-Q_sS)+\tr\left(Q^{-1}_sR^{\epsilon}+Q_sS^{\epsilon}\right)
\right]\,ds+\epsilon\,\int_0^t~\tr\left(Q^{-1}_sdM_s\right)
\end{array}\end{equation}
with the collection of matrices $(R^{\epsilon},S^{\epsilon})$ defined in (\ref{Hyp-RS-repsilon}). 
\end{lem}
The proof of this lemma is technical, and is thus given in the Appendix.

\subsection{A Dyson-Type Equation}\label{dyson-section}

We assume $\epsilon\leq \varepsilon_0$, and $\varpi=0$ and thus $(U,V)\,=\,(R,\kappa S)$ as in (\ref{cond-UV-ref-intro-uv-ref}). Now let $(\mathbf{q}_{t,i})_{1\leq i\leq r}$ be the orthonormal eigenvectors associated with the eigenvalues $0< \lambda_r(t)< \ldots<\lambda_1(t)$ of the matrix Riccati diffusion $Q_t\in\Sa_r^+$. For any $H\in\{A,R,S,U,V\}$ we set
$$
H_{t,i} \,:=\, \mathbf{q}_{t,i}^{\prime}H\,\mathbf{q}_{t,i}
$$
We then have the following general Dyson-type eigenvalue equation.

\begin{prop}
 Up to a change of probability space the eigenvalues 
 \begin{equation}\label{dyson-Ric}
\begin{array}{l}
\displaystyle
d\lambda_i(t)
=
\left[\Theta_{t,i}(\lambda_i(t))+\frac{\epsilon^2}{4}~\sum_{j\not=i}~\frac{\lambda_i(t)~\left(U_{t,j}+\lambda_j(t)^2~V_{t,j}\right)+\lambda_j(t)~\left(U_{t,i}+\lambda_i(t)^2~V_{t,i}\right)}{\lambda_i(t)-\lambda_j(t)}\right]dt\\
\\
\hskip5cm\displaystyle+
\epsilon~\sqrt{\lambda_i(t)\left(U_{t,i}+\lambda_i(t)^2~V_{t,i}\right)}~dW_{t,i}
\end{array}
\end{equation}
for some sequence $W_{t,i}$ of independent Brownian motions and the Riccati drift function
$$
\Theta_{t,i}(\lambda)=2A_{t,i}\,\lambda+R_{t,i}-\lambda^2\,S_{t,i}
$$
\end{prop}

\proof
Using the second order Hadamard variational formula we have
$$
d\lambda_i(t)\,=\,\left[\mathbf{q}_{t,i}^{\prime}\Theta(Q_t)\mathbf{q}_{t,i}+\epsilon^2\,\sum_{j\not=i}\,\frac{1}{\lambda_i(t)-\lambda_j(t)}\,\partial_t\,\langle M_{\cdot,j,i} \vert M_{\cdot,j,i}\rangle_t\right]~dt+
\epsilon~dM_{t,i,i}
$$
with the collection of martingale $$
\begin{array}{l}
\displaystyle
dM_{t,j,i} \,:=\, \mathbf{q}_{t,j}^{\prime}\,dM_t\,\mathbf{q}_{t,i}
\\
\\
\displaystyle\Longrightarrow\qquad
4\,\partial_t\,\langle M_{\cdot,j,i}  \vert M_{\cdot,j,i} \rangle_t \,=\, 1_{i=j}\,\lambda_i(t)\left(U_{t,i}+\lambda_i(t)^2\,V_{t,i}\right)\\
\\
\hskip6cm+\,
\lambda_i(t)\left(U_{t,j}+\lambda_j(t)^2\,V_{t,j}\right)+\lambda_j(t)\left(U_{t,i}+\lambda_i(t)^2\,V_{t,i}\right)
\end{array}$$
Also observe that for any $i\not=j$ we have $\partial_t\langle M_{\cdot,i,i}\vert M_{\cdot,j,j}\rangle_t:=0$.This yields the formula (\ref{dyson-Ric}).
\qed

We consider the diffusion function (\ref{cond-UV-ref-intro-uv-ref}) and we assume that
$$
(A,R,S,U,V)=(\mathfrak{a}\,I,\mathfrak{r}\,I,\mathfrak{s}\,I, \mathfrak{u}\,I, \mathfrak{v}\,I)\quad \mbox{\rm for some}\quad \mathfrak{a}\in\RR\qquad
\mathfrak{r},\mathfrak{s}\in]0,\infty[\quad\mbox{\rm and}\quad \mathfrak{u},\mathfrak{v}\geq 0
$$ 
In this special case, the eigenvalues $0< \lambda_r(t)< \ldots<\lambda_1(t)$ of the matrix Riccati diffusion $Q_t\in\Sa_r^+$ satisfy the Dyson-type diffusion equation 
\begin{equation}\label{dyson-mod-intro}
\begin{array}{l}
\displaystyle
d\lambda_i(t)
=
\left[\Theta(\lambda_i(t))+\frac{\epsilon^2}{4}~\sum_{j\not=i}~\frac{\lambda_i(t)\Sigma_{\kappa,0}(\lambda_j(t))+\lambda_j(t)\Sigma_{\kappa,0}(\lambda_i(t))}{\lambda_i(t)-\lambda_j(t)}\right]dt+
\epsilon\,\sqrt{\lambda_i(t)}\,\Sigma_{\kappa,0}^{1/2}(\lambda_i(t))\,dW^{i}_t
\end{array}
\end{equation}
with the (re-defined here) one-dimensional Riccati drift and diffusion functions
$$
\Theta(\lambda) \,:=\, 2\mathfrak{a}\lambda+\mathfrak{r}-\lambda^2\mathfrak{s}\qquad \mbox{\rm and }\qquad\Sigma_{\kappa,0}(\lambda)\,:=\,\mathfrak{u}+\lambda^2~\mathfrak{v} 
$$

When $\epsilon=0$ the equation (\ref{dyson-mod-intro}) resumes to a univariate Riccati equation; that is we have that $\lambda_i(t)=\lambda(t)$ for any $1\leq i\leq r$. In this situation it is rather well known that for any $t\geq \upsilon>0$
$$
 \left\vert\lambda(t)-\frac{\mathfrak{a}+\sqrt{\mathfrak{a}^2+\mathfrak{r}\mathfrak{s}}}{\mathfrak{s}}\right\vert \,\leq\, c_{\upsilon}\,\exp{\left(-2t\sqrt{\mathfrak{a}^2+\mathfrak{r}\mathfrak{s}}\right)}\quad \mbox{\rm for some finite constant}~c_{\upsilon}<\infty
$$
A proof of the above assertion can be found for instance in~\cite{2017arXiv171110065B}. Clearly, the very special case in (\ref{dyson-mod-intro-ref}) follows from the above.

\addtocontents{toc}{\protect\setcounter{tocdepth}{1}}
\section{Proofs of the Main Theorems}\label{proof-main-theo}

\subsection{Proof of Theorem~\ref{theo-existence-s-ric}}\label{theo-existence-s-ric-proof}

The proof of the first assertion follows the arguments provided in Section 3 of~\cite{Khasminskii}. More precisely, consider the sets
$$
	\Omega_n:=\{P\in \Sa^0_r~:~ \tr(P)\leq n\} 
$$
and the exit time
$$
	\tau_n:=\inf{\left\{t\geq 0~:~Q_t\not\in \Omega_n\right\}}
$$

Up to a change of probability space the process $Q_{t\wedge \tau_n}$ when $\epsilon=2/\sqrt{N}$ coincides with the evolution of sample covariance matrices of an associated system of particles interacting with their internal sample covariance matrices; see~\cite{Bishop/DelMoral/Niclas:2017,DelMoral/Tugaut:2016} and Section~\ref{EnKF-sec} in the present article. Notice that this system of interacting diffusions is well defined on $[0,\tau_n]$. Up to a time-rescaling of the Brownian motions in (\ref{f21}), this result is also met for any $\epsilon\geq 0$, so that $Q_{t\wedge \tau_n}$ cannot exit the set $\Sa^0_r$.
 
For any $m\geq n$ we have 
$$
	Q_{t\wedge \tau_m} \,:=\, Q_{t\wedge \tau_n} \,=\, Q_t\quad \mbox{\rm for any} \quad t\in [0,\tau_n]
$$
Let $\tau^\star$ be the finite or infinite limit of the monotone increasing sequence $\tau_n$. The stochastic process,
$$
\mathbf{Q}_t \,=\, Q_{t\wedge \tau_n}\,1_{[0,\tau_n[}(t)
$$
is a well-defined Markov process for any $t<\tau_n$. Finally, observe that
$$
\partial_t\EE(\tr(Q_{t\wedge \tau_n})) \,\leq\, 2\mathfrak{a}\, \EE(\tr(Q_{t\wedge \tau_n}))+\mathfrak{r}-\mathfrak{s}\,\EE(\tr(Q_{t\wedge \tau_n}))^2 \,\leq\, 2\mathfrak{a}\, \EE(\tr(Q_{t\wedge \tau_n}))+\mathfrak{r}
$$
with the parameters 
$$
(\mathfrak{a},\mathfrak{r},\mathfrak{s})\,:=\,(\mu(A),\, \tr(R),\, r^{-1}\lambda_{r}(S))
$$
This implies that
$$
n\,\PP\left(\tau_n\leq t\right)\,\leq\, \EE(\tr(\mathbf{Q}_{t\wedge \tau_n})) \,=\, \EE(\tr(\mathbf{Q}_{t})1_{t<\tau_n})+\EE(\tr(\mathbf{Q}_{\tau_n})1_{\tau_n\leq t}) \,\leq\, e^{2\mathfrak{a}t}\,(\tr(Q_{0})+\mathfrak{r}/(2\mathfrak{a}))
$$ 
from which we check that
$$
\PP\left(\tau_n\leq t\right) \,\leq\, \frac{1}{n}\left(e^{2\mathfrak{a}t}\,(\tr(Q_{0})+\mathfrak{r}/(2\mathfrak{a}))\right) \quad\Longrightarrow\quad \PP(\tau^\star=\infty)=1
$$
We conclude that (\ref{f21}) has an unique weak solution.

The proof of second assertion is a consequence of  Lemma~\ref{lem-det} combined with the McKean argument developed  in Proposition 4.3 in~\cite{mayerhofer2011strong}. To check this claim, we set
$$
Z_t(Q):=\mbox{\rm det}\left(\mathcal{E}_{t}^{\epsilon}(Q)^{-1}\phi^{\epsilon}_t(Q)\left(\mathcal{E}_{t}^{\epsilon}(Q)^{\prime}\right)^{-1}\right)
\quad\mbox{\rm and}\quad
\tau_{ Q}^{\epsilon}:=\inf{\left\{t> 0~:~\phi^{\epsilon}_t(Q)\in \partial \Sa_r^+\right\}}
$$
Using (\ref{trace-det-formula}) we have
$$
\mbox{\rm det}\left(\mathcal{E}_{t}^{\epsilon}(Q)^{-1}\left(\mathcal{E}_{t}^{\epsilon}(Q)^{\prime}\right)^{-1}\right) \,=\, \exp{\left[
-2\int_0^t\tr(A-\phi^{\epsilon}_s(Q)S)~ds\right]}
$$
By Lemma~\ref{lem-det} we have the decomposition
\begin{eqnarray*}
\log{Z_t(Q)} &=&  \log{Z_0(Q)}+\mathfrak{m}_t(Q) + \int_0^t
\tr\left(Q^{-1}_s\left(\Sigma_{1,0}-\frac{\epsilon^2}{2}\frac{r+1}{2}\,\Sigma_{\kappa,\varpi}\right)(Q_s)\right)\,ds \\ 
&\geq& \log{Z_0(Q)}+\mathfrak{m}_t(Q)
\end{eqnarray*}
with the continuous local martingale $\mathfrak{m}_t(Q)$ on $[0,\tau^{\epsilon}_{ Q}[$ defined by
$$
\mathfrak{m}_t(Q) \,:=\, \epsilon\,\int_0^t\tr\left(Q^{-1}_s\,dM_s\right)
$$
and the noting that the following positive mapping satisfies
\begin{eqnarray*}
\int_0^t
\tr\left(Q^{-1}_s\left(\Sigma_{1,0}-\frac{\epsilon^2}{2}\frac{r+1}{2}\,\Sigma_{\kappa,\varpi}\right)(Q_s)\right)\,ds
&\geq &\displaystyle \int_0^t\left[\tr\left(Q^{-1}_sR^{\epsilon}\right)+\tr\left(Q_sS^{\epsilon}\right)\right]\,ds \\
&\geq& 2t\,\sqrt{\tr\left(R^{\epsilon}S^{\epsilon}\right)}
\end{eqnarray*}
The end of the proof of the existence and uniqueness of a strong solution of the matrix Riccati diffusion (\ref{f21}) on $\Sa^+_r$ is now a consequence of Proposition 4.3 in~\cite{mayerhofer2011strong}. Specifically, if $\tau_{ Q}^{\epsilon}<\infty$ on some event  with positive probability, then on this event set we have
$$
\lim_{t\rightarrow \tau_{ Q}^{\epsilon}}\log{Z_t(Q)}=-\infty \qquad\Longrightarrow\qquad
\lim_{t\rightarrow \tau_{ Q}^{\epsilon}}\mathfrak{m}_t(Q)=-\infty
$$
This contradicts the fact that either $\lim_{t\rightarrow \tau_{ Q}^{\epsilon}}\mathfrak{m}_t(Q)\in \RR$ or 
$$
\limsup_{t\rightarrow \tau_{ Q}^{\epsilon}}\mathfrak{m}_t(Q) \,=\, +\infty \,=\, -\liminf_{t\rightarrow \tau_{ Q}^{\epsilon}}\mathfrak{m}_t(Q)
$$
This ends the proof of the second assertion.

Now we come to the proof of (\ref{equivalence-intro}). We set
$$
d\widetilde{\Wa}_t \,:=\, (Q_t\frownotimes \Sigma_{\kappa,\varpi}(Q_t))^{1/2}\,d\Wa_t
$$
Thus, the angle bracket of the matrix-valued martingale $\widetilde{\Wa}_t$ is given by
$$
\begin{array}{l}
\partial_t\langle\, \widetilde{\Wa}(i,j)~\vert~\widetilde{\Wa}(k,l) \,\rangle_t\\
\\
\displaystyle\qquad=~\sum_{1\leq i^{\prime},k^{\prime}\leq r}(Q_t\frownotimes \Sigma_{\kappa,\varpi}(Q_t))^{1/2}((i,j),(i^{\prime},j^{\prime}))~1_{(i^{\prime},j^{\prime})=(k^{\prime},l^{\prime})}~
(Q_t\frownotimes \Sigma_{\kappa,\varpi}(Q_t))^{1/2}((k,l),(k^{\prime},l^{\prime}))\\
\\
\displaystyle\qquad=~(Q_t\frownotimes \Sigma_{\kappa,\varpi}(Q_t))((i,j),(k,l))
\end{array}$$
Using (\ref{ref-M-frownotimes}) we conclude that
$$
dQ_t \,\stackrel{ law}{=}\, \Theta(Q_t)\,dt+\epsilon\,(Q_t\frownotimes \Sigma_{\kappa,\varpi}(Q_t))^{1/2}~d\Wa_t
$$
For any matrix $H\in \Ma_r$ we have
$$
 (P_1\,\overline{\otimes}\,P_2)(H) \,=\, P_1\,H^{\prime}\,P_2 \,=\, P_2\,H \, P_1 \,=\, (P_2\otimes P_1)(H) \,=\, (P_1\otimes P_2)(H^{\prime}) \,=\, (P_2\,\overline{\otimes}\,P_1)(H^{\prime})
$$
Also observe that 
$$
\begin{array}{l}
(P_1\,\otimes\,P_2)(H) \,=\, P_1\, H \, P_2 \,=\, P_2\,H^{\prime}\,P_1 \,=\, (P_2~\otimes~ P_1)(H^{\prime})\\
\\
\qquad\qquad\displaystyle\Longrightarrow\qquad (P_1\frownotimes P_2)(H) \,=\, \frac{1}{2}\left[(P_1\otimes P_2)+(P_1\,\overline{\otimes}\,P_2)\right]\left(\frac{H+H^{\prime}}{2}\right)\\
\\ \qquad\qquad\displaystyle\Longrightarrow\qquad (P_1\frownotimes P_2)(H) \,=\, \frac{1}{4}\left[(P_1\otimes P_2)+(P_1\,\overline{\otimes}\,P_2)\right]\left[{(I\,{\otimes}\,I)+(I\,\overline{\otimes}\,I)}\right](H)
\end{array}
$$
This shows that
$$
H^{\prime}=-H \quad\Longrightarrow\quad (P_1\frownotimes P_2)(H)=0 \quad\Longrightarrow\quad (P_1\frownotimes P_2)^{1/2}(H)=0
$$
Additionally, we have
$$
H=H^{\prime} \quad\Longrightarrow\quad (P_1\frownotimes P_2)(H)=(P_1\otimes_{ s} P_2)(H)
$$
By Doob's representation theorem (see Theorem 4.2~\cite{karatzas}, and the original work of Doob~\cite{doob}), the proof of (\ref{equivalence-intro}) is now a consequence of the fact that
$$
\frac{\Wa_t+\Wa^{\prime}_t}{2}  ~\stackrel{ law}{=}~ \Va_{t,\mathrm{sym}}
$$
The proof of (\ref{equivalence-intro-2}) comes from the fact that
$
 v_t=\varsigma (\Va_{t,\mathrm{sym}})$ is an $\overline{r}$-dimensional Brownian motion, and we have that
$$
dq_t \,=\, \theta(q_t)\,dt+\epsilon\,\sum_{1\leq i\leq \overline{r}}\,\sigma_i(q_t)\,dv^i_t
$$
In Stratonovitch form we have
$$
d q_t=\theta_{\epsilon}(q_t)~d t+\epsilon~\sum_{1\leq i\leq \overline{r}}~\sigma_i(q_t)\circ dv^i_t
\quad\mbox{\rm 
with the drift}\quad
\theta_{\epsilon}^i=\theta^{i}-\frac{\epsilon^2}{2}\sum_{1\leq k,l\leq \overline{r}}~\sigma_l^k~\partial_{q_k}\sigma_l^i
$$
The notation $\sigma_i(q_t)\,\circ\, dv^i_t$ implies that It\^o integrals are replaced by Stratonovitch integrals. We also recall that
$$
\epsilon\leq \varepsilon_0 \quad\Longrightarrow\quad \forall t>0,~~ Q_t\in\Sa_r^+ \quad\Longrightarrow\quad \forall t>0,~~ q_t\in \Da_{\overline{r}}:=\varsigma(\Sa_r^+)
$$
This shows that for any $t>0$ the process $q_t$ never visits the boundary $\partial \Da_{\overline{r}}=\varsigma(\partial \Sa_r^+)$, even when we start at some state $q_0\in \partial \Da_{\overline{r}}$. On the other hand, we have
$$
(\ref{eq-ap-2}) \quad\Longrightarrow\quad \forall q\in \Da_{\overline{r}},~~ \left\{\varsigma^{-1}(q)\otimes_s \Sigma_{\kappa,\varpi}\left(\varsigma^{-1}(q)\right)\right\}\,>\,0 
$$
This shows that the linear span of the $\overline{r}$-vector fields $q\in \Da_{\overline{r}}\mapsto\sigma_i(q)\in \RR^{\overline{r}}$ of the diffusion is all $\RR^{\overline{r}}$.  Also notice that the set of point $q\in \varsigma(\Sa^0_r)$ for which $\det(\left\{\varsigma^{-1}(q)\otimes_s \Sigma_{\kappa,\varpi}\varsigma^{-1}(q))\right\})=0$ coincides with $\partial \Da_{\overline{r}}$ which is of null measure in $\varsigma(\Sa_r)$. In other words the elliptic degeneracies of the diffusion $q_t$ are of null Lebesgue measure.

The generator of the diffusion  $q_t\in\Da_{\overline{r}}$ can be expressed  in H\"ormander form  by the formula
\begin{eqnarray*}
L&=&\mathfrak{X}_{\epsilon,0}+\frac{1}{2}~\sum_{1\leq i\leq \overline{r}}~\mathfrak{X}_{\epsilon,i}^2
\end{eqnarray*}
with the first order $C^{\infty}$-vector fields on $\Da_{\overline{r}}$ given by
$$
\mathfrak{X}_{\epsilon,0}:=\sum_{1\leq i\leq \overline{r}}~\theta_{\epsilon}^i~\partial_{q_i}\quad\mbox{\rm and}\quad \mathfrak{X}_{\epsilon,i}:=\epsilon~\sum_{1\leq k\leq \overline{r}}~\sigma_i^k~\partial_{q_k}
$$
The operator $L$ is hypo-elliptic, since the Lie algebra generated by the $\overline{r}$ vector fields $(\mathfrak{X}_{\epsilon,i})_{1\leq i\leq  \overline{r}}$ span the entire Euclidian space $\RR^{ \overline{r}}$ at any state $q\in \Da_{\overline{r}}$. By H\"ormander's theorem, it follows that the transition semigroup $\pi^{\epsilon}_t(p,dq)$ of $q_t$ has smooth positive densities $\rho^{\epsilon}\in C^{\infty}(]0,\infty[\times\Da^2_{\overline{r}})$; see e.g.~\cite{arnaudon,bell}, and the reference by Bramanti~\cite{Bramanti} dedicated to hypo-elliptic operators and H\"ormander vector fields.

This ends the proof of the theorem.\qed

\subsection{Proof of Theorem~\ref{theo-existence-s-ric-proof-bis}}\label{theo-existence-s-ric-bis-proof}

Let $(\mathfrak{a},\mathfrak{r},\mathfrak{s}):=(\mu(A), \tr(R), r^{-1}\lambda_r(S))$. Also define the collection of parameters
$$
\mathfrak{r}^\epsilon_{n}\,:= \,\mathfrak{r}+\frac{\epsilon^2}{2}~(n-1)~\lambda_1(U)\quad \mbox{\rm and}\quad
\mathfrak{s}^\epsilon_{n}\,:=\,\mathfrak{s}-\frac{\epsilon^2}{2}~(n-1)~\lambda_1(V)
$$
Observe that
$$
\kappa=0\quad\Longrightarrow\quad  \mathfrak{s}^\epsilon_{n}\,:=\,\mathfrak{s}
$$
For any $n\geq 1$ we let $\varepsilon_n$ be the largest (finite) parameter $\epsilon\geq0$ such that $\mathfrak{s}^\epsilon_{n}>0$ and we set
$$
(\mathfrak{r}_n,\mathfrak{s}_n):=(\mathfrak{r}^{\varepsilon_n}_{n},\mathfrak{s}^{\varepsilon_n}_{n})
$$
Let $\mathfrak{p}_{t,n}$ be the one-dimensional Riccati flow associated with the differential equation
$$
\partial_t\mathfrak{p}_{t,n}\,=\,2\mathfrak{a}\mathfrak{p}_{t,n}+\mathfrak{r}_{n}-\mathfrak{s}_{n}\,\mathfrak{p}_{t,n}^2 
$$
and let $\mathfrak{p}_{0,n}=\tr(Q)$. In this notation, for any $n\geq 1$ and any $\epsilon\in [0,\varepsilon_n]$ we have the estimate
\begin{equation}\label{trace-estimates}
\EE[\tr(\phi^{\epsilon}_t(Q))^n]^{1/n} \,\leq\, \mathfrak{p}_{t,n} \,\leq\, \mathfrak{p}_{\infty,n}\vee \tr(Q)
\quad\mbox{
with}
\quad
\mathfrak{p}_{\infty,n}\,:=\,\frac{\mathfrak{a}+\sqrt{\mathfrak{a}^2+\mathfrak{r}_{n}\mathfrak{s}_{n}}}{\mathfrak{s}_{n}}
\end{equation}
Observe that
$$
\kappa=0\quad\Longrightarrow\quad \mathfrak{p}_{\infty,n}\,:=\,\frac{\mathfrak{a}+\sqrt{\mathfrak{a}^2+\mathfrak{r}\mathfrak{s}+\frac{\epsilon^2}{2}\,(n-1)\,\mathfrak{r}\,\lambda_1(U)}}{\mathfrak{s}} \,\leq\,
\mathfrak{p}_{\infty,1}+\frac{\epsilon}{\mathfrak{s}}\,\sqrt{n-1}\,\sqrt{\frac{\mathfrak{r}\,\lambda_1(U)}{2}}
$$
To check (\ref{trace-estimates}), observe that
$$
\mathfrak{q}_t:=\tr(Q_t) ~~\Rightarrow~~
\tr(\Theta(Q_t)) \,\leq\, 2\mathfrak{a}\mathfrak{q}_t+\mathfrak{r}-\mathfrak{s}\mathfrak{q}_t^2\quad \mbox{\rm and}\quad
\tr(Q_t\Sigma_{\kappa,\varpi}(Q_t)) \,\leq\, \mathfrak{q}_t\lambda_1(U)+\mathfrak{q}_t^3\lambda_1(V)
$$
This yields the formula,
$$
d\mathfrak{q}_t^n\,=\, n\left[\mathfrak{q}^{n-1}_t\,\tr(\Theta(Q_t))+~\frac{\epsilon^2}{2}\,(n-1)\,\mathfrak{q}^{n-2}_t\,\tr(Q_t\Sigma_{\kappa,\varpi}(Q_t))\right]dt+\epsilon\,n\,\mathfrak{q}^{n-1}_t\,d\tr(M_t)
$$
from which we check the differential inequalities
\begin{eqnarray*}
n^{-1}\,\partial_t\EE(\mathfrak{q}_t^n)&\leq &2\mathfrak{a}~\EE(\mathfrak{q}_t^n)+\mathfrak{r}^\epsilon_{n}~\EE(\mathfrak{q}^{n-1}_t)-\mathfrak{s}^\epsilon_{n}~~\EE(\mathfrak{q}^{n+1}_t)\\
&\leq &2\mathfrak{a}~\EE(\mathfrak{q}_t^n)+\mathfrak{r}^\epsilon_{n}~\EE(\mathfrak{q}^{n}_t)^{1-1/n}-\mathfrak{s}^\epsilon_{n}~\EE(\mathfrak{q}^{n}_t)^{1+1/n}
\end{eqnarray*}
The last lines follows from the fact that,
$$
\EE(\mathfrak{q}^{n-1}_t) \,\leq\, \EE(\mathfrak{q}^{n}_t)^{1-1/n}\quad \mbox{\rm and}\quad
\EE(\mathfrak{q}^{n+1}_t) \,\leq\, \EE(\mathfrak{q}^{n}_t)^{1+1/n}
$$
We conclude that
$$
\partial_t\EE(\mathfrak{q}_t^n)^{1/n}\,=\,n^{-1}~\EE(\mathfrak{q}_t^n)^{-(1-1/n)}~
\partial_t\EE(\mathfrak{q}_t^n) ~\leq~ 2\mathfrak{a}~\EE(\mathfrak{q}_t^n)^{1/n}+\mathfrak{r}^\epsilon_{n}-\mathfrak{s}^\epsilon_{n}~\EE(\mathfrak{q}^{n}_t)^{2/n}
$$
Now (\ref{trace-estimates}) is a direct consequence of Lemma~\ref{Lemma-1} and the estimates on one-dimensional Riccati flows presented in~\cite{2017arXiv171110065B}. 

In addition, using the uniform estimates presented in~\cite{2017arXiv171110065B}  for any $t\geq \upsilon>0$ we have
$$
\EE(\tr(\phi^{\epsilon}_t(Q))^n)^{1/n} \,\leq\, \mathfrak{p}_{t,n} \,\leq\, c_{\upsilon}\,\mathfrak{p}_{\infty,n}^{\star}
\quad\mbox{
with}
\quad
\mathfrak{p}_{\infty,n}^{\star}\,:=\,\frac{\mathfrak{a}+3\sqrt{\mathfrak{a}^2+\mathfrak{r}_{n}\mathfrak{s}_{n}}}{\mathfrak{s}_{n}}
$$
Observe again that
$$
\kappa=0 \quad\Longrightarrow\quad \mathfrak{p}_{\infty,n}^{\star}\,\leq \,
\mathfrak{p}_{\infty,1}^{\star}+\frac{3\,\epsilon}{\mathfrak{s}}\,\sqrt{n-1}\,\sqrt{\frac{\mathfrak{r}\,\lambda_1(U)}{2}}
$$
This completes the proof of the Riccati diffusion moment estimates in (\ref{trace-Phi-inverse}), (\ref{trace-Phi-inverse-bis}) and (\ref{ref-trace-unif-intro}).

Now we come to the proof of the trace estimates of the inverse stochastic flow $\phi^{-\epsilon}_t(Q)$ stated in (\ref{trace-Phi-inverse}) and (\ref{trace-Phi-inverse-bis}). The approach follows the preceding discussion but is more notationally and computationally burdensome, given the form of the inverse flow; e.g. see (\ref{eq-inverse}). We set
$$
\mathfrak{a}_{-}:=-\lambda_r(A_{\mathrm{sym}})\qquad \mathfrak{r}_-:=\tr\left(S\right)\qquad \mathfrak{s}_-:=r^{-1}\lambda_r\left(R\right)
$$
Note the exchanged roles of $R$ and $S$ in $\mathfrak{s}_-$ and $\mathfrak{r}_-$.

For any $n\geq 1$ we let $\varepsilon_{n,-}$ be the largest parameter $\epsilon\geq0$ such that
$$
	\mathfrak{s}^\epsilon_{n,-} \,:=\, \mathfrak{s}_- -\frac{\epsilon^2}{2}\,\left[ (n+r^{-1})\,\lambda_1(U)+\,\frac{\lambda_1(V)}{4}\right] \,>\,0
$$
Also consider the collection of parameters
$$
\mathfrak{r}^\epsilon_{n,-}(Q) \,=\,\mathfrak{r}_-
+\frac{\epsilon^2}{2}\,\left[
\left(1+\frac{r}{2}\right)\,\tr(V)+(n-1)\,\lambda_1(V)+~\frac{\lambda_1(V)}{4}\,\left(\mathfrak{p}_{\infty,2n}^2 \vee \tr(Q)^2\right)\right]
$$
with the sequence of parameters $\mathfrak{p}_{\infty,n}$ introduced in (\ref{trace-estimates}).

In this notation, for any $n\geq 1$, $Q\in\Sa_r^+$ and any $\epsilon\in [0,\varepsilon_{n,-}]$ we have the uniform estimate
\begin{equation}\label{trace-estimates-inverse}
\sup_{t\geq 0}\,{\EE(\tr(\phi^{-\epsilon}_t(Q))^n)^{1/n}} \,\leq\, \mathfrak{p}_{\infty,n,-}(Q)\vee \tr(Q^{-1})
\end{equation}
with the collection of parameters
$$
\mathfrak{p}_{\infty,n,-}(Q) \,:=\,\frac{\mathfrak{a}_-+\sqrt{\mathfrak{a}_-^2+\mathfrak{r}_{n,-}(Q)\mathfrak{s}_{n,-}}}{\mathfrak{s}_{n,-}}\qquad \mbox{\rm with}\qquad (\mathfrak{r}_{n,-}(Q),\mathfrak{s}_{n,-}):=(\mathfrak{r}^{\varepsilon_{n,-}}_{n,-}(Q),\mathfrak{s}^{\varepsilon_{n,-}}_{n,-})
$$
To check this claim, observe that
$$
\mathfrak{q}_{t,-} \,:=\, \tr(Q^{-1}_t) \quad\Longrightarrow\quad \tr\left(
\Theta^\epsilon_{-}(Q^{-1}_t)\right) \,\leq\, 2\mathfrak{a}_-\mathfrak{q}_{t,-}   +\mathfrak{r}^\epsilon_{t,-}- \mathfrak{s}^\epsilon_{1,-}\,\mathfrak{q}_{t,-}^2
$$
with the functions
$$
 \mathfrak{r}^\epsilon_{t,-}:=\mathfrak{r}_-+\frac{\epsilon^2}{4}\,\left[
\left(2+r\right)\,\tr\left(V\right)+\frac{\lambda_1(V)}{2}\,\tr\left(Q_t\right)^2\right]
$$
In the last display we have used the fact that
$$
\tr\left(Q^{-2}_t\right)\,\geq\, r^{-1} (\tr\left(Q^{-1}_t\right))^2\quad\mbox{\rm and}\quad
\tr\left(VQ_t\right)\tr\left(Q^{-1}_t\right) \,\leq\, \frac{\lambda_1(V)}{2}\,\left(\tr\left(Q_t\right)^2+\tr\left(Q_t^{-1}\right)^2\right)
$$
On the other hand, we have
$$
\displaystyle d\mathfrak{q}_{t,-} \,=\, \tr\left(
\Theta^\epsilon_{-}(Q^{-1}_t)\right)dt+\epsilon\,d\mathfrak{m}_{t,-}
$$
with
$$
\partial_t\langle \mathfrak{m}_{\cdot,-}~\vert~\mathfrak{m}_{\cdot,-}\rangle_t\,=\,\tr\left(Q^{-1}_t\Sigma_{\kappa,\varpi,-}\left(Q^{-1}_t\right)\right) \,\leq\, \mathfrak{q}_{t,-}\,\lambda_1(V)+\mathfrak{q}_{t,-}^3\,\lambda_1(U)
$$
Thus, for any $n\geq 1$ we have
\begin{equation*}
\begin{array}{l}
\displaystyle n^{-1}\,\partial_t\,\EE\left(\mathfrak{q}_{t,-}^n\right) \,\leq\, 2 \mathfrak{a}_{-}\,\EE\left(\mathfrak{q}_{t,-}^n\right) + \EE\left(\mathfrak{r}^\epsilon_{t,n,-}\,\mathfrak{q}_{t,-}^{n-1}\right)- \mathfrak{s}^\epsilon_{n,-}\,\EE\left(\mathfrak{q}_{t,-}^{n+1}\right)
\end{array}
\end{equation*}
and the collection of stochastic processes
\begin{eqnarray*}
\mathfrak{r}^\epsilon_{t,n,-} &:=& \mathfrak{r}^\epsilon_{t,-}+\frac{\epsilon^2}{2}\,(n-1)\,\lambda_1(V) \\
 &=& \tr\left(S\right)
+\frac{\epsilon^2}{2}~\left[
\left(1+\frac{r}{2}\right)\,\tr\left(V\right)+(n-1)\,\lambda_1(V)+\frac{\lambda_1(V)}{4}\,\tr\left(Q_t\right)^2\right]
\end{eqnarray*}
On the other hand, using (\ref{trace-estimates}) we check that
$$
\EE\left(\mathfrak{r}^\epsilon_{t,n,-}\,\mathfrak{q}_{t,-}^{n-1}\right) \,\leq\, \EE\left(\mathfrak{q}_{t,-}^{n}\right)^{1-1/n}\,\mathfrak{r}^\epsilon_{\star,n,-}\qquad\mbox{\rm
with}\qquad
\mathfrak{r}^\epsilon_{\star,n,-} \,:=\, \sup_{t\geq 0}{\EE((\mathfrak{r}^\epsilon_{t,n,-})^{\,n})^{1/n}}\,\leq\,
\mathfrak{r}^\epsilon_{n,-}(Q) 
$$
This yields the estimate
$$
\displaystyle n^{-1}\partial_t\EE\left(\mathfrak{q}_{t,-}^n\right) \,=\, 2\mathfrak{a}_{-}\,\EE\left(\mathfrak{q}_{t,-}^{n} \right) +\mathfrak{r}^\epsilon_{\star,n,-} \,\EE\left(\mathfrak{q}_{t,-}^{n}\right)^{1-1/n}- \mathfrak{s}^\epsilon_{n,-}\,\EE\left(\mathfrak{q}_{t,-}^{n}\right)^{1+1/n}
$$
The end of the proof of (\ref{trace-estimates-inverse}) now follows the same lines of arguments as the proof of the trace estimates (\ref{trace-estimates}), thus it is skipped. This ends the proof of (\ref{trace-estimates-inverse}) and thus the proof of the inverse Riccati diffusion moment estimates in (\ref{trace-Phi-inverse}). The uniform estimates on the inverse flow in (\ref{trace-Phi-inverse-bis}) follow the same line of arguments as in the proof of the l.h.s in (\ref{trace-Phi-inverse-bis}) given the inverse moment estimates already proved in (\ref{trace-Phi-inverse}).

This ends the proof of the theorem.\qed

\subsection{Proof of Theorem~\ref{theo-3-intro}}\label{proof-theo-3-intro}

Consider the Gramian matrix,
$$
\GG_t(Q)=\int_0^t\,\mathcal{E}_{s}(Q)^{\prime}S~\mathcal{E}_{s}(Q)~ds
$$
and the non-negative matrix function
$$
G_t(Q)\,=\,\mathcal{E}_{t}(Q)\left[\left(Q\,\GG_{t}(Q)\Sigma_{\kappa,\varpi}(Q)\right)_{\mathrm{sym}}+\tfrac{1}{2}\left[
Q\,\tr\left(\Sigma_{\kappa,\varpi}(Q)\,\GG_{t}(Q)\right)+\Sigma_{\kappa,\varpi}(Q)\,\tr\left(Q\,\GG_{t}(Q)\right)\right]\right]\mathcal{E}_{t}(Q)^{\prime}
$$
In this notation we have the second order decomposition 
$$
\phi^{\epsilon}_t(Q)=\phi_{t}(Q)+\epsilon~\MM^{\epsilon}_{t}(Q)-\frac{\epsilon^2}{2}~\BB_{t}^{\epsilon}(Q)
$$
with the processes
$$
\MM^{\epsilon}_{t}(Q) \,:=\, \int_0^t\,\Ea_{t-u}(\phi^{\epsilon}_u(Q))~ dM_{u}(Q)~\Ea_{t-u}(\phi^{\epsilon}_u(Q))^{\prime}\quad
\mbox{\rm and}\quad
\BB_{t}^{\epsilon}(Q)\,:=\,\int_0^t\, G_{t-u}\left(\phi^{\epsilon}_u(Q)\right)\,du
$$
A proof of the above decomposition can be found in~\cite[in the proof of Theorem 1.3]{Bishop/DelMoral/Niclas:2017}. The above forward-backward perturbation formula can be thought as an extended version of the Alekseev-Gr\"obner lemma to diffusion flows in matrix spaces~\cite{singh}.

Using (\ref{ref-E-1}) we have
$$
\Vert \GG_t(Q)\Vert \,\leq\, c\,(1+\Vert Q\Vert^2) \quad\Longrightarrow\quad
 \Vert  G_t\left(Q\right)\Vert \,\leq\, c\,(1+\Vert Q\Vert^5)\,(\lambda_1(U)+\lambda_1(V)\,\Vert Q\Vert^2)\,\exp{\left(-2\beta t\right)}
$$
Using the generalized Minkowski inequality we check the estimate
$$
\vertiii{\BB_{t}^{\epsilon}(Q)}_n \,\leq\, c\,\int_0^t \left[1+\vertiii{\phi^{\epsilon}_u(Q)}_{10n}^5\right]\left[\lambda_1(U)+\lambda_1(V)\,\vertiii{\phi^{\epsilon}_u(Q)}_{4n}^2\right]\,\exp{\left(-2\beta (t-u)\right)}\,du
$$
By (\ref{trace-Phi-inverse}) for any $\epsilon\leq \varepsilon_{10n}(V)$ we have
$$
\vertiii{\BB_{t}^{\epsilon}(Q)}_n \,\leq\, c_n\, (1+\Vert Q\Vert^5)\,(\lambda_1(U)+\lambda_1(V)\,\Vert Q\Vert^2)
$$
This yields the uniform bias estimate
$$
0~\leq\, \phi_{t}\left(Q\right)- \EE\left[\phi^{\epsilon}_t(Q)\right] \,\leq\, c\,\epsilon^2\,
(1+\Vert Q\Vert^5)\,(\lambda_1(U)+\lambda_1(V)\,\Vert Q\Vert^2)\, I \quad\Longrightarrow\quad (\ref{ref-vp-max})
$$
In addition, using (\ref{ref-trace-unif-intro}) when $\kappa=0$ we have for any $\epsilon\geq 0$,
$$
\vertiii{\BB_{t}^{\epsilon}(Q)}_n \,\leq\, c\,(1+ \Vert Q\Vert^5)\,(1+\epsilon~\sqrt{n})^5
$$
which yields
$$
0~\leq\, \phi_{t}\left(Q\right)- \EE\left[\phi^{\epsilon}_t(Q)\right]\,\leq\, c\,\epsilon^2\,(1+ \Vert Q\Vert^5)\,(1+\epsilon~\sqrt{n})^5\, I
$$
and we may refine (\ref{ref-vp-max}) appropriately in this case, $\kappa=0$, for any $\epsilon\geq 0$.

The trace of the martingale
$$
\MM^{\epsilon}_{s,t}(Q) \,:=\, \int_0^s\,\Ea_{t-u}(\phi^{\epsilon}_u(Q))\, dM_{u}(Q)\,\Ea_{t-u}(\phi^{\epsilon}_u(Q))^{\prime}
$$
is a martingale with angle bracket
$$
\begin{array}{l}
\displaystyle
4^{-1} \partial_s\,\langle  \tr\left(\MM^{\epsilon}_{{\,\cdot},t}(Q)\right)~\vert~\tr\left(\MM^{\epsilon}_{{\,\cdot},t}(Q)\right)\rangle_s\\
\\
\qquad\qquad\displaystyle= ~\tr\left[
\phi^{\epsilon}_s(Q)\Ea_{t-s}(\phi^{\epsilon}_s(Q))^{\prime}\,
\Ea_{t-s}(\phi^{\epsilon}_s(Q))\,\Sigma_{\kappa,\varpi}\left(\phi^{\epsilon}_s(Q)\right)\,\Ea_{t-s}(\phi^{\epsilon}_s(Q))^{\prime}\,\Ea_{t-s}(\phi^{\epsilon}_s(Q))
\right]\\
\\
\qquad\qquad\leq~ \Vert \Ea_{t-s}(\phi^{\epsilon}_s(Q))\Vert_{\mathrm{Frob}}^4\,\tr(\phi^{\epsilon}_s(Q))~\tr\left(\Sigma_{\kappa,\varpi}\left(\phi^{\epsilon}_s(Q)\right)\right)
\end{array}
$$

Applying the Burkholder-Davis-Gundy inequality presented in~\cite{Bishop/DelMoral/Niclas:2017} we find
$$
\vertiii{ \MM^{\epsilon}_{t}(Q)}_{2n}^2 ~\leq\, c\,n\,\int_0^t\EE\left[\, \Vert \Ea_{t-s}(\phi^{\epsilon}_s(Q))\Vert^{4n}~\tr(\phi^{\epsilon}_s(Q))^n~\tr\left(\Sigma_{\kappa,\varpi}\left(\phi^{\epsilon}_s(Q)\right)\right)^n\,\right]^{1/n}\,ds
$$
The estimates (\ref{trace-Phi-inverse}) imply that
\begin{eqnarray*}
\vertiii{ \MM^{\epsilon}_{t}(Q)}_{2n}^2  &\leq& c\,n\,\int_0^t\Vert \Ea_{t-s}(\mathscr{P}_{\infty})\Vert^{4}\,\EE\left[ (1+\Vert \phi^{\epsilon}_s(Q)\Vert^{5n})\,(\lambda_1(U)^n+\lambda_1(V)^n\,\Vert \phi^{\epsilon}_s(Q)\Vert^{2n})
~\right]^{1/n}\,ds\\~\\
 &\leq& c_n\,(1+\Vert Q\Vert^7)
\end{eqnarray*}
We conclude that
$$
\vertiii{ \MM^{\epsilon}_{t}(Q)}_{2n-1} \,\leq\, \vertiii{ \MM^{\epsilon}_{t}(Q)}_{2n} \,\leq\, c_n\,(1+\Vert Q\Vert^{7/2})
$$
and therefore
$$
\epsilon^{-1}\vertiii{\phi^{\epsilon}_t(Q)-\phi_{t}(Q)}_n \,\leq\, c_n\,(1+\Vert Q\Vert^7)
$$
This ends the proof of (\ref{ref-phi-1-est}). 

Observe when $\kappa=0$, for any $\epsilon\geq 0$ the estimates (\ref{ref-trace-unif-intro}) implies that
$$
\EE\left(\Vert \MM^{\epsilon}_{t}(Q)\Vert^{n}\right)^{1/n}
\displaystyle \,\leq\, c\,n^{1/2}\,(1+ \Vert Q\Vert^{5/2})\,(1+\epsilon~\sqrt{n})^{5/2}
$$
and therefore
\begin{eqnarray*}
\epsilon^{-1}~\vertiii{\phi^{\epsilon}_t(Q)-\phi_{t}(Q)}_n&\leq& c\,(1+\epsilon~\sqrt{n})^{5/2}\left[
\,(1+\epsilon\,\sqrt{n})^{5/2}+\epsilon\,\sqrt{n}\right]\,(1+ \Vert Q\Vert^5)\\
&\leq&c\,(1+ \Vert Q\Vert^5)\,(1+\epsilon~\sqrt{n})^5
\end{eqnarray*}
This ends the proof of (\ref{ref-V-0-1}).
 
The proof of the theorem is now complete.\qed

\subsection{Proof of Theorem~\ref{theo-stab-intro}}\label{theo-stab-intro-proof}

Using (\ref{trace-Phi-inverse-bis}) we have for any $t\geq \upsilon$ and $Q\in \Sa^+_r$ we have the uniform estimate
$$
\Pi_t^{\epsilon}(\Lambda)(Q)~\leq c_{\upsilon}
$$
as soon as $\epsilon\leq \varepsilon_{1}(U,V)\wedge \varepsilon_{1}(V)$, for some constant $c_{\upsilon}$ whose values only depend on $\upsilon$. This implies that $\Lambda$ is a Lyapunov function with compact level sets. Also note, for any bounded measurable function $F$ on $\Sa^0_r$, any $t>0$, and any $P\in \Sa^0_r$ we have
$$
\int_{\Sa^+_r}\Pi_t^{\epsilon}(P,dQ)\,F(Q) \,=\, \int_{\Da_{\overline{r}}}\,\pi_t^{\epsilon}(\varsigma(P),dq)\,(F\circ\varsigma^{-1})(q)
$$
Recalling that continuous images of compact sets are compact, and the density $(p,q)\mapsto \rho_t^{\epsilon}(p,q)$ is continuous for any $t>0$,
for any compact set $K\subset\Sa_r^+$ we have
$$
\inf_{(p,q)\in \varsigma(K)^2} \,\rho_t^{\epsilon}(p,q) \,:=\, \rho_{t,K}^{\epsilon} \,>\, 0
$$
We conclude that for any compact $K\subset\Sa_r^+$, $P\in K$ and $F\geq 0$ we have
$$
\int_{\Sa^+_r}\Pi_t^{\epsilon}(P,dQ)\,F(Q) \,=\, \int_{\Da_{\overline{r}}}\,\rho_t^{\epsilon}(\varsigma(P),q)\,
(F\circ\varsigma^{-1})(q)\,\gamma_{\overline{r}}(dq) \,\,\geq\,\, \varrho_{t,K}^{\epsilon}\int_{\Da_{\overline{r}}}\, \overline{\gamma}_{\varsigma(K)}\,(F\circ\varsigma^{-1})(q)
$$
with the uniform probability measure $\overline{\gamma}_{\varsigma(K)}$ on $\varsigma(K)$ defined by
$$
 \overline{\gamma}_{\varsigma(K)}(dq) \,:=\, \frac{\gamma_{\overline{r}}(dq)\,1_{\varsigma(K)}(q)}{\gamma_{\overline{r}}(\varsigma(K))}
 \quad \mbox{\rm and the parameter}\quad
\varrho_{t,K}^{\epsilon} \,:=\, \rho_{t,K}^{\epsilon}\,\gamma_{\overline{r}}(\varsigma(K)) \,>\, 0
$$
Then, for any compact $K$ and any $t>0$ we have the minorisation condition
 $$
 \forall P\in K, \qquad \Pi^{\epsilon}_t(P,dQ) \,\geq\, \varrho_{t,K}^{\epsilon}\,\overline{\Gamma}_{K}(dQ)
 $$
with the uniform probability measure $\overline{\Gamma}_{K}$ on $K$. This condition, combined with the existence of a Lyapunov function with compact level sets, ensures that the law $Q_t$ converges exponentially fast to a unique invariant measure $\Gamma_{\infty}^{\epsilon}=\Gamma_{\infty}^{\epsilon}\Pi_t^{\epsilon}$, as the time horizon $t\rightarrow\infty$. The contraction estimates are now a consequence of Theorem 8.2.21 and Theorem 17.4.1 in~\cite{penev}; see also~\cite{hairer-2}. This completes the proof of the theorem. \qed

\subsection{Proof of Theorem~\ref{det-E-theo}}\label{det-E-theo-proof}

Using (\ref{trace-detQt}), for any $\zeta\in\RR$, we have
\begin{eqnarray}
\mbox{\rm det}(\Ea_t^{\epsilon}(Q)\Ea_t^{\epsilon}(Q)^{\prime})^{\zeta}
&=&\exp{\left[2\zeta\int_0^t\tr(A-Q_s\,S)\,ds\right]}\nonumber\\
&\leq &\mbox{\rm det}(Q_tQ^{-1})^{\zeta}\exp{\left[
-\zeta\int_0^t
\tr\left(Q^{-1}_sR^{\epsilon}+Q_s\,S^{\epsilon}\right)\,ds+\epsilon\,\zeta\,\mathfrak{m}_t\right]} \label{ref-girsanov}
\end{eqnarray}
with the martingale
$$
d\mathfrak{m}_t \,:=\, -\tr\left(Q^{-1}_tdM_t\right) \quad\Longrightarrow\quad \partial_t\langle \mathfrak{m}~\vert~\mathfrak{m}\rangle_t \,=\, \tr\left(Q^{-1}_t\Sigma_{\kappa,\varpi}(Q_t)\right) \,\leq\, \tr\left(Q^{-1}_tU+Q_tV\right)
$$ 
This implies that
$$
\begin{array}{l}
\EE\left[\mbox{\rm det}(\Ea_t^{\epsilon}(Q)\Ea_t^{\epsilon}(Q)^{\prime})^{\zeta}\right]\\
\\
\qquad\qquad\displaystyle\leq~ \EE\left[\mbox{\rm det}(Q_tQ^{-1})^{2\zeta}\right]^{1/2}~\EE\left[\exp{\left[
-2\zeta\int_0^t
\tr\left(Q^{-1}_sR^\epsilon_{\zeta}+Q_s\,S^\epsilon_{\zeta}\right)
\,ds\right]}\,\Za_{t,\zeta}^{\epsilon}\right]^{1/2} 
\end{array}
$$
with the parameters  $(R^\epsilon_{\zeta},S^\epsilon_{\zeta})$ introduced in (\ref{def-RS-repsilon}) and the exponential martingale
$$
\Za_{t,\zeta}^{\epsilon} \,:=\,\exp{\left[2\epsilon\,\zeta\,\mathfrak{m}_t-\frac{(2\epsilon\,\zeta)^2}{2}\,\langle \mathfrak{m}~\vert~\mathfrak{m}\rangle_t\right]}
$$
By Friedland's inequality (\ref{krause-ref}) for any $4\zeta\geq 1$ we have
$$
\vertiii{\mbox{\rm det}(Q_t)}_{2\zeta}
\,\leq\, 
\mbox{\rm det}(\phi_t(Q))+
c\,\left(\Vert \phi_t(Q)\Vert^{r-1}+\vertiii{Q_t}_{4\zeta(r-1)}^{r-1}
\right)
\vertiii{ Q_t-\phi_t(Q)}_{4\zeta}
$$
Recalling that $\mbox{\rm det}(Q)\leq r^{-r}~\tr(Q)^r\leq c~\Vert Q\Vert^r$ and using (\ref{trace-Phi-inverse}) we then check that
$$
\epsilon\,\leq\, \varepsilon_{4\zeta r}(V) \qquad\Longrightarrow\qquad
\vertiii{\mbox{\rm det}(\phi^{\epsilon}_t(Q))}_{2\zeta}
\,\leq\,
c_{\zeta}\,(1+\Vert Q\Vert^{r+1})
$$
In this case, using (\ref{QQ-1}) we conclude that
 $$
\displaystyle\mbox{\rm det}(Q)\,\EE\left[\mbox{\rm det}(\Ea_t^{\epsilon}(Q)\Ea_t^{\epsilon}(Q)^{\prime})^{\zeta}\right]^{1/\zeta}
\displaystyle \,\leq\,  c_{\zeta}\,(1+\Vert Q\Vert^{r+1})
\exp{\left[
-2t\,\sqrt{\tr\left(R^\epsilon_{\zeta}S^\epsilon_{\zeta}\right)}
\right]}
$$
This ends the proof of the estimate (\ref{ars-0}). 

Now we come to the proof of (\ref{ars}). Using (\ref{ref-girsanov}) for any $(\zeta,\zeta^{\prime})\in\RR^2$ such that 
$$
2\zeta^{\prime}>\zeta>0 \quad\Longleftrightarrow\quad -1< \,\xi:=\frac{\zeta}{\zeta^{\prime}}-1 \,<1
$$ we have
\begin{eqnarray*}
\exp{\left[2\zeta^{\prime}~\int_0^t~\tr(A-Q_sS)~ds\right]}
&\leq &\Za_{t,\zeta/2}^{\epsilon}~\mbox{\rm det}(Q_tQ^{-1})^{\zeta}
\displaystyle\exp{\left[-\int_0^t~F_{\zeta,\zeta^{\prime}}^{\epsilon}(Q_s)
~ds\right]}~
\end{eqnarray*}
with the functional
\begin{eqnarray*}
F_{\zeta,\zeta^{\prime}}^{\epsilon}(Q)&=&2(\zeta-\zeta^{\prime})\,\tr(A-QS)+
\zeta\,
\tr\left(Q^{-1}R^{\epsilon}+QS^{\epsilon}\right)-\frac{(\epsilon\,\zeta)^2}{2}\,\tr\left(Q^{-1}U+QV\right)\\
&=&2(\zeta-\zeta^{\prime})\,\tr(A)+(2\zeta^{\prime}-\zeta)\,
\tr\left(Q\left[S-\frac{\zeta }{2\zeta^{\prime}-\zeta}
\frac{\epsilon^{2}}{2}
\left[\frac{r+1}{2}+\zeta\right]\,V\right]\right)\\
&&\qquad\qquad\qquad\qquad\qquad\qquad\qquad\qquad\qquad+\zeta\,\tr\left(Q^{-1}\left[R-
\frac{\epsilon^{2}}{2}
\left[\frac{r+1}{2}+\zeta\right]\,U\right]\right)
\end{eqnarray*}
with the matrices $(U,V)$ defined in (\ref{cond-UV-ref-intro-uv}). Rewritten in terms of the parameters $(\xi,\zeta^{\prime})$ we have,
\begin{eqnarray*}
\frac{1}{\zeta^{\prime}}\,F_{\zeta^{\prime}(1+\xi),\zeta^{\prime}}^{\epsilon}(Q)&=&2\xi\,\tr(A)+(1-\xi)\,
\tr\left(Q\left[S-\frac{1+\xi }{1-\xi}
\frac{\epsilon^{2}}{2}
\left[\frac{r+1}{2}+\zeta^{\prime}(1+\xi)\right]\,V\right]\right)\\
&&\qquad\qquad\qquad\qquad\qquad+(1+\xi)\,\tr\left(Q^{-1}\left[R-
\frac{\epsilon^{2}}{2}
\left[\frac{r+1}{2}+\zeta^{\prime}(1+\xi)\right]\,U\right]\right)
\end{eqnarray*}
for any $\vert\xi\vert\leq 1$ and $\zeta^{\prime}>0$.
We let 
$$
	\xi_0:=\frac{\tr(A)}{\sqrt{\tr(A)^2+\tr(RS)}}
$$
and we choose $\epsilon$ such that
$$
\widetilde{S}_{\zeta^{\prime}}^{\epsilon}:=\,S-\frac{\epsilon^{2}}{2}\,\frac{1+\xi_0 }{1-\xi_0}
\left[\frac{r+1}{2}+\zeta^{\prime}(1+\xi_0)\right]V\geq 0\quad \mbox{\rm and}\quad
\widetilde{R}_{\zeta^{\prime}}^{\epsilon} :=\, R-
\frac{\epsilon^{2}}{2} \left[\frac{r+1}{2}+\zeta^{\prime}(1+\xi_0)\right]U\geq 0
$$
Now, using (\ref{QQ-1}) we check that
\begin{eqnarray*}
\frac{1}{2\zeta^{\prime}}\,F_{(1+\xi_0)\zeta^{\prime},\zeta^{\prime}}^{\epsilon}(Q)&\geq &\xi_0\,\tr(A)+\sqrt{1-\xi_0^2}
\,\sqrt{\tr\left(\widetilde{R}_{\zeta^{\prime}}^{\epsilon}\widetilde{S}_{\zeta^{\prime}}^{\epsilon}\right)}\\
&=&\frac{\tr(A)^2}{\sqrt{\tr(A)^2+\tr(RS)}}+\frac{\sqrt{\tr(RS)}}{\sqrt{\tr(A)^2+\tr(RS)}}
\,\sqrt{\tr\left(\widetilde{R}_{\zeta^{\prime}}^{\epsilon}\widetilde{S}_{\zeta^{\prime}}^{\epsilon}\right)}
\end{eqnarray*}
which yields the uniform estimate
$$
\begin{array}{l}
\displaystyle
\frac{1}{2\zeta^{\prime}}\,F_{(1+\xi_0)\zeta^{\prime},\zeta^{\prime}}^{\epsilon}\\
\\
\qquad\qquad\geq~ \displaystyle \sqrt{\tr(A)^2+\tr(RS)}\,\left[1-\frac{\sqrt{\tr(RS)}}{\tr(A)^2+\tr(RS)}\left(\sqrt{\tr(RS)}-
\,\sqrt{\tr\left(\widetilde{R}_{\zeta^{\prime}}^{\epsilon}\widetilde{S}_{\zeta^{\prime}}^{\epsilon}\right)}\right)\right]
\end{array}
$$
We conclude that for any $\zeta^{\prime}\geq 0$ there exists some $\varepsilon_0$ and some function $\hbar_{\zeta^{\prime}}(\epsilon)\in [0,1]$ such that $\lim_{\epsilon\rightarrow 0}\hbar_{\zeta^{\prime}}(\epsilon)=0$ such that for any time horizon $t\geq 0$ and any $Q>0$ we have the almost sure estimate
$$
\begin{array}{l}
\displaystyle
\exp{\left[2\zeta^{\prime} \int_0^t\tr(A-Q_sS)\,ds\right]}\\
\\
\qquad\displaystyle\leq~\Za_{t,(1+\xi_0)\zeta^{\prime}/2}^{\epsilon}\,\mbox{\rm det}(Q_tQ^{-1})^{(1+\xi_0)\zeta^{\prime}}
\displaystyle\exp{\left[-2\zeta^{\prime}\left(\sqrt{\tr(A)^2+\tr(RS)}\,(1-\hbar_{\zeta^{\prime}}(\epsilon))\right)t
\right]}~
\end{array}
$$

Moreover, for any non-negative parameters $(\epsilon,\zeta)$ the exponential martingale $\Za_{t,\zeta}^{\epsilon}$
can be interpreted as a change of probability measure. Let $\Fa_t$ be the filtration generated by the diffusion $Q_t$ and let $\PP^{\epsilon}_{\zeta}$ be the probability defined by 
$$
\Za_{t,\zeta}^{\epsilon}\,:=~\exp{\left[\int_0^t \tr\left(H^{\epsilon}_{s,\zeta}d\Wa_s\right)-\frac{1}{2}\int_0^t\tr\left(H^{\epsilon}_{s,\zeta}\left(
H^{\epsilon}_{s,\zeta}\right)^{\prime}\right)ds\right]} \,=\,\frac{d\PP^{\epsilon}_{\zeta}}{d\PP}{|~ \Fa_t}
$$
with the stochastic process
 $$
 \begin{array}{l}
\displaystyle H_{t,\zeta}^{\epsilon}\,=\,-2\epsilon\,\zeta\,\Sigma_{\kappa,\varpi}(Q_t)^{1/2}\,Q^{-1/2}_t\\
 \\
\qquad\quad\displaystyle \Longrightarrow \quad\tr\left(H^{\epsilon}_{s,\zeta}\,d\Wa_s\right) \,=\,-2\epsilon\,\zeta\,
 \tr\left(Q^{-1}_t\left[Q^{1/2}_t\,d\Wa_s\,\Sigma_{\kappa,\varpi}(Q_t)^{1/2}\right]\right) \,=\,-2\epsilon\,\zeta\,\tr\left(Q^{-1}_sdM_s\right)
 \end{array}$$
 By Girsanov's theorem, under $\PP^{\epsilon}_{\zeta}$ the process
 $$
d\Wa^{\epsilon}_{t,\zeta} \,=\, d\Wa_t+2\epsilon\,\zeta\,Q^{-1/2}_t\Sigma_{\kappa,\varpi}(Q_t)^{1/2}\,dt
 $$
 is an $(r\times r)$-Brownian motion. Thus, under $\PP^{\epsilon}_{\zeta}$, the matrix Riccati diffusion $Q_t$ is the solution of the equation
 $$
 dQ_t \,=\, \Theta_{\epsilon,\zeta}(Q_t)dt+\epsilon\,dM_t
 $$
 with the drift function
 $$
 \Theta_{\epsilon,\zeta}(Q) \,=\,AQ+QA^{\prime}+\left(R-2\epsilon^2\zeta U\right)-Q\left[S+2\epsilon^2\zeta V\right]Q \,\leq\, \Theta(Q)
 $$
We conclude that
$$
\begin{array}{l}
\displaystyle
\EE\left[\exp{\left[2\zeta^{\prime}\int_0^t \tr(A-Q_sS)\,ds\right]}\right]\\
\\
\qquad\qquad\displaystyle\leq~ \EE\left(\mbox{\rm det}(Q_{t,\zeta^{\prime}}^\epsilon Q^{-1})^{(1+\xi_0)\zeta^{\prime}}\right)
\exp{\left[-2\zeta^{\prime}\left(\sqrt{\tr(A)^2+\tr(RS)}\,(1-\hbar_{\zeta^{\prime}}(\epsilon))\right)t
\right]}
\end{array}
$$
where $Q_{t,\zeta^{\prime}}^\epsilon$ is a matrix Riccati diffusion defined similarly to $Q_t$ but with the replacement 
$$
(R,S) \quad\longleftarrow\quad \left[(R,S) - (U,-V)\epsilon^2(1+\xi_0)\zeta^{\prime} \right]
$$ 
and with $(U,V)$ defined in (\ref{cond-UV-ref-intro-uv}). This ends the proof of (\ref{ars}).

The proof of the theorem is complete. \qed

\appendix
\section{Appendix}

In this appendix we first derive (\ref{eq-ap-1}) and (\ref{eq-ap-2}). Then we prove the estimate in (\ref{ref-app-1}). Finally, we prove the Liouville formula stated in Lemma \ref{lem-det}.

\subsection{Proof of (\ref{eq-ap-1}) and (\ref{eq-ap-2})}\label{proof-eq-ap-1}
We have
$$
\left\{P_1\otimes_{ s} P_2\right\}=\varsigma \circ(P_1\otimes_{ s} P_2)\circ \varsigma^{-1} \qquad
\Longleftrightarrow\qquad \left\{P_1\otimes_{ s} P_2\right\}\varsigma (H)=\varsigma \left((P_1\otimes_{ s} P_2)(H)\right)
$$
Observe that
$$
 \langle H_1,(P_1\otimes_{ s} P_2)(H_2)\rangle_{\mathrm{Frob}}=\langle
\varsigma (H_1),\left\{P_1\otimes_{ s} P_2\right\}\varsigma (H_2)\rangle_{\overline{r}}
$$
We also have
$$
\left\{P_1\otimes_{ s} P_2\right\}^{1/2}\varsigma (H)=\varsigma \left((P_1\otimes_{ s} P_2)^{1/2}(H)\right)
\quad\Longleftrightarrow\quad
\left\{P_1\otimes_{ s} P_2\right\}^{1/2}=\varsigma \circ(P_1\otimes_{ s} P_2)^{1/2}\circ \varsigma^{-1}
$$
To check this claim notice that
$$
 \begin{array}{l}
\Ta\varsigma (H)\,:=\,
\varsigma \left((P_1\otimes_{ s} P_2)^{1/2}(H)\right)\\
\\
\qquad\Longrightarrow\qquad
\Ta (\Ta\varsigma (H)) \,=\, T\varsigma \left((P_1\otimes_{ s} P_2)^{1/2}(H)\right)
\,=\, \varsigma \left((P_1\otimes_{ s} P_2)^{1/2}(P_1\otimes_{ s} P_2)^{1/2}(H)\right)\\
\\
\qquad\Longrightarrow\qquad \Ta^2 \,=\, \left\{P_1\otimes_{ s} P_2\right\}
\end{array}
$$
This ends the proof of (\ref{eq-ap-1}).

When $P_1,P_2>0$ we have
$$
\begin{array}{l}
 \varsigma (H)^{\prime}\left\{P_1\otimes_{ s} P_2\right\}\varsigma (H)\\
 \\
 \qquad=~\tr\left(H(P_1\otimes_{ s} P_2)H \right)\\
\\\qquad=~
\tr\left(HP_1HP_2 \right) \,=\, \tr\left([H^{1/2}P_1H^{1/2}] [H^{1/2}P_2H^{1/2}] \right)\,>\,0, \qquad \forall H\in \Sa_r-\{0\}
\end{array}
$$
Then, we also have
$$
\begin{array}{l}
 \varsigma (H)^{\prime}\left\{P_1\otimes_{ s} P_2\right\}\varsigma (H) \,\leq\, 
\tr( H^{1/2}P_1H^{1/2})\, \tr(H^{1/2}P_2H^{1/2}) \,\leq\, \lambda_1(P_1)\lambda_2(P_2)\,\Vert H\Vert^2
\\
\\
\qquad\Longrightarrow\qquad \lambda_1(\left\{P_1\otimes_{ s} P_2\right\}) \,\leq\, \lambda_1(P_1)\lambda_1(P_2)
\end{array}
$$
Similarly, we have
$$
\begin{array}{l}
 \varsigma (H)^{\prime}\left\{P_1\otimes_{ s} P_2\right\}\varsigma (H) \,\geq\, \lambda_r(P_2)\,\tr\left(HP_1H \right)
 \,\geq\,  \lambda_r(P_1)\lambda_r(P_2)\,\Vert H\Vert^2\\
 \\
 \qquad\Longrightarrow\qquad \lambda_r(\left\{P_1\otimes_{ s} P_2\right\}) \,\geq\, \lambda_r(P_1)\lambda_r(P_2)
\end{array}
$$
This ends the proof of (\ref{eq-ap-2}). \qed

\subsection{Proof of the Estimate (\ref{ref-app-1})}\label{ref-app-1-proof}

By Lemma 4.7 in~\cite{Bishop/DelMoral:2016} we have the uniform estimate
$$
\Vert \phi_{t}(Q)^{-1}\Vert \,\leq\, c\,(1+\Vert Q^{-1}\Vert)
$$
Using (\ref{trace-Phi-inverse}) and (\ref{ref-phi-1-est}) for any $\epsilon\leq \varepsilon_{2n}(U,V)\wedge \varepsilon_{20n}(V)$ we check that
$$
\begin{array}{l}
\phi^{-\epsilon}_t(Q)- \phi_{t}\left(Q\right)^{-1}\,=\,\phi^{-\epsilon}_t(Q)\left[\phi^{\epsilon}_t(Q)- \phi_{t}\left(Q\right)\right]\phi_{t}\left(Q\right)^{-1}\\
\\
\qquad\Longrightarrow\quad \vertiii{\phi^{-\epsilon}_t(Q)- \phi_{t}\left(Q\right)^{-1}}_n \,\leq\, c\,(1+\Vert Q^{-1}\Vert)\,\vertiii{\phi^{\epsilon}_t(Q)-\phi_{t}(Q)}_{2n}\vertiii{\phi^{-\epsilon}_t(Q)}_{2n}\\
\\
\qquad\Longrightarrow\quad \vertiii{\phi^{-\epsilon}_t(Q)- \phi_{t}\left(Q\right)^{-1}}_n \,\leq\, c_n\,\epsilon\,(1+\Vert Q^{-1}\Vert)\,(1+ \Vert Q\Vert^8)\end{array}
$$
This ends the proof of (\ref{ref-app-1}). \qed

\subsection{Proof of Lemma~\ref{lem-det}}\label{lem-det-proof}

Fix some matrix $Q\in\Sa_r^+$ and set
$$
Q_t\,=\,\phi^{\epsilon}_t(Q)\quad\mbox{\rm and}\quad
\widetilde{Q}_t \,:=\,
\mathcal{E}_{t}^{\epsilon}(Q)^{-1}Q_t\left(\mathcal{E}_{t}^{\epsilon}(Q)^{\prime}\right)^{-1} \qquad\Longleftrightarrow\qquad
\widetilde{Q}_t^{-1} \,:=\, \mathcal{E}_{t}^{\epsilon}(Q)^{\prime}Q^{-1}_t\mathcal{E}_{t}^{\epsilon}(Q)
$$
Note that
$$
\begin{array}{l}
\displaystyle
d\mathcal{E}_{t}^{\epsilon}(Q)^{-1}\,=\,
-\mathcal{E}_{t}^{\epsilon}(Q)^{-1}\left(d\mathcal{E}_{t}^{\epsilon}(Q)\right)\mathcal{E}_{t}^{\epsilon}(Q)^{-1}\\
\\
\qquad\displaystyle\Longleftrightarrow \qquad \partial_t\mathcal{E}_{t}^{\epsilon}(Q)^{-1}\,=\,
-\mathcal{E}_{t}^{\epsilon}(Q)^{-1}\left(\partial_t\mathcal{E}_{t}^{\epsilon}(Q)\right)\mathcal{E}_{t}^{\epsilon}(Q)^{-1} \,=\, -\mathcal{E}_{t}^{\epsilon}(Q)^{-1}(A-Q_tS)
\end{array}$$
This implies that
\begin{eqnarray*}
d\widetilde{Q}_t&=&\mathcal{E}_{t}^{\epsilon}(Q)^{-1}\left[dQ_t-(A-Q_tS)\,Q_t-Q_t\,
(A-Q_tS)^{\prime}\right]\left(\mathcal{E}_{t}^{\epsilon}(Q)^{\prime}\right)^{-1}\\
&=&\mathcal{E}_{t}^{\epsilon}(Q)^{-1}\left[R+Q_t\,S\,Q_t\right]\left(\mathcal{E}_{t}^{\epsilon}(Q)^{\prime}\right)^{-1}dt+\epsilon~d\widetilde{M}_t\\
&=&\mathcal{E}_{t}^{\epsilon}(Q)^{-1}Q^{1/2}_t\left[Q^{-1/2}_tR\,Q^{-1/2}_t+Q^{1/2}_tS\,Q^{1/2}_t\right]Q^{1/2}_t\left(\mathcal{E}_{t}^{\epsilon}(Q)^{\prime}\right)^{-1}dt+\epsilon\,d\widetilde{M}_t\\
&&\\
\Longrightarrow && ~~~~~~~\,\widetilde{Q}_t^{-1}\,d\widetilde{Q}_t ~=~
\mathcal{E}_{t}^{\epsilon}(Q)^{\prime}\left[Q^{-1}_tR+S\,Q_t\right]\left(\mathcal{E}_{t}^{\epsilon}(Q)^{\prime}\right)^{-1}dt+\epsilon\,  \mathcal{E}_{t}^{\epsilon}(Q)^{\prime}Q^{-1}_t\,dM_t\left(\mathcal{E}_{t}^{\epsilon}(Q)^{\prime}\right)^{-1}\\
\Longrightarrow&& \tr\left(\widetilde{Q}_t^{-1}\,d\widetilde{Q}_t\right)~=~\tr\left(Q^{-1}_tR+S\,Q_t\right)\,dt+\epsilon\,\tr\left(Q^{-1}_t\,dM_t\right)
\end{eqnarray*}
with the martingale
$$
\begin{array}{l}
\displaystyle
d\widetilde{M}_t \,:=\, \mathcal{E}_{t}^{\epsilon}(Q)^{-1}dM_t\left(\mathcal{E}_{t}^{\epsilon}(Q)^{\prime}\right)^{-1}\\
\\
\Longrightarrow \quad
d\widetilde{M}_t~\widetilde{Q}_t^{-1} =\mathcal{E}_{t}^{\epsilon}(Q)^{-1}~dM_t~Q^{-1}\mathcal{E}^{\epsilon}_{t}(Q)\\
\\
\hskip3cm\quad\mbox{\rm and}\quad
d\widetilde{M}_t~\widetilde{Q}_t^{-1}d\widetilde{M}_t~\widetilde{Q}_t^{-1}=\mathcal{E}_{t}^{\epsilon}(Q)^{-1}~dM_t~Q^{-1}~dM_t~Q^{-1}\mathcal{E}^{\epsilon}_{t}(Q)\\
\\
\Longrightarrow \quad \tr\left(d\widetilde{M}_t\,\widetilde{Q}_t^{-1} \right) \,=\, \tr(dM_t\,Q^{-1})\quad\mbox{\rm and}\quad
\tr\left(d\widetilde{M}_t\,\widetilde{Q}_t^{-1}d\widetilde{M}_t\,\widetilde{Q}_t^{-1}\right) \,=\,
\tr\left(dM_t\,Q^{-1}\,dM_t\,Q^{-1}\right)
\end{array}
$$
For a more rigorous derivation of the angle bracket of matrix-valued martingales we refer the reader to Section 3 in~\cite{Bishop/DelMoral/Niclas:2017}.

The determinant function $f(\cdot):=\mbox{\rm det}(\cdot)$ is smooth on the space of invertible matrices. The first and second Fr\'echet derivatives are given for any $H,H_1,H_2\in\Ma_r$ by the Jacobi formulae
\begin{eqnarray*}
\nabla f(A)\cdot H&=&f(A)\,\tr(HA^{-1})\\
\nabla^2 f(A)\cdot (H_1,H_2)&=& -f(A)\,\left[\tr(H_1A^{-1}H_2A^{-1})-\tr(H_1A^{-1})\tr(H_2A^{-1})\right]
\end{eqnarray*}
Using the Ito differential calculus for stochastic matrix diffusions developed in~\cite{Bishop/DelMoral/Niclas:2017}, with a slight abuse of notation we find the formula,
\begin{eqnarray*}
df(\widetilde{Q}_t)&=&f(\widetilde{Q}_t)\,\left[
\tr(\widetilde{Q}_t^{-1} d\widetilde{Q}_t)-\frac{\epsilon^2}{2}\left[\tr(d\widetilde{M}_t\widetilde{Q}_t^{-1}d\widetilde{M}_t\widetilde{Q}_t^{-1})-\tr(d\widetilde{M}_t\widetilde{Q}_t^{-1})\tr(d\widetilde{M}_t\widetilde{Q}_t^{-1})\right]
\right]\\
&=&f(\widetilde{Q}_t)\,\bigg[
\tr(Q^{-1}_tR+SQ_t) \\
&&\qquad -\frac{\epsilon^2}{2}\left[\tr\left(dM_t\,Q_t^{-1}\,dM_t\,Q_t^{-1}\right)-\tr(dM_t\,Q_t^{-1})\tr(dM_t\,Q_t^{-1})\right]\bigg]\,dt+\epsilon\,dM_t(f)
\end{eqnarray*}
with the martingale
$$
dM_t(f) \,:=\, f(\widetilde{Q}_t)\,\tr\left(Q^{-1}_t\,dM_t\right)
$$
Recalling that
$$
\begin{array}{l}
\displaystyle
2\,dM_t\,Q_t^{-1} \,=\, Q_t^{1/2}\,d\Wa_t\,\Sigma^{1/2}_{\kappa,\varpi}\left(Q_t\right)Q_t^{-1}+\Sigma^{1/2}_{\kappa,\varpi}\left(Q_t\right)\,d\Wa^{\prime}_t\,Q_t^{-1/2}
\end{array}
$$
we check that
$$
\begin{array}{l}
\displaystyle
\tr\left(dM_t\,Q_t^{-1}\right)\,=\,\tr\left(d\Wa_t\,\Sigma^{1/2}_{\kappa,\varpi}\left(Q_t\right)Q_t^{-1/2}\right)\\
\\
\qquad\Longrightarrow\qquad\begin{array}[t]{rcl}
\tr\left(dM_t\,Q_t^{-1}\right)\tr\left(dM_t\,Q_t^{-1}\right)&=&\tr\left(\Sigma^{1/2}_{\kappa,\varpi}\left(Q_t\right)Q_t^{-1}\,\Sigma^{1/2}_{\kappa,\varpi}\left(Q_t\right)\right)dt\\
&\leq &~\tr(Q_t^{-1}U+Q_tV)\,dt
\end{array}
\\
\\
\qquad\Longrightarrow\qquad dM_t(f)\,dM_t(f)=f(\widetilde{Q}_t)^2\,\tr\left(Q_t^{-1}\,\Sigma_{\kappa,\varpi}\left(Q_t\right)\right) \,\leq\, f(\widetilde{Q}_t)^2\,
\tr\left(Q_t^{-1}U+Q_tV\right)dt
\end{array}
$$
The first implication follows from the fact that
$$
\tr(Q\,d\Wa_t)\,\tr(d\Wa_t\,Q)=\tr(QQ^{\prime})\,dt
$$
Similarly, we have
$$
\begin{array}{l}
\displaystyle
4\,\tr\left(dM_t~Q_t^{-1}\,dM_t\,Q_t^{-1}\right)\\
\\
\qquad\displaystyle=\,\tr\left[\left(Q_t^{1/2}\,d\Wa_t\,\Sigma^{1/2}_{\kappa,\varpi}\left(Q_t\right)Q_t^{-1}+\Sigma^{1/2}_{\kappa,\varpi}\left(Q_t\right)\,d\Wa^{\prime}_t\,Q_t^{-1/2}\right)\right.\\
\\
\hskip3cm\left.\times \left(Q_t^{1/2}\,d\Wa_t\,\Sigma^{1/2}_{\kappa,\varpi}\left(Q_t\right)Q_t^{-1}+\Sigma^{1/2}_{\kappa,\varpi}\left(Q_t\right)\,d\Wa^{\prime}_t\,Q_t^{-1/2}\right)\right]\\
\\
\qquad\displaystyle=\,2\,\tr\left(d\Wa_t\,\Sigma^{1/2}_{\kappa,\varpi}\left(Q_t\right)Q_t^{-1/2}\,d\Wa_t\,\Sigma^{1/2}_{\kappa,\varpi}\left(Q_t\right)\,Q_t^{-1/2}\right)\\
\\
\hskip3cm+\,2\,\tr\left(\Sigma^{1/2}_{\kappa,\varpi}\left(Q_t\right)Q_t^{-1}\Sigma^{1/2}_{\kappa,\varpi}\left(Q_t\right)\,d\Wa^{\prime}_t\,d\Wa_t\right)
\end{array}
$$
Recalling the standard identities,
$$
d\Wa_t\,Q\,d\Wa_t\,=\,Q^{\prime}\,dt\qquad\mbox{\rm and}\qquad 
d\Wa_t\,d\Wa_t^{\prime} \,=\, r\,I\,dt \,=\, d\Wa_t^{\prime}\,d\Wa_t
$$
we check that
$$
\begin{array}{l}
\displaystyle
\tr\left(dM_t\,Q_t^{-1}\,dM_t\,Q_t^{-1}\right)
\displaystyle\,=\,\frac{r+1}{2}\,\tr\left(\Sigma^{1/2}_{\kappa,\varpi}\left(Q_t\right)Q_t^{-1}\,\Sigma^{1/2}_{\kappa,\varpi}\left(Q_t\right)\right) \,\leq\, \frac{r+1}{2}\,\tr\left(
Q_t^{-1}U+VQ_t
\right)
\end{array}
$$
For a more rigorous derivation of the angle bracket of matrix-valued martingales we refer the reader to Section 3 in~\cite{Bishop/DelMoral/Niclas:2017}. In summary, we have proved that
\begin{eqnarray*}
df(\widetilde{Q}_t)&=&f(\widetilde{Q}_t)\,\left[\tr(Q^{-1}_tR+SQ_t)-\epsilon^{2}\,\frac{r-1}{4}\,\tr\left(Q_t^{-1}\Sigma_{\kappa,\varpi}\left(Q_t\right)\right)
\right]\,dt+\epsilon\,f(\widetilde{Q}_t)\,\tr\left(Q^{-1}_t\,dM_t\right)\\
&\geq &f(\widetilde{Q}_t)\,\left[\tr(Q^{-1}_tR+S\,Q_t)-\epsilon^{2}\,\frac{r-1}{4}\,\tr\left(Q_t^{-1}U+VQ_t\right)
\right]\,dt+\epsilon\,f(\widetilde{Q}_t)\,\tr\left(Q^{-1}_t\,dM_t\right)\\
&=&f(\widetilde{Q}_t)\,\left[\tr\left(Q^{-1}_t\left(R-\epsilon^{2}\,\frac{r-1}{4}\,U\right)\right)+\tr\left(Q_t\left(S-\epsilon^{2}\,\frac{r-1}{4}\,V\right)\right)
\right]\,dt+\epsilon\, dM_t(f)
\end{eqnarray*}

Now let $g(\cdot)\, :=\, \log{f(\cdot)}$. Applying Ito's formula we conclude that
\begin{eqnarray*}
\displaystyle dg(\widetilde{Q}_t)&=&\left[\tr(Q^{-1}_tR+SQ_t)-\frac{\epsilon^{2}}{2}\,\frac{r+1}{2}\,\tr\left(Q_t^{-1}\Sigma_{\kappa,\varpi}\left(Q_t\right)\right)
\right]\,dt+\epsilon\,\tr\left(Q^{-1}_t\,dM_t\right)\\
\\
&\geq& \left[\tr\left(Q^{-1}_t\left(R-\frac{\epsilon^{2}}{2}\,\frac{r+1}{2}\,U\right)\right)+\tr\left(Q_t\left(S-
\frac{\epsilon^{2}}{2}\,\frac{r+1}{2}\,V\right)\right)\right]\,dt+\epsilon\,\tr\left(Q^{-1}_t\,dM_t\right)
\end{eqnarray*}
This ends the proof of the lemma.\qed

\end{document}